\documentclass[12pt,letterpaper, final]{amsart}
\usepackage{amssymb, amsmath, amsthm}
\usepackage{mathrsfs}  
\usepackage{graphicx}
\usepackage[colorlinks=true, citecolor=blue, linkcolor=red, urlcolor=blue]{hyperref}
\usepackage[text={6in,8in},centering]{geometry}
 \usepackage{xcolor}
\allowdisplaybreaks % allow break in formulas

\usepackage{mathrsfs}  

\usepackage{bm}
\usepackage{ifdraft}
\ifoptionfinal{
\usepackage[disable]{todonotes}
}{
\usepackage[norefs, nocites]{refcheck}

\usepackage[notref, notcite]{showkeys}
\usepackage[bordercolor=white, color=white]{todonotes}
}

\makeatletter\providecommand\@dotsep{5}\def\listtodoname{List of Todos}\def\listoftodos{\hypersetup{linkcolor=black}\@starttoc{tdo}\listtodoname\hypersetup{linkcolor=blue}}\makeatother
\usepackage{hyperref}
\usepackage{etoolbox}

\newtheorem{theorem}{Theorem}[section]
\newtheorem{lemma}[theorem]{Lemma}
\newtheorem{corollary}[theorem]{Corollary}
\newtheorem{proposition}[theorem]{Proposition}
\theoremstyle{definition}
\newtheorem{definition}[theorem]{Definition}

\theoremstyle{remark}
\newtheorem{remark}[theorem]{Remark}
\newtheorem{example}{Example}
\numberwithin{equation}{section}

\def\N{\mathbb N}

\renewcommand{\leq}{\leqslant}
\renewcommand{\geq}{\geqslant}

\def\p{\partial}

\newcommand*\xbar[1]{%
   \hbox{%
     \vbox{%
       \hrule height 0.5pt % The actual bar
       \kern0.5ex%         % Distance between bar and symbol
       \hbox{%
         \ensuremath{#1}%
       }%
     }%
   }%
} 

\title[]{Inverse Spectral Problems with Sparse Data and Applications to Passive Imaging on manifolds}
\author[Feizmohammadi]{Ali Feizmohammadi}
\address{A. Feizmohammadi, Department of Mathematics,
		University of Toronto, 3359 Mississauga Road, Mississauga, ON L5L1C6, Canada}
\email{ali.feizmohammadi@utoronto.ca}

\author[Krupchyk]{Katya Krupchyk}
\address
        {K. Krupchyk, Department of Mathematics\\
University of California, Irvine\\
CA 92697-3875, USA }

\email{katya.krupchyk@uci.edu}

\begin{document}

%\begin{titlepage}
%\maketitle
%\end{titlepage}

\maketitle

\begin{abstract}

Motivated by inverse problems with a single passive measurement, we introduce and analyze a new class of inverse spectral problems on closed Riemannian manifolds. Specifically, we establish two general uniqueness results for the recovery of a potential in the stationary Schr\"odinger operator from partial spectral data, which consists of a possibly sparse subset of its eigenvalues and the restrictions of the corresponding eigenfunctions to a nonempty open subset of the manifold. Crucially, the eigenfunctions are not assumed to be orthogonal, and no information about global norming constants is required.

The partial data formulation of our inverse spectral problems is naturally suited to the analysis of inverse problems with passive measurements, where only limited observational access to the solution is available. Leveraging this structure, we establish generic uniqueness results for a broad class of evolutionary PDEs, in which both the coefficients and the initial or source data are to be recovered from knowledge of the solution restricted to a subset of spacetime. These results introduce a spectral framework for passive imaging and extend inverse spectral theory into a regime characterized by highly incomplete, yet physically realistic, data.

\end{abstract} 

%\tableofcontents

\section{Introduction}

The interplay between geometry and spectral theory has long been a central theme in mathematical analysis, dating back to the seminal question posed by Kac~\cite{Kac69}: ``Can one hear the shape of a drum?'' This question, formalized within the framework of inverse spectral problems, seeks to determine to what extent the geometric or topological structure of a Riemannian manifold $(M,g)$ is encoded in the spectrum of a natural differential operator, typically the canonical Laplace--Beltrami operator $-\Delta_g$.
We refer to~\cite{milnor1964, GWW} for striking classical counterexamples, to the recent breakthrough~\cite{HZ22}, and to~\cite{Zelditch_2004, DH13, G00} for comprehensive surveys on inverse spectral problems.

In this article, we investigate a variant of the classical inverse spectral problem in which, alongside (a possibly sparse portion of) the spectral data, supplementary information is available in the form of the knowledge of the associated eigenfunctions, restricted to an open subset $O$ of the manifold referred to as the observation region. Within the framework of Kac's question--``Can one hear the shape of a drum?''-- this modified problem may be interpreted as only hearing a (possibly small) subset of the natural frequencies of the drum together with observing a corresponding subset of vibrational modes on a localized region $O$ of the surface. This corresponds to the scenario in which an observer is able to "see" a part of the drum as it vibrates. 
\begin{itemize}
  \item[{\bf (Q)}] {\it How much global information in a physical system can be recovered from partial knowledge of the spectrum and measurements of an associated set of eigenfunctions at certain locations?}
\end{itemize}
Before stating our two main results on the above inverse spectral problem, namely Theorem~\ref{thm_spectral} and Theorem~\ref{thm_anosov}, let us briefly clarify the notion of {\em partial} or {\em sparse} knowledge of the spectral data appearing in \textbf{(Q)}, as well as the title of this work. This aspect constitutes a key novelty of our formulation and is what makes it suitable for applications to the so-called inverse problems with passive measurements.

Our aim is to consider two Schr\"odinger operators $P_1=-\Delta_g+V_1$ and $P_2=-\Delta_g+V_2$ on a closed Riemannian manifold 
$(M,g)$ and to show their equivalence under the following minimal assumptions: (i) the spectrum of $P_1$ is a subset of the spectrum of $P_2$, meaning that, a priori, our spectral knowledge of $P_2$ may be severely limited, as it could potentially have a much larger spectrum; (ii) the restrictions to a fixed set $O\subset M$ of a corresponding set of linearly independent eigenfunctions for $P_1$ and $P_2$ are equal, with no assumed knowledge of norming constants or orthogonality properties; and
(iii) a portion (possibly infinite) of the spectrum of $P_1$ may be also missing from our dataset altogether. As we will see, this notion of sparsity, encoded in conditions (i)--(iii), is fundamental in applications to a broad class of inverse problems with passive measurements, where a single passive observation is used to recover information about an unknown medium or structure.

In the recent work~\cite{FKU24}, motivated by the Calder\'on problem for the fractional Schr\"odinger operator, a variant of question~\textbf{(Q)} was affirmatively resolved in the setting of the Laplacian on a closed Riemannian manifold. Specifically,~\cite[Theorem 1.11]{FKU24} establishes that a smooth, closed, connected Riemannian manifold can be uniquely determined, up to isometry, from the knowledge of the Laplacian's eigenvalues and a set of linearly independent (though not necessarily orthonormal) eigenfunctions restricted to an open subset, provided certain geometric conditions on the subset are satisfied. This result represents the strongest known recovery theorem in the setting where the manifold itself is unknown and must be reconstructed from spectral data, without any a priori information about global norming constants or orthogonality relations among the eigenfunctions.

Motivated by inverse problems for partial differential equations with single passive measurements, the present work addresses question~\textbf{(Q)} in a different setting: determining the potential in the stationary Schr\"odinger operator from sparse spectral data, in the sense described two paragraphs above, and assuming the underlying Riemannian manifold is known.

\vspace*{-1.2\baselineskip}

\subsection{Inverse spectral results}
\label{sec_ip_spectral}

Let $(M, g)$ be a smooth, closed, and connected Riemannian manifold of dimension $n \ge 2$. By closed, we mean that $M$ is compact and without boundary. 

Let $-\Delta_g$ denote the positive Laplace--Beltrami operator on $M$, and let $V \in C^\infty(M)$ be a real-valued potential. Consider the Schr\"odinger operator $-\Delta_g + V$. Equipped with the domain $C^\infty(M)$, this operator is unbounded, symmetric, and essentially self-adjoint on $L^2(M)$. Its unique self-adjoint extension has domain $\mathcal{D}(-\Delta_g + V) = H^2(M)$, the standard Sobolev space on $M$; see \cite[Section~12]{Grigis_Sjostrand_book}.

The spectrum of $-\Delta_g + V$ is discrete and consists of real eigenvalues, each isolated and of finite multiplicity:
\[
\inf_{x \in M} V(x) \le \mu_1 \le \mu_2 \le \dots \nearrow +\infty.
\]
The eigenspaces corresponding to distinct eigenvalues are orthogonal in $L^2(M)$, and their direct sum is dense in $L^2(M)$.

We now proceed to state our inverse spectral results for the Schr\"odinger operator $-\Delta_g + V$ with incomplete spectral data. To that end, we begin by introducing the  Rauch--Taylor / Bardos--Lebeau--Rauch geometric control condition; see \cite{Rauch_Taylor_1974, BLR0, BLR}.

\begin{definition}
Let $O \subset M$ be a nonempty open set, and let $T > 0$. We say that the pair $(O, T)$ satisfies the \emph{geometric control condition (GCC)} if every unit-speed, inextendible (maximal) geodesic in $(M, g)$ intersects $O$ at some time $t \in (0, T)$. We also say that $O$ satisfies GCC if there exists $T > 0$ such that $(O, T)$ satisfies the geometric control condition.
\end{definition}

We also need the following definition.
\begin{definition}
	\label{def_antipodal}
	Given any point $p \in M$, the \emph{antipodal set} of $p$, denoted by $\mathcal{A}_{M,g}(p)$, is defined as
	\[
	\mathcal{A}_{M,g}(p) = \left\{ q \in M \,:\, \mathrm{dist}_g(p, q) = \max_{p' \in M} \mathrm{dist}_g(p, p') \right\}.
	\]
\end{definition}
Since the manifold $M$ is compact, we have $\mathcal{A}_{M,g}(p) \neq \emptyset$ for all $p \in M$.  

Let $O \subset M$ be a nonempty open set. We now impose the following assumption on the observation set $O$:
\begin{enumerate}
  \item[\textbf{(H)}] The set $O$ satisfies GCC, and $\exists$ $p \in O$ such that $\mathcal{A}_{M,g}(p) \subset O$.
\end{enumerate}

Our first spectral result is stated as follows.

\begin{theorem}
	\label{thm_spectral}
	Let $(M,g)$ be a smooth, closed, and connected Riemannian manifold of dimension $n \geq 2$, and let $O \subset M$ be a nonempty, open, and connected subset such that $M \setminus \overline{O} \neq \emptyset$. Suppose that assumption~\textbf{(H)} is satisfied. For $j = 1, 2$, let $V_j \in C^{\infty}(M)$ be real-valued potentials, and denote by
	\[
	\mu_1^{(j)} \leq \mu_2^{(j)} \leq \cdots
	\]
	the eigenvalues of the Schr\"odinger operator $-\Delta_g + V_j$ on $M$, listed in increasing order and counted with multiplicity. Let $\{\psi_k^{(j)}\}_{k=1}^{\infty} \subset C^{\infty}(M)$ be a linearly independent set of corresponding eigenfunctions, satisfying
	\[
	(-\Delta_g + V_j) \psi_k^{(j)} = \mu_k^{(j)} \psi_k^{(j)} \quad \text{on } M.
	\]
	Suppose that
	\begin{equation} 
		\label{psi_eq}
		\mu_k^{(1)} = \mu_{b_k}^{(2)} \quad \text{and} \quad \psi_k^{(1)}(x) = \psi_{b_k}^{(2)}(x) \quad \text{for all } x \in O \text{ and all } k \geq N,
	\end{equation}
	for some $N \in \mathbb{N}$ and some strictly increasing sequence $b_N < b_{N+1} < \cdots$ of positive integers. Then $V_1 = V_2$ on $M$.
\end{theorem}

\begin{remark}
We emphasize that the eigenfunctions $\{\psi_k^{(j)}\}_{k=1}^\infty$ are not assumed to be orthonormal. This is a fundamental aspect of the theorem that, as we will see, makes it suitable for applications in many inverse problems for PDEs subject to passive measurements. 

\end{remark}

\begin{remark}
Theorem~\ref{thm_spectral} does not require full spectral data from either operator. In particular, a significant amount of the low-frequency spectra of both $-\Delta_g + V_1$ and $-\Delta_g + V_2$ may remain entirely unknown. Moreover, even in the high-frequency regime, only a subsequence of the spectrum of $-\Delta_g + V_2$ is used in the identification. By Weyl's law and the matching condition \eqref{psi_eq}, the indexing sequence $b_k$ satisfies $b_k = k + o(k)$ as $k \to \infty$. However, this does not preclude the possibility that $b_k - k \to \infty$. For example, one may take
\[
b_k = k + \left\lfloor \frac{k}{\log k} \right\rfloor, \quad k>1.
\]
In this case, among the first $b_k$ eigenvalues of $-\Delta_g + V_2$, only $k$ are used, and exactly $\left\lfloor \frac{k}{\log k} \right\rfloor$ eigenvalues (counting multiplicity) are skipped. Thus, the map $k \mapsto b_k$ may omit infinitely many indices, and a substantial portion of the high-frequency spectrum of $-\Delta_g + V_2$ may also remain unused.
\end{remark}

\begin{remark}
\label{rem_extension_basis}
In Theorem~\ref{thm_spectral}, the assumption that the full system of eigenfunctions $\{\psi_k^{(j)}\}_{k = 1}^\infty$ of the operator $-\Delta_g + V_j$ is linearly independent for $j = 1, 2$ can be relaxed. It suffices to assume linear independence only for the sets $\{\psi_k^{(1)}\}_{k \ge N}$ and $\{\psi_{b_k}^{(2)}\}_{k \ge N}$ of eigenfunctions appearing in~\eqref{psi_eq}. Indeed, if the sequence $\{\mu_{b_k}^{(2)}\}_{k \ge N}$ exhausts the spectrum of $-\Delta_g + V_2$ (counted with multiplicities), then the set $\{\psi_{b_k}^{(2)}\}_{k \ge N}$ already forms a linearly independent system of eigenfunctions for $-\Delta_g + V_2$. Otherwise, we can complete this system as follows: if the multiplicity of some eigenvalue $\mu_{b_k}^{(2)}$ for $k \ge N$ exceeds one, we supplement the system with additional eigenfunctions from the corresponding eigenspace to obtain a full basis; we also include bases for any remaining eigenspaces not yet represented. The resulting set is a linearly independent system of eigenfunctions of $-\Delta_g + V_2$, spanning the entire spectrum (with multiplicities accounted for). Similarly, the system $\{\psi_k^{(1)}\}_{k \ge N}$ can be extended to a linearly independent system of eigenfunctions corresponding to the full spectrum of $-\Delta_g + V_1$, with multiplicities correctly taken into account.
\end{remark}

\begin{remark}
	\label{O_connected}
In Theorem~\ref{thm_spectral}, the connectedness assumption on the observation set $O$ is only required when $N > 1$. In that case, the proof relies on unique continuation of eigenfunctions from a neighborhood of a point $q\in \mathcal{A}_{M,g}(p)$ in assumption~\textbf{(H)} to the entire set $O$, which necessitates its connectedness.
\end{remark}

We next make several remarks concerning assumption~\textbf{(H)} on the observation set~$O$ in Theorem~\ref{thm_spectral}, which plays a critical role in ensuring the identifiability of the potential. As a first point, we emphasize that recovery from the spectral data~\eqref{psi_eq} is not, in general, possible if the measurements are taken on an arbitrary nonempty open set~$O$, as demonstrated already in the one-dimensional case; see~\cite{HL78}.

\begin{remark}
\label{rem_non-trapping}
In assumption \textbf{(H)} of Theorem~\ref{thm_spectral}, the requirement that the set $O$ satisfies the GCC can be replaced by the condition that the manifold $(M \setminus O, g)$ is nontrapping. Indeed, the nontrapping assumption on $(M \setminus O, g)$ implies that, for all sufficiently large $T > 0$, the pair $(O, T)$ satisfies the GCC; see \cite[Lemma~5.1]{FKU24}.

Moreover, if one further assumes that $(M \setminus O, g)$ is a simple Riemannian manifold, meaning that it is smooth, simply connected, has a smooth and strictly convex boundary, and contains no conjugate points, then assumption \textbf{(H)} is satisfied, as shown in \cite[Corollary~1.5]{FKU24}.

As explained in \cite{FKU24}, this framework provides many concrete examples of manifolds $(M, g)$ to which Theorem~\ref{thm_spectral} applies, where the observation set $O$ may be arbitrarily small compared to the volume of $M \setminus O$. For instance, one may start with any compact, simple Riemannian manifold $(N, g)$ with smooth boundary, and then smoothly extend $N$ to a closed, smooth, connected manifold $M$, extending the metric accordingly. In this construction, the set $O = M \setminus N$ satisfies assumption \textbf{(H)}.
\end{remark}

\begin{remark}
Let $(M, g)$ be a smooth, closed, and connected Riemannian manifold of dimension $n \geq 2$. Then there always exists an observation set $O \subset M$ that satisfies assumption~\textbf{(H)}. Indeed, for any point $x_0 \in M$, the closure $\overline{B_r(x_0)}$ of a sufficiently small geodesic ball centered at $x_0$ is a simple Riemannian manifold with boundary. Therefore, by Remark~\ref{rem_non-trapping}, the complement $M \setminus \overline{B_r(x_0)}=:O$ satisfies assumption~\textbf{(H)}.
\end{remark}

\begin{remark}
Let $(M, g)$ be the standard sphere $\mathbb{S}^{n-1}$ equipped with its canonical round metric. Then any open spherical cap strictly larger than the northern hemisphere satisfies assumption~\textbf{(H)} and may be used as the set $O$.
\end{remark}

To proceed, we introduce the following definition; see \cite{Macia_2021}.
\begin{definition}
\label{de_observability_eigenfunctions}
Let $V \in C^\infty(M)$ be real-valued. We say that the eigenfunctions of the Schr\"odinger operator $-\Delta_g + V$ are \emph{observable from an open set} $O \subset M$ if there exists a constant $C > 0$, depending only on $(M, g)$, $O$, and $V$, such that
\begin{equation}
\label{eq_observability_eigenfunctions_def}
\|\phi\|_{L^2(M)} \le C \|\phi\|_{L^2(O)}
\end{equation}
for every eigenfunction $\phi \in C^\infty(M)$ of $-\Delta_g + V$.
\end{definition}

It is known that \eqref{eq_observability_eigenfunctions_def} follows from observability estimates for the wave equation, as well as from observability estimates for the time-dependent Schr\"odinger equation from an open set $O$ at a fixed time $T$; see, for example, \cite{Macia_2021}, and also Lemma~\ref{lem_nontrapping} and Remark~\ref{rem_nontrapping} below. The geometric control condition (GCC) on the pair $(O, T)$ is sufficient for such observability estimates to hold; see \cite{Rauch_Taylor_1974, BLR0, BLR, Lebeau_1992}, and also \cite{Macia_2021}. Furthermore, the GCC is necessary for the observability of the time-dependent Schr\"odinger equation with $V = 0$ on compact Riemannian manifolds with periodic geodesic flow, such as the standard sphere; see~\cite[Theorem~2]{Macia_2021_new}.

\begin{remark}
In assumption \textbf{(H)} of Theorem~\ref{thm_spectral}, the requirement that the set $O$ satisfies the GCC can be replaced by the assumption that the eigenfunctions of the Schr\"odinger operator $-\Delta_g + V$ are observable from $O$. Indeed, the observability estimate \eqref{eq_observability_eigenfunctions_def} in Definition~\ref{de_observability_eigenfunctions} is precisely what is used in the proof of Theorem~\ref{thm_spectral}.
\end{remark}

We next present examples of manifolds where the eigenfunctions of the Schr\"odinger operator $-\Delta_g + V$ are observable from any nonempty open set $O \subset M$, in the sense of Definition~\ref{de_observability_eigenfunctions}. Before doing so, we note that, as explained in \cite{Dyatlov_Jin_2018}, the observability estimate \eqref{eq_observability_eigenfunctions_def} with $V = 0$ fails on the round sphere for certain nonempty open subsets.

\begin{example}
\label{eq_ex_torus}
Let $(M, g)$ be the torus $\mathbb{T}^n = (\mathbb{R} / 2\pi \mathbb{Z})^n$ equipped with the flat metric, and let $V \in C^\infty(\mathbb{T}^n)$ be real-valued. Then the eigenfunctions of $-\Delta_g + V$ are observable from any nonempty open set $O \subset M$, in the sense of Definition~\ref{de_observability_eigenfunctions}. This result was established in \cite{Burq_Zworski_2012},  see also \cite{Bourgain_Burq_Zworski_2013}, in the case $n = 2$ for both rational and irrational tori, and \cite{Anantharaman_Macia_2014} for general $n$. For the case $V = 0$, we refer to \cite{Jaffard_1990} and \cite{Haraux_1989} (for $n = 2$), and \cite{Komornik_1992} (for higher dimensions).
\end{example}

\begin{example}
\label{eq_ex_Anosov}
Let $(M, g)$ be a smooth, closed, and connected Riemannian surface whose geodesic flow has the Anosov property, and let $V \in C^\infty(M)$ be real-valued. Then the eigenfunctions of $-\Delta_g + V$ are observable from any nonempty open set $O \subset M$, in the sense of Definition~\ref{de_observability_eigenfunctions}. This follows from a more general result established in \cite[Theorem~2]{DJN}; see also \cite{Dyatlov_Jin_2018} for the case of compact connected hyperbolic surfaces. It is conjectured in~\cite{ADM_2024} that if $(M, g)$ is a closed, connected Riemannian manifold with negative sectional curvature, then the eigenfunctions of $-\Delta_g$ are observable from any nonempty open set $O \subset M$ in the sense of Definition~\ref{de_observability_eigenfunctions}.
\end{example}

We now briefly comment on inverse spectral results analogous to Theorem~\ref{thm_spectral} in the setting of Riemannian manifolds with boundary.

\begin{remark}
	\label{rmk_manifold_boundary}
As will be evident from following the proof of Theorem~\ref{thm_spectral}, the theorem may also be stated on compact Riemannian manifolds with boundary under an analogous formulation of the condition {\bf (H)}, taking into account the notion of GCC on Riemannian manifolds with boundary in \cite{BLR}. As an example, if $(M,g)$ is either a Euclidean domain with a smooth boundary or if it is a simple Riemannian manifold and if $O$ is an open neighborhood of $\p M$, then Theorem~\ref{thm_spectral} remains valid. However, for the sake of clarity of proof we have only considered manifolds without boundary in this work. %(replacing condition {\bf (H)} with the latter simplicity assumption and by utilizing the fact that $O$ is a neighborhood of $\p M$). 
\end{remark}

Our second inverse spectral result concerns the reconstruction of a potential from partial spectral data, in the sense that an infinite number of spectral pairs are missing from the dataset. To quantify this, we introduce the following notion from interpolation theory; see \cite{BC89}, \cite{K57}, as well as \cite[Section~1.3]{OU09}. See also Theorem~\ref{thm_PW_interpolation} below.

\begin{definition}
\label{def_u_d}
A set $\Gamma = \{\gamma_k\}_{k=1}^\infty \subset \mathbb{R}$ is called \emph{uniformly discrete} if
\[
\inf_{j \ne k} |\gamma_j - \gamma_k| > 0.
\]
The \emph{upper uniform density} of a uniformly discrete set $\Gamma$ is defined by
\[
D^+(\Gamma) := \lim_{l \to \infty} \sup_{s \in \mathbb{R}} \frac{\#(\Gamma \cap (s, s + l))}{l}.
\]
\end{definition}

We also introduce the following notion.
\begin{definition}
\label{def_adm}
Let $\Lambda = \{\lambda_k\}_{k=1}^\infty \subset \mathbb{R}$ with $\lambda_1 \leq \lambda_2 \leq \ldots$ and $\lim\limits_{k \to \infty} \lambda_k = +\infty$. A proper subset $\widetilde{\Lambda} \subset \Lambda$ is said to be \textbf{$\Lambda$-sparse} if either $\widetilde{\Lambda} = \emptyset$, or, if $\widetilde{\Lambda} \ne \emptyset$, then the set
\[
\Gamma = \left\{
\begin{array}{l}
\sqrt{\lambda_j - \lambda_1 + 1} - \sqrt{\lambda_k - \lambda_1 + 1}, \\
\pm\left(\sqrt{\lambda_j - \lambda_1 + 1} + \sqrt{\lambda_k - \lambda_1 + 1} \,\right)
\end{array}
: \lambda_j, \lambda_k \in \widetilde{\Lambda} \right\}
\]
is uniformly discrete and has zero upper uniform density, i.e., $D^+(\Gamma) = 0$.
\end{definition}

Let us remark that throughout this paper, the set $\Lambda$ will essentially denote the spectrum of the operator $-\Delta_g + V$ on a Riemannian manifold. Our notion of $\Lambda$-sparsity serves as a quantitative measure of how much spectral data, consisting of eigenvalues and the restrictions to $O$ of the corresponding eigenspaces, can be entirely missing from the dataset while still allowing for unique recovery of the potential $V$ on $M$. Importantly, the definition of sparsity depends only on the values of the eigenvalues, and \textbf{not} on their multiplicities. This has the following implication: whenever an eigenvalue is missing from the dataset, we treat the entire associated eigenspace (restricted to $O$) as missing. Conversely, if an eigenvalue is present, we require that the full information of its eigenspace (restricted to $O$) is available in the dataset.

In addition to the $\Lambda$-sparsity condition, we also require a uniform observability property on a nonempty open set $O \subset M$:

\begin{enumerate}
  \item[\textbf{(UO)}] Let $V \in C^\infty(M)$ be real-valued. Then the eigenfunctions of the operator $-\Delta_g + V$ are observable from every nonempty open subset of $O$ (including $O$ itself), in the sense of Definition~\ref{de_observability_eigenfunctions}.
\end{enumerate}

We now state our second inverse spectral result, in which the identification of eigenpairs is restricted to subsequences, and the spectral data for both operators may be incomplete, allowing countably many eigenvalues and eigenfunctions to be entirely missing.

\begin{theorem}
\label{thm_anosov}
Let $(M, g)$ be a smooth, closed, connected Riemannian manifold of dimension $n \geq 2$, and let $O \subset M$ be a nonempty, open, and connected subset with smooth boundary such that $M \setminus \overline{O} \neq \emptyset$. Suppose that there exists a point $p \in O$ such that $\mathcal{A}_{M, g}(p) \subset O$.

For $j = 1, 2$, let $V_j \in C^\infty(M)$ be real-valued potentials such that assumption~\textbf{(UO)} holds for the set $O$ and both potentials $V_1$ and $V_2$. Denote by
\[
\mu_1^{(j)} \leq \mu_2^{(j)} \leq \cdots
\]
the eigenvalues of the Schr\"odinger operator $-\Delta_g + V_j$ on $M$, listed in increasing order and counted with multiplicity. Let $\{\psi_k^{(j)}\}_{k=1}^{\infty} \subset C^{\infty}(M)$ be a linearly independent set of corresponding eigenfunctions satisfying
\[
(-\Delta_g + V_j)\, \psi_k^{(j)} = \mu_k^{(j)}\, \psi_k^{(j)} \quad \text{on } M.
\]

Let $\mathbb{A} = \{a_k\}_{k=1}^\infty \subset \mathbb{N}$ be a nonempty subset with $a_1 < a_2 < \cdots$, and suppose that the set $\{\mu_k^{(1)} : k \in \mathbb{N} \setminus \mathbb{A} \}$ is $\{\mu_k^{(1)}\}_{k=1}^\infty$-sparse in the sense of Definition~\ref{def_adm}. Assume further that there exists a strictly increasing sequence $\{b_k\}_{k=1}^\infty \subset \mathbb{N}$ such that
\begin{equation}
\label{psi_eq_1}
\mu_{a_k}^{(1)} = \mu_{b_k}^{(2)} \quad \text{and} \quad \psi_{a_k}^{(1)}(x) = \psi_{b_k}^{(2)}(x) \quad \text{for all } x \in O \text{ and all } k \in \mathbb{N}.
\end{equation}
Then $V_1 = V_2$ on $M$.
\end{theorem}

\begin{remark}
\label{rem_extension_basis_2}
As in Theorem~\ref{thm_spectral}, the assumption in Theorem~\ref{thm_anosov} that the full system of eigenfunctions $\{\psi_k^{(j)}\}_{k = 1}^\infty$ of the operator $-\Delta_g + V_j $, for $ j = 1, 2$, is linearly independent can be relaxed. It suffices to assume linear independence only for the subsets  $\{\psi_{a_k}^{(1)}\}_{k=1}^\infty$ and $\{\psi_{b_k}^{(2)}\}_{k=1}^\infty$ of eigenfunctions appearing in~\eqref{psi_eq_1}; see Remark~\ref{rem_extension_basis}.
\end{remark}

\begin{remark}
Although Theorem~\ref{thm_anosov} requires the uniform observability condition on a nonempty open set $O \subset M$ to hold for both potentials $V_1$ and $V_2$, this assumption is justified by explicit examples. For instance, it is satisfied when $(M, g)$ is either a flat torus of arbitrary dimension or an Anosov surface, and $O$ is any nonempty open subset, for all smooth potentials; see Examples~\ref{eq_ex_torus} and~\ref{eq_ex_Anosov}, and the references given there.
\end{remark}

\begin{remark}
In Theorem~\ref{thm_anosov}, the observation set $O$ is required to satisfy the uniform observability condition, which is a much stronger assumption than the observability of eigenfunctions of $-\Delta_g + V$ on $O$. The latter condition is sufficient for Theorem~\ref{thm_spectral}. However, the inverse spectral result in Theorem~\ref{thm_anosov} is also stronger, as it allows for a countably infinite subset of the spectral data to be missing for both operators.
\end{remark}

Next, we present a simple example of a smooth, closed, connected Riemannian manifold and an observation set that together satisfy all the geometric assumptions of Theorem~\ref{thm_anosov}.

\begin{example}
Let $(M, g)$ be the flat torus $\mathbb{T}^n = (\mathbb{R} / \mathbb{Z})^n$ equipped with the standard metric, and let $V_1, V_2 \in C^{\infty}(\mathbb{T}^n)$ be real-valued potentials. Then the eigenfunctions of the Schr\"odinger operators $-\Delta_g + V_j$, for $j = 1, 2$, are observable from any nonempty open set $O \subset \mathbb{T}^n$. Consider the point $p = [(0, \ldots, 0)] \in \mathbb{T}^n$, and observe that
\[
q = \left[\left(\tfrac{1}{2}, \ldots, \tfrac{1}{2}\right)\right] \in \mathcal{A}_{\mathbb{T}^n, g}(p).
\]
Let $O \subset \mathbb{T}^n$ be a connected open set that contains neighborhoods of both $p$ and $q$. Then $(M, g)$ and $O$ satisfy the assumptions of Theorem~\ref{thm_anosov}.
\end{example}

We now show that on any smooth, closed Riemannian manifold, it is possible to construct an infinite subsequence $\{a_k\}_{k=1}^\infty$ such that the set $\{\mu^{(1)}_{a_k}\}_{k=1}^\infty$ is $\{\mu_k^{(1)}\}_{k=1}^\infty$-sparse in the sense of Definition~\ref{def_adm}, as considered in Theorem~\ref{thm_anosov}.

\begin{example}
Let $(M, g)$ be a smooth, closed, connected Riemannian manifold of dimension $n \geq 2$, and let $V_1 \in C^\infty(M)$ be a real-valued potential. Let $\mu_1^{(1)} \leq \mu_2^{(1)} \leq \cdots$ denote the eigenvalues of the Schrödinger operator $-\Delta_g + V_1$ on $M$, listed in non-decreasing order and counted with multiplicity. By Weyl's law, there exist constants $k_0 \in \mathbb{N}$ and $c > 1$, depending only on $(M, g)$ and $V_1$, such that
\[
c^{-1} \mu^{n/2} \leq N(\mu) \leq c \mu^{n/2} \quad \text{for all } \mu > \mu_{k_0},
\]
where $N(\mu)$ denotes the number of eigenvalues of $-\Delta_g + V_1$ less than or equal to $\mu$, counted with multiplicity; see~\cite[Theorem 15.2, p.~130]{Shubin_book}. Let $A \gg 1$ be a sufficiently large constant, depending only on $c$. In view of the two-sided Weyl bound above, it is straightforward to see that there exists an eigenvalue in each interval $(A^{2k},A^{2k+1}]$ for each $k>k_1$ where $k_1>k_0$ is a sufficiently large integer. To construct our sparse set $\{\mu_{a_k}^{(1)}\}_{k=1}^\infty$, we choose exactly one eigenvalue from each of the above intervals (we can include all of the eigenvalues that are equal to it according to the multiplicity as multiplicity of an eigenvalue is irrelevant in the notion of sparsity). The resulting sequence $\{\mu_{a_k}^{(1)}\}_{k=1}^\infty$ is then $\{\mu_k^{(1)}\}_{k=1}^\infty$-sparse.
\end{example}

Finally, in this section, we outline the main ideas behind the proofs of Theorem~\ref{thm_spectral} and Theorem~\ref{thm_anosov}. Both arguments rely on establishing a connection between the non-normalized spectral data and the corresponding wave equations. This connection allows us to leverage key properties of the wave equation, such as finite speed of propagation, energy estimates, and a global version of Tataru’s unique continuation theorem (see \cite[Theorem~3.24]{LL23}, \cite{KKL}, and also \cite{Tataru}). This is precisely why the antipodal set plays a crucial role and appears in both assumption~\textbf{(H)} and the assumptions of Theorem~\ref{thm_anosov}. Another essential ingredient enabling this reduction to the wave setting is the observability estimate~\eqref{eq_observability_eigenfunctions_def} for eigenfunctions of the Schr\"odinger operator. The proof of Theorem~\ref{thm_anosov} further relies on a deep Paley--Wiener type interpolation theorem due to Beurling~\cite{BC89} and Kahane~\cite{K57}, which is instrumental in allowing for a countably infinite subset of the spectral data to be missing for both operators.

\subsection{Previous literature on {\bf (Q)}} 

The majority of the literature in the multi-dimensional setting addresses a well-known, simplified variant of question~\textbf{(Q)}, commonly 
known as the Gel'fand inverse spectral problem, in which global norming constants for the eigenfunctions are also assumed to be given to us.
Returning to the example of "hearing the shape of the drum", Gel'fand's inverse problem assumes that we can actively hit the drum with arbitrary a priori known intensities at a certain subset of its surface while simultaneously measuring the ensuing vibrations there.  For Schr\"odinger operators on bounded domains in Euclidean space, this form of the inverse spectral problem, recovering a potential from the Dirichlet eigenvalues and the boundary traces of the normal derivatives of the normalized eigenfunctions, was resolved in the works \cite{Nachman_Sylveser_Uhlmann_88} and \cite{Nov_1988}.

In the geometric setting of the Laplace--Beltrami operator (or more general second-order invariant differential operators) on a Riemannian manifold $(M, g)$, Gel'fand’s problem asks whether both the operator and the isometry class of $M$ can be determined from the knowledge of the spectrum and an $L^2(M)$-orthonormal basis of eigenfunctions, measured either on an open interior subset $O \subset M$ or on the boundary of $M$. The availability of a complete orthonormal system in this context enables a reduction to an associated inverse problem for the wave equation.

This setting has been extensively studied; for example, in the context of a closed manifold, the spectral data can be formulated as the eigenvalues and the restrictions to $O$ of an $L^2(M)$-orthonormal basis of eigenfunctions of $-\Delta_g$; see \cite{BKL, KrKaLa, HLOS_2018, Lu_Jinpeng_2025}. Inverse spectral problems of this kind are foundational in the study of inverse problems for evolution equations with active measurements, where one can initiate controlled experiments and observe the system’s response. A prominent application is geophysical exploration, in which artificial seismic or electromagnetic signals are introduced into the Earth and the reflected waves are recorded at the surface, with the goal of reconstructing internal structures. Analogous inverse problems are studied for the heat equation, Schr\"odinger equations, and other evolutionary PDEs. A comprehensive treatment of such inverse problems on Riemannian manifolds with boundary is provided in the monograph~\cite{KKL}. Remarkably, all these problems have been shown to be equivalent in the case of manifolds with boundary~\cite{KKLM}.

Gel'fand’s inverse spectral problem for the acoustic operator $-c^2(x)\Delta$ was first solved in the pioneering work of Belishev~\cite{Bel1} for domains in $\mathbb{R}^n$, $n \geq 2$, and later extended to general Riemannian manifolds by Belishev and Kurylev~\cite{BK92}. Their solution is based on the development of a boundary-based control theory for the wave equation, known as the Boundary Control (BC) method. A crucial ingredient in this approach is a sharp unique continuation theorem for wave equations due to Tataru~\cite{Tataru} which was followed by related contributions from Robbiano and Zuily~\cite{RZ} as well as further developments by Tataru~\cite{Tataru2}. We refer to~\cite{Bel2, KKL} for comprehensive expositions of the BC method. 

The work~\cite{Isozaki91} provides a solution to Gel'fand’s inverse spectral problem in the Euclidean setting when a finite number of spectral data are missing, under the assumption that the eigenfunctions form an $L^2$-orthonormal basis. See also~\cite{Katchalov_Kurylev_1998} for the corresponding result in the case of Riemannian manifolds with boundary. A solution to the inverse spectral problem with asymptotically close spectral data was established in~\cite{Choulli_Stefanov_2013} and~\cite{Bellassoued_Choulli_DF_Kian_Stefanov_2021}, under the assumption that the eigenfunctions are globally orthonormalized in $L^2$. For recent developments and additional references concerning classical Gel'fand inverse spectral problems, we refer the reader to~\cite{Liu_Quan_Saksala_Yan_2025, Lu_Jinpeng_2025}.

We conclude this section by emphasizing that, to the best of our knowledge, the result~\cite[Theorem 1.11]{FKU24} remains the only available affirmative answer to question~\textbf{(Q)} in the challenging case where the Riemannian manifold is unknown and must be reconstructed from spectral data without any a priori normalization of the eigenfunctions.

Furthermore, to the best of our knowledge, in the absence of global norming constants, no prior results exist in dimensions two or higher for the case when part of the spectral data, whether finite or infinite, is missing. The present work advances the theory by addressing this question in a different important setting: the recovery of a potential in the stationary Schr\"odinger operator, under the assumption that the manifold is known, but again without assuming any orthogonality or normalization of the eigenfunctions.

\subsection{Applications to inverse problems with a single passive measurement}
\label{sec_ip_passive}

We will now discuss several applications of our inverse spectral results in relation to inverse problems for PDEs with a single passive measurement. Broadly speaking, the motivation for such inverse problems stems from real-world situations in which the input cannot be controlled, but the system's response can be observed—e.g., medical imaging (e.g., EEG, fMRI), geophysical exploration, or astrophysics. These are problems characterized by limited data: the source may be unknown or partially known, and the measurements may be confined to a fixed location or a small number of detectors. Mathematically, such problems fall outside the classical framework of inverse problems with active measurements, where multiple (often infinite) inputs and abundant data allow for more comprehensive reconstruction. Instead, one is confronted with the challenge of extracting maximal geometric or analytic information from minimal (and often noisy) data. This raises fundamental questions: What can be recovered from a single passive measurement? Under what conditions are uniqueness or stability results possible? And how are such constraints governed by the underlying physics? 

A prototypical model for such an inverse problem may be described as follows; suppose that $u$ satisfies a certain PDE with unknown coefficients, governing a physical law in nature, subject to some unknown source or initial data. Can one uniquely determine the coefficients as well as the source or initial data from the measurement of $u$ on some subset of the spacetime? Due to the highly limited nature of the data, the inverse problem described above is notoriously challenging, and most existing results rely on rather restrictive assumptions. For instance, some works assume that the coefficients in the PDE are all a priori known, and focus solely on recovering the source or initial data; such problems are known as inverse source problems, see e.g., \cite{BLMM19,OU13,SU09, Liu_Uhlmann_Yang_2024}. In contrast, other results consider the recovery of an unknown coefficient in the PDE from a single passive measurement, but under the assumption that the source or initial data is a priori known and specially designed (often highly singular); see e.g., \cite{BK81,FK23,HLO13,HLOS_2018,IY24,KLLY20,SU11}. Within the latter framework and starting from the works \cite{K95,K97}, Klibanov and his coauthors have developed a convexification method via Carleman estimates that yields suitable numerically convergent algorithms  based on the stability properties of the continuum problem and allowing construction of numerical methods which would not be affected by the problem of local minima in such inverse problems, see \cite{KL21} for more details. 

In the more general setting, where both the initial data or source and the coefficients in the PDE are all unknown, the available results in multidimensions typically impose further structural constraints. These include assumptions that the unknowns satisfy explicit equations \cite{KM20,LU15} or that they obey certain monotonicity properties \cite{KU25} or that they are finite dimensional \cite{Kian_H_Liu_2025}. It should be noted that in the one dimensional setting, there are optimal results that solve the problem for generic classes of initial data \cite{AMR14,Pie79,Suz86}, see for instance \cite[Theorem 0]{Suz86} that shows that the genericity assumption on the initial data is indeed optimal. We also refer the reader to recent work of the first author in \cite{Fei25} that solves inverse 1-D problems subject to a passive measurement under some support conditions for the initial data \cite{Fei25}. The latter 1-D results are all based on the comprehensive understanding of inverse Sturm-Liouville theory that was developed in the 1940s, see e.g. \cite{Borg_1946,HL78,Gelfand_Levitan_1951,GS99,LG64}. One of the main goals in this paper is to show that generalizations of the generic (optimal) 1-D results to multidimensions is indeed possible, thanks to our multidimensional inverse spectral results. We will illustrate this by considering several evolutionary models. Whether or not these generic assumptions are optimal in the multidimensional setup is not known to us.

For our first example, we study heat equations. Let $V \in C^\infty(M)$ and $f \in C^\infty(M)$. Consider the initial value problem for the heat equation,
\begin{equation}
\label{int_heat_1}
\begin{cases}
\partial_t u(t,x) - \Delta_g u(t,x) + V(x) u(t,x) = 0, & \text{on } (0,\infty) \times M, \\
u(0,x) = f(x), & \text{on } M.
\end{cases}
\end{equation}
It is well known that the problem \eqref{int_heat_1} admits a unique solution $u = u^f \in C^\infty([0,\infty) \times M)$. The existence of the solution follows from its spectral representation based on the eigenvalues and eigenfunctions of the self-adjoint operator $-\Delta_g + V$; see~\eqref{eq_11_1}. Uniqueness follows from standard energy estimates; see \cite[Chapter 6, pp.~482–483]{Taylor_book_I}.

Our first result concerns the inverse problem of simultaneously determining the unknown potential $V$ and the initial data $f$ on $M$ from a single measurement of the solution $u^f$ to \eqref{int_heat_1} on the set $(0,\varepsilon) \times O$, where $\varepsilon > 0$ is arbitrary and $O \subset M$ is open.
\begin{theorem}
\label{thm_main_heat}
Let $(M,g)$ be a smooth, closed, and connected Riemannian manifold of dimension $n \geq 2$, and let $O \subset M$ be a nonempty, open, and connected subset such that $M \setminus \overline{O} \neq \emptyset$. Suppose that assumption~\textbf{(H)} is satisfied. Let $f_1, f_2 \in C^\infty(M)$ and $V_1, V_2 \in C^\infty(M)$. 
Assume that there exists $k_0 \in \mathbb{N}$ such that the eigenvalues $\mu_k^{(1)}$ of the operator $-\Delta_g + V_1$ are simple for all $k \geq k_0$, and that $f_1$ has nontrivial Fourier coefficients with respect to the $L^2(M)$-orthonormal basis $\{\psi_k^{(1)}\}$ of eigenfunctions corresponding to $\mu_k^{(1)}$, meaning
\begin{equation}
\label{int_heat_2_1}
(f_1, \psi_k^{(1)})_{L^2(M)} \ne 0 \quad \text{for all } k \geq k_0.
\end{equation}
Assume furthermore that for some $\varepsilon > 0$,
\begin{equation}
\label{int_heat_3}
u_1^{f_1}(t,x) = u_2^{f_2}(t,x), \quad \text{for all } (t,x) \in (0,\varepsilon) \times O,
\end{equation}
where $u_j^{f_j} \in C^\infty([0,\infty) \times M)$ is the solution to \eqref{int_heat_1} with $V = V_j$ and $f = f_j$, for $j = 1,2$. Then $V_1 = V_2$ and $f_1 = f_2$ on $M$.
\end{theorem}

\begin{remark}
It is a generic property that $C^\infty$ potential perturbations of the Laplacian $-\Delta_g$, associated with a $C^\infty$ Riemannian metric on a compact manifold, have only simple eigenvalues; see \cite{Uhlenbeck_1976, Albert_1975}.  Moreover, the property that the initial data $f_1$ has nontrivial Fourier coefficients with respect to the $L^2(M)$-orthonormal basis of eigenfunctions $\{\psi_k^{(1)}\}$ of the operator $-\Delta_g + V_1$ is also generic. Indeed, the set
\[
U = \left\{ f \in L^2(M) : (f, \psi_k^{(1)})_{L^2(M)} \ne 0 \quad \text{for all } k \in \mathbb{N} \right\}
\]
is residual in $L^2(M)$, since it can be written as a countable intersection $U = \bigcap_{k=1}^\infty U_k$, where
\[
U_k = \left\{ f \in L^2(M) : (f, \psi_k^{(1)})_{L^2(M)} \ne 0 \right\}.
\]
Each $U_k$ is open and dense in $L^2(M)$ for $k = 1, 2, \dots$.
\end{remark}

\begin{remark}
We note that the generic assumption~\eqref{int_heat_2_1} on the initial data, which appears in Theorem~\ref{thm_main_heat}, is necessary even in the one-dimensional case; see~\cite{Suz86, Murayama_1981}. Moreover, to the best of our knowledge, all previous works in one dimension on inverse problems with passive measurements have assumed that the spectra of both operators are simple.
\end{remark}

Our next result addresses the inverse problem of simultaneously determining the unknown Riemannian manifold $(M, g)$ up to an isometry, and the initial data $f$ up to the corresponding gauge transformation, from a single measurement of the solution $u^f$ to \eqref{int_heat_1} with $V = 0$, observed on the set $(0, \varepsilon) \times O$, where $\varepsilon > 0$ is arbitrary and $O \subset M$ is open.

\begin{theorem}
\label{thm_main_heat_metric}
Let $(M_j, g_j)$ be smooth, closed, and connected Riemannian manifolds of dimension $n \geq 2$, for $j = 1, 2$. Let $O \subset M_1 \cap M_2$ be a nonempty open set such that $M_j \setminus \overline{O} \neq \emptyset$ for $j = 1, 2$, and suppose that assumption~\textbf{(H)} is satisfied for both $(M_j, g_j)$ and $O$, and that $g_1|_O = g_2|_O$.
Assume that all eigenvalues of the Laplace--Beltrami operators $-\Delta_{g_j}$ are simple, and denote them by
\begin{equation}
\label{int_heat_2_1_metric_eigen}
0 = \mu_0^{(j)} < \mu_1^{(j)} < \mu_2^{(j)} < \dots,
\end{equation}
for $j = 1, 2$. Let $f_j \in C^\infty(M_j)$, and suppose that $f_j$ has nontrivial Fourier coefficients with respect to the $L^2(M_j)$-orthonormal basis $\{\psi_k^{(j)}\}_{k=0}^\infty$ of eigenfunctions corresponding to $\mu_k^{(j)}$, that is,
\begin{equation}
\label{int_heat_2_1_metric}
(f_j, \psi_k^{(j)})_{L^2(M_j)} \neq 0 \quad \text{for all } k = 0, 1, 2, \dots,\quad j=1,2.
\end{equation}

Assume furthermore that for some $\varepsilon > 0$,
\begin{equation}
\label{int_heat_3_metric}
u_1^{f_1}(t,x) = u_2^{f_2}(t,x) \quad \text{for all } (t,x) \in (0, \varepsilon) \times O,
\end{equation}
where $u_j^{f_j} \in C^\infty([0, \infty) \times M_j)$ is the solution to \eqref{int_heat_1} with $g = g_j$, $V = 0$, and $f = f_j$, for $j = 1, 2$. Then there exists a smooth diffeomorphism $\Phi : M_1 \to M_2$ such that
\[
\Phi|_O = \mathrm{Id}, \quad g_1 = \Phi^* g_2, \quad \text{and} \quad f_1 = f_2 \circ \Phi \quad \text{on } M_1.
\]
\end{theorem}

\begin{remark} 
The simplicity of the spectrum is a generic property of the Laplacian $-\Delta_g$ associated with a $C^\infty$ Riemannian metric; see \cite{Uhlenbeck_1976}, \cite{Tanikawa_1979}, and \cite{Bando_Urakawa_1983}.
\end{remark}

\begin{remark}
Note that in Theorem~\ref{thm_main_heat_metric}, the simplicity of the spectra of both Laplacians $-\Delta_{g_j}$ and the condition that both initial data $f_j$ for the heat equation \eqref{int_heat_1}, $j = 1, 2$, have nontrivial Fourier coefficients are assumed. This contrasts with the stronger result of Theorem~\ref{thm_main_heat}, where such generic assumptions are imposed only on the first operator. The need to impose these generic conditions on both operators in Theorem~\ref{thm_main_heat_metric} arises from the requirement to deduce, from the measurement data \eqref{int_heat_3_metric}, the equality of all eigenvalues $\mu_k^{(1)} = \mu_k^{(2)}$, $k = 0, 1, 2, \dots$, of the Laplacians $-\Delta_{g_1}$ and $-\Delta_{g_2}$, as well as the equality on $O$ of the corresponding eigenfunctions, which are not necessarily orthonormalized but are linearly independent. This step relies on standard techniques involving the spectral representation of solutions to the heat equation and the use of Laplace transform methods; see \cite{KKL, KKLM}. The recovery of all eigenvalues and the associated non-orthonormalized eigenfunctions on $O$ is crucial for applying \cite[Theorem 1.11]{FKU24}, which enables the reconstruction of Riemannian manifolds up to isometry in Theorem~\ref{thm_main_heat_metric}. We emphasize that, in the case where the Riemannian manifolds are unknown, \cite[Theorem 1.11]{FKU24} provides the strongest available inverse spectral result in the setting where the eigenfunctions in the spectral data are not orthonormalized.
\end{remark}

\begin{remark}
In the case where the Riemannian manifold is known and the goal is to recover the potential from spectral data measured on the observation domain $O$, our Theorem~\ref{thm_spectral} does not require knowledge of all eigenvalues and the associated non-orthonormalized eigenfunctions on $O$. Specifically, we require only the equality \eqref{psi_eq}, where the spectral data for the first operator $-\Delta_g + V_1$ may correspond to a subset of the spectral data for the second operator. This relaxation in the known-manifold setting is precisely what enables the stronger results for inverse problems with a single passive measurement, in which one simultaneously recovers the potential and the initial data; see Theorems~\ref{thm_main_heat}, \ref{thm_main_Schrodinger}, and \ref{thm_main_wave}. Notably, in these theorems, the generic assumptions are imposed only on the first operator and its corresponding initial data.
\end{remark}

\begin{remark}
We note that results, similar to those in Theorem~\ref{thm_main_heat_metric}, for the simultaneous recovery of a potential and the initial data in equation \eqref{int_heat_1} can be derived from known results in inverse spectral theory, provided that genericity assumptions are imposed on both operators and both initial data, and that certain additional assumptions are made on the Riemannian manifold. For instance, if $(M, g)$ is a closed, negatively curved Riemannian manifold with simple length spectrum, and if both operators $-\Delta_g + V_j$ have simple spectra, and both initial data $f_j$, $j = 1, 2$, have nontrivial Fourier coefficients, then, relying on results from \cite{Guillemin_Kazhdan_1980, Croke_Sharafutdinov_1998}, it follows that the measurement data \eqref{int_heat_3} determines $V_1 = V_2$ and $f_1 = f_2$. We emphasize, however, that Theorems~\ref{thm_main_heat}, \ref{thm_main_Schrodinger}, and \ref{thm_main_wave} establish significantly stronger results, as they require the generic assumptions only for the first operator and its corresponding initial data. This improvement is made possible by Theorem~\ref{thm_spectral}.
\end{remark}

\begin{remark}
The generic assumptions on the first operator and its corresponding initial data in Theorem~\ref{thm_main_heat} appear to be necessary for our approach. If the multiplicity of an eigenvalue $\mu_k^{(1)}$ for some $k \geq k_0$ is greater than one, then working with only a single fixed initial data does not allow one to recover a basis of the corresponding eigenspace, which is required in order to apply Theorem~\ref{thm_spectral}.
\end{remark}

Our next result concerns inverse problems in quantum mechanics. Recall that the evolution of a particle under the influence of an external potential $V$ is described by the time-dependent Schr\"odinger equation
\begin{equation}
\label{schrodinger_ip_pf}
\begin{aligned}
\begin{cases}
\textrm{i}\,\partial_t u(t,x) = -\Delta_g u(t,x) + V(x)u(t,x)
& \text{on } (0,\infty) \times M, \\
u(0,x) = f(x) & \text{on } M,
\end{cases}
\end{aligned}
\end{equation}
where $u(t,x)$ is the wave function of the system, and $M$ is the spatial domain of interest. Let $V \in C^\infty(M)$ and $f \in C^\infty(M)$. It is well known that the problem \eqref{schrodinger_ip_pf} admits a unique solution $u = u^f \in C^\infty([0,\infty) \times M)$.

As described in \cite{LS11}, it is possible to measure the wave function on an open \emph{observation} subset  $O \subset M$ of the spatial domain using physical experiments, thus fitting again into our framework of an inverse problem with a single passive measurement. Similar to the case of the heat equation, we study the inverse problem of simultaneously determining the unknown potential $V$ and the initial quantum state $f$ on $M$ from a single measurement of the wave function on the set $(0,\infty)\times O$. That is, we aim to invert the map
\[
(f, V) \quad \mapsto \quad u^f \quad \text{on } (0,\infty) \times O,
\]
where $u^f$ is the solution to \eqref{schrodinger_ip_pf}.

Our main result for this model is as follows.

\begin{theorem}
\label{thm_main_Schrodinger}
Let $(M,g)$ be a smooth, closed, and connected Riemannian manifold of dimension $n \geq 2$, and let $O \subset M$ be a nonempty, open, and connected subset such that $M \setminus \overline{O} \neq \emptyset$. Suppose that assumption~\textbf{(H)} is satisfied. Let $f_1, f_2 \in C^\infty(M)$ and $V_1, V_2 \in C^\infty(M)$. 
Assume that there exists $k_0 \in \mathbb{N}$ such that the eigenvalues $\mu_k^{(1)}$ of the operator $-\Delta_g + V_1$ are simple for all $k \geq k_0$, and that $f_1$ has nontrivial Fourier coefficients with respect to the $L^2(M)$-orthonormal basis $\{\psi_k^{(1)}\}$ of eigenfunctions corresponding to $\mu_k^{(1)}$, that is,
\begin{equation}
\label{int_schrodinder_2_1}
(f_1, \psi_k^{(1)})_{L^2(M)} \ne 0 \quad \text{for all } k \geq k_0.
\end{equation}
Assume furthermore that
\begin{equation}
\label{int_schrodinder_3}
u_1^{f_1}(t,x) = u_2^{f_2}(t,x), \quad \text{for all } (t,x) \in (0,\infty) \times O,
\end{equation}
where $u_j^{f_j} \in C^\infty([0,\infty) \times M)$ denotes the solution to \eqref{schrodinger_ip_pf} with $V = V_j$ and $f = f_j$, for $j = 1,2$. Then $V_1 = V_2$ and $f_1 = f_2$ on $M$.
\end{theorem}

Another application we discuss arises in seismology, where the goal is to probe the Earth's interior using seismic waves generated by natural events such as earthquakes. To be precise, let 
$V \in C^\infty(M)$, $f \in C^\infty(M)$, and $h \in C^\infty(M)$, and consider the initial value problem for the wave equation:
\begin{equation}
	\label{wave_ip_pf}
	\begin{aligned}
		\begin{cases}
			\partial_t^2 u(t,x) - \Delta_g u(t,x) + V(x)u(t,x) = 0
			& \text{on } (0,\infty) \times M, \\
			u(0,x) = f(x) & \text{on } M, \\
			\partial_t u(0,x) = h(x) & \text{on } M.
		\end{cases}
	\end{aligned}
\end{equation}
It is well known that the problem \eqref{wave_ip_pf} admits a unique solution $u = u^{f,h} \in C^\infty([0,\infty) \times M)$; see \cite[Theorem 2.13]{Saksala_Shedlock_2025} and \cite[Chapter 6, pp.\ 485–486]{Taylor_book_I} and \cite[Section 2.3]{KKL}.

We study the inverse problem of simultaneously determining the unknown potential $V$ and the initial data $f$ and $h$ on $M$ from a single measurement of the solution $u^{f,h}$ on the subset $(0,\infty) \times O$. That is, we aim to invert the map
\[
(f, h, V) \quad \mapsto \quad u^{f,h} \quad \text{on } (0,\infty) \times O.
\]

Our main result for the wave model is as follows.

\begin{theorem}
\label{thm_main_wave}
Let $(M,g)$ be a smooth, closed, and connected Riemannian manifold of dimension $n \geq 2$, and let $O \subset M$ be a nonempty, open, and connected subset such that $M \setminus \overline{O} \neq \emptyset$. Suppose that assumption~\textbf{(H)} is satisfied. Let $f_1, f_2 \in C^\infty(M)$, $h_1, h_2 \in C^\infty(M)$, and $V_1, V_2 \in C^\infty(M)$. 
Assume that there exists $k_0 \in \mathbb{N}$ such that the eigenvalues $\mu_k^{(1)}$ of the operator $-\Delta_g + V_1$ are simple for all $k \geq k_0$, and that either $f_1$ or $h_1$ has nontrivial Fourier coefficients with respect to the $L^2(M)$-orthonormal basis $\{\psi_k^{(1)}\}$ of eigenfunctions corresponding to $\mu_k^{(1)}$. That is,
\begin{equation}
\label{int_wave_2_1}
(f_1, \psi_k^{(1)})_{L^2(M)} \ne 0 \quad \text{for all } k \geq k_0, \quad \text{or} \quad (h_1, \psi_k^{(1)})_{L^2(M)} \ne 0 \quad \text{for all } k \geq k_0.
\end{equation}
Assume furthermore that
\begin{equation}
\label{int_wave_3}
u_1^{f_1, h_1}(t,x) = u_2^{f_2, h_2}(t,x), \quad \text{for all } (t,x) \in (0,\infty) \times O,
\end{equation}
where $u_j^{f_j, h_j} \in C^\infty([0,\infty) \times M)$ denotes the solution to \eqref{wave_ip_pf} with $V = V_j$, $f = f_j$, and $h = h_j$, for $j = 1, 2$. Then $V_1 = V_2$, $f_1 = f_2$, and $h_1 = h_2$ on $M$.
\end{theorem}

\begin{remark}
As in Theorem~\ref{thm_main_heat_metric} for the initial value problem for the heat equation, one can similarly recover the unknown Riemannian manifold $(M, g)$ up to isometry, as well as the initial data up to the corresponding gauge transformation, for both the time-dependent Schr\"odinger equation \eqref{schrodinger_ip_pf} and the wave equation \eqref{wave_ip_pf} with zero potentials. This recovery is possible from a single measurement of the solution on the set $(0, \infty) \times O$, under generic assumptions on both operators and initial data.  The arguments are analogous to those used in the heat case, and we do not state these results explicitly, as they are similar to Theorem~\ref{thm_main_heat_metric}.
\end{remark}

\subsection{Organization of the paper} 
The remainder of the paper is organized as follows. In Sections~\ref{proof_thm_1} and~\ref{proof_thm_2}, we present the proofs of our two main inverse spectral results, Theorems~\ref{thm_spectral} and~\ref{thm_anosov}, respectively. Section~\ref{proof_passive_measurements} is devoted to inverse problems with single passive measurements and contains the proofs of Theorems~\ref{thm_main_heat},~\ref{thm_main_heat_metric},~\ref{thm_main_Schrodinger}, and~\ref{thm_main_wave}. The Appendix collects auxiliary results used in the proof of Theorem~\ref{thm_spectral}.

\section{Proof of Theorem~\ref{thm_spectral}}

\label{proof_thm_1}

We begin by stating the following known result concerning the observability of eigenfunctions of the Schr\"odinger operator from an open subset $O \subset M$ that satisfies the geometric control condition (GCC) of Rauch--Taylor / Bardos--Lebeau--Rauch \cite{Rauch_Taylor_1974, BLR0, BLR}. This lemma justifies the assumption that $O$ satisfies the GCC in hypothesis~\textbf{(H)} of Theorem~\ref{thm_spectral}.
\begin{lemma}
	\label{lem_nontrapping}
	Let $(M, g)$ be a smooth, closed Riemannian manifold of dimension $n \geq 2$. Let $O \subset M$ be a nonempty open set satisfying the GCC, and let $V \in C^{\infty}(M)$. Then there exists a constant $C > 0$, depending only on $(M, g)$, $O$, and $V$, such that
	$$
	\|\phi\|_{L^2(M)} \leq C \|\phi\|_{L^2(O)},
	$$
	for every eigenfunction $\phi\in C^\infty(M)$ of the operator $-\Delta_g + V$ on $M$.
\end{lemma}

\begin{proof}
The proof is standard and is included here for completeness and the convenience of the reader. We follow \cite{Macia_2021}. Since $O$ satisfies the GCC, there exists a time $T > 0$ such that the Schr\"odinger evolution satisfies an observability estimate from $O$. That is, there exists a constant $C > 0$, depending only on $(M, g)$, $O$, $T$, and $V$, such that for every initial datum $u_0 \in L^2(M)$, the solution to
\begin{equation}
\label{eq_10_1}
\begin{cases}
i \partial_t u(t,x) = -\Delta_g u(t,x) + V(x) u(t,x), & (t,x) \in \mathbb{R} \times M, \\
u(0,x) = u_0(x), & x \in M,
\end{cases}
\end{equation}
satisfies the observability estimate
\begin{equation}
\label{eq_10_2}
\|u_0\|_{L^2(M)}^2 \leq C \int_0^T \|u(t,\cdot)\|_{L^2(O)}^2 \, dt,
\end{equation}
see \cite{Lebeau_1992}; see also \cite[Theorem 1]{Macia_2021}.

Let $\phi \in C^\infty(M)$ be an eigenfunction of $-\Delta_g + V$ on $M$, that is,
\[
(-\Delta_g + V)\phi = \lambda \phi \quad \text{on } M
\]
for some eigenvalue $\lambda \in \mathbb{R}$. Then the function $u(t,x) = e^{-i\lambda t} \phi(x)$ solves \eqref{eq_10_1} with initial condition $u_0 = \phi$. Applying the observability estimate \eqref{eq_10_2} to this solution completes the proof.
\end{proof}

\begin{remark}
\label{rem_nontrapping}
In the proof of Lemma~\ref{lem_nontrapping}, instead of relying on observability estimates for the Schr\"odinger equation, one could use observability estimates for the wave equation; see \cite{Rauch_Taylor_1974}, \cite[Theorem~1.5]{Laurent_Leautaud_2016}, and \cite[Proposition~1.2]{Le_Rousseau_Lebeau_TT} for the use of other energy spaces. We refer to the proof of \cite[Lemma~5.1]{FKU24} for such an argument in the case $V = 0$.
\end{remark}

Let $s \in \mathbb{R}$, and let $H^s(M)$ denote the standard Sobolev space on $M$, equipped with the norm $\|u\|_{H^s(M)}$; see \cite[page 100]{Sogge_book_2014}. We now recall an equivalent characterization of the Sobolev norm that will be used in the proof of Theorem~\ref{thm_spectral}. 

Let $V \in C^\infty(M)$ be real-valued. Then the Schr\"odinger operator $-\Delta_g + V$ is a self-adjoint, elliptic operator on $L^2(M)$ with domain $\mathcal{D}(-\Delta_g + V) = H^2(M)$. Let $\tau > 0$ be such that the operator $P := -\Delta_g + V + \tau\ge I$ in the sense of self-adjoint operators. For any $s \in \mathbb{R}$, the operator $P^{s/2}$ is defined via the functional calculus for self-adjoint operators. According to Seeley's theorem, $P^{s/2}$ belongs to the class $\Psi^s_{\mathrm{cl}}(M)$ of classical pseudodifferential operators of order $s$; see \cite{Seeley_1997} and \cite[Theorem 5.3.1, page 73]{Agranovich_1994}. Consequently, for all $s \in \mathbb{R}$, the operator $P^{s/2}$ defines a continuous isomorphism from $H^s(M)$ to $L^2(M)$, and therefore, the following norm equivalence holds:
\begin{equation}
\label{eq_12_4}
\|u\|_{H^s(M)} \asymp \|P^{s/2} u\|_{L^2(M)}, 
\end{equation}
see \cite[Corollary 5.3.2, page 73]{Agranovich_1994}.

The following remark will be used later in the proof. 
\begin{remark}
	\label{rmk_V_1_normal}
Without loss of generality, we may assume that $\{\psi_k^{(1)}\}_{k=1}^\infty$ forms an $L^2(M)$-orthonormal basis. This can be achieved by applying the Gram--Schmidt orthonormalization procedure to $\{\psi_k^{(1)}\}_{k=N}^\infty$ and simultaneously performing the same operations on $\{\psi_{b_k}^{(2)}\}_{k=N}^\infty$ at each step, thereby preserving the equalities in \eqref{psi_eq}; see \cite[Section 5]{FKU24}. If $N \geq 2$, we complete the orthonormalization by applying Gram--Schmidt to $\{\psi_k^{(1)}\}_{k=1}^{N-1}$ after projecting onto the orthogonal complement of $\mathrm{span}\{\psi_k^{(1)}\}_{k=N}^\infty$. This yields an $L^2(M)$-orthonormal basis consisting of eigenfunctions, with the identities in \eqref{psi_eq} still preserved.
\end{remark}

We will first prove Theorem~\ref{thm_spectral} in the case of $N=1$.

\begin{proof}[Proof of Theorem~\ref{thm_spectral} in the case $N=1$]

The proof in this case is similar to that of \cite[Theorem 1.11]{FKU24}, but requires some modifications due to the fact that, in general, the index set $\{b_k\}_{k=1}^\infty$ may not coincide with $\mathbb{N}$, and also due to the presence of potentials. We introduce the relevant notation and provide a proof of the necessary modifications here, while the remainder of the proof is included in Appendix~\ref{appendix_first_thm}. This is done both because similar arguments are used in the proofs of Theorem~\ref{thm_spectral} for $N > 1$ and Theorem~\ref{thm_anosov}, and for completeness and the reader's convenience.

Let $p \in O$ be as in {\bf (H)} and consider a fixed $q \in \mathcal A_{M,g}(p)\subset O$. As shown in \cite[Section 5]{FKU24}, there exists $\varepsilon_0 \in (0,1)$ such that for all $0 < \varepsilon \le \varepsilon_0$, there holds:
\begin{equation}
\label{eq_condition_P}
\quad \quad {\textbf{(P)}} \quad \text{if $x \in \{y \in M \,:\, \textrm{dist}_{g}(p,y) \geq \textrm{dist}_{g}(p,q) - \varepsilon\}$ then $x \in O$.}
\end{equation}

Next, given any $x\in O$ and any $r>0$ sufficiently small, let $\mathbb B_r(x)$ be the open geodesic ball centred at the point $x\in O$ with radius $r$. For the remainder of this proof, we will fix a $\varepsilon \in (0,1)$ sufficiently small so that 
	\begin{equation}\label{balls_cond} \overline{\mathbb B_{\varepsilon}(q) \cup \mathbb B_{\varepsilon}(p)} \subset O,\end{equation}
	and also so that property {\bf (P)} above is satisfied. Let us define open connected sets $\mathcal V,\mathcal W \subset O$ by 
\begin{equation}
\label{eq_W_V}	
\mathcal W= \mathbb B_{\frac{\varepsilon}{2}}(p) \quad \text{and} \quad  \mathcal V= \mathbb B_{\frac{\varepsilon}{2}}(q).
\end{equation}

Let $f_1 \in C^{\infty}_0(\mathcal{V})$. Recalling from Remark~\ref{rmk_V_1_normal} that $\{\psi^{(1)}_k\}_{k=1}^{\infty}$ may be assumed to form an orthonormal basis for $L^2(M)$, we write
\begin{equation}\label{f_1_exp}
f_1(x) = \sum_{k=1}^\infty a_k \psi_k^{(1)}(x), \quad x \in M,
\end{equation}
where
\begin{equation}\label{a_k_def}
a_k = (f_1, \psi^{(1)}_k)_{L^2(\mathcal{V})} = \int_{\mathcal{V}} f_1(x)\, \overline{\psi_k^{(1)}(x)}\, dV_g = \int_{\mathcal{V}} f_1(x)\, \overline{\psi_{b_k}^{(2)}(x)}\, dV_g,
\end{equation}
and the convergence in \eqref{f_1_exp} is to be understood in the $L^2(M)$ sense. Here $dV_g$ denotes the Riemannian volume form.

Let $\tau > 0$ be such that, for both operators $P_j := -\Delta_g + V_j + \tau$, $j = 1, 2$, we have $P_j \ge I$ in the sense of self-adjoint operators. Since $f_1 \in H^s(M)$ for all $s \in \mathbb{R}$, it follows from \eqref{eq_12_4} with $P = P_1$ that the series in \eqref{f_1_exp} converges in the $H^s(M)$ norm for all $s \in \mathbb{R}$. In particular, by the Sobolev embedding theorem, the series in \eqref{f_1_exp} converges uniformly on $M$.

Since $f_1\in C^\infty_0(\mathcal{V})$, its Fourier coefficients satisfy the estimate
\begin{equation}
\label{eq_11_1_1_proof_1}
|a_k| \le C_L\, (\mu_k^{(1)} + \tau)^{-L}, \quad \text{for all } L \in \mathbb{N},
\end{equation}
where $C_L>0$. 

Let us (for now, formally) define
\begin{equation}
    \label{tilde_f}
    f_2(x) := \sum_{k=1}^{\infty} a_k \,\psi^{(2)}_{b_k}(x), \qquad x \in M,
\end{equation}
where the coefficients $a_k$ are as in \eqref{a_k_def}. We claim that the definition above is well-posed in the sense that the right-hand side converges in the $H^{s}(M)$ norm for any $s \geq 0$, and consequently that $f_2 \in C^{\infty}(M)$.

To justify this, we first observe from \eqref{psi_eq} that
\begin{equation}
    \label{eq_500_1} 
    \|\psi^{(2)}_{b_k}\|_{L^2(O)} = \|\psi^{(1)}_k\|_{L^2(O)} \leq \|\psi^{(1)}_k\|_{L^2(M)} = 1, \qquad k \in \mathbb{N}.
\end{equation}
Applying Lemma~\ref{lem_nontrapping} with $V = V_2$ and using \eqref{eq_500_1}, we obtain
\begin{equation}
    \label{eq_500_2}
    \|\psi^{(2)}_{b_k}\|_{L^2(M)} \leq C, \qquad k \in \mathbb{N},
\end{equation}
for some constant $C > 0$ independent of $k$.

Next, fix $m \in \mathbb{N}$. Using \eqref{eq_12_4} with $P = P_2$ and $s = 2m$, and noting that
\[
P_2^m \psi^{(2)}_{b_k} = (\mu_{b_k}^{(2)} + \tau)^m \psi^{(2)}_{b_k} \quad \text{on } M,
\]
we obtain
\begin{equation}
    \label{eq_500_3}
    \|\psi^{(2)}_{b_k}\|_{H^{2m}(M)} \leq C ( \mu_{b_k}^{(2)} + \tau)^m \|\psi^{(2)}_{b_k}\|_{L^2(M)} 
    \leq C (\mu_{b_k}^{(2)} + \tau)^m, \qquad k \in \mathbb{N},
\end{equation}
where we used \eqref{eq_500_2} in the last inequality.

Interpolating between the bounds \eqref{eq_500_2} and \eqref{eq_500_3}, we obtain that for any $s \ge 0$,
\begin{equation}
    \label{psi_est_1}
    \|\psi^{(2)}_{b_k}\|_{H^s(M)} \leq C (\mu^{(2)}_{b_k} + \tau)^{s/2}, \qquad k \in \mathbb{N},
\end{equation}
for some constant $C > 0$ independent of $k$. 

Recalling that $\mu_k^{(1)} = \mu_{b_k}^{(2)}$ and combining \eqref{psi_est_1} with \eqref{eq_11_1_1_proof_1}, we obtain, for any $s \ge 0$ and any $L > 0$,
\begin{equation}
	\label{eq_500_4}
	\|a_k \psi_{b_k}^{(2)}\|_{H^s(M)} \leq \frac{C}{(\mu_k^{(1)} + \tau)^{L - s/2}}, \qquad k \in \mathbb{N},
\end{equation}
where $C > 0$ is independent of $k$.

Next, we recall the following consequence of Weyl's law for the elliptic operator $P_1\ge I$: there exists a constant $C > 0$ such that
\begin{equation}
\label{eq_11_4}
N(\mu) \le C \mu^{\frac{n}{2}}, \quad \text{for } \mu \gg 1,
\end{equation}
where $N(\mu)$ is the number of eigenvalues of $P_1$ less than or equal to $\mu$, counted with multiplicity; see \cite[Theorem 15.2, p.~130]{Shubin_book}. It follows from \eqref{eq_11_4} that
\begin{equation}
\label{eq_11_5}
\mu_k^{(1)} + \tau \ge C^{-2/n} k^{2/n}, \quad \text{for } k \gg 1.
\end{equation}

Fixing $s \ge 0$ and choosing $L > n/2 + s/2$, we combine \eqref{eq_500_4} and \eqref{eq_11_5} to obtain
\begin{equation}
	\label{eq_500_4_new_a}
	\|a_k \psi_{b_k}^{(2)}\|_{H^s(M)} \leq \frac{C}{k^{\frac{2}{n}(L - s/2)}}, \quad \text{for } k \gg 1,
\end{equation}
where $C > 0$ is independent of $k$. The estimate \eqref{eq_500_4_new_a} shows that the series in \eqref{tilde_f} converges in the $H^s(M)$ norm for any $s \ge 0$. Therefore, the Sobolev embedding theorem implies that $f_2 \in C^\infty(M)$. In particular, the series also converges uniformly on $M$.

Now Lemma~\ref{lemma_app_eq_f} implies that $f_1(x) = f_2(x)$ for all $x \in M$. Therefore, for any function $f \in C^{\infty}_0(\mathcal{V})$, the following holds:
\begin{equation}
	\label{eigen_expansion}
	f(x) = \sum_{k=1}^{\infty} (f, \psi^{(2)}_{b_k})_{L^2(\mathcal{V})} \, \psi^{(2)}_{b_k}(x), \qquad x \in M,
\end{equation}
where we recall that the convergence of the series above holds in any $H^s(M)$ norm for $s \geq 0$.

We claim that equation \eqref{eigen_expansion} implies that $b_k = k$ for all $k \in \mathbb{N}$, which is equivalent to proving that $\{b_k\}_{k=1}^\infty = \mathbb{N}$. Suppose, for the sake of contradiction, that there exists some $m \in \mathbb{N} \setminus \{b_k\}_{k=1}^\infty$. 
There are two cases to consider. 

\textbf{Case 1:} $\mu_m^{(2)} \neq \mu_{b_k}^{(2)}$ for all $k \in \mathbb{N}$. In this case, since any two eigenfunctions of $-\Delta_g + V_2$ corresponding to distinct eigenvalues are orthogonal, it follows from \eqref{eigen_expansion} that 
\[
(f, \psi_m^{(2)})_{L^2(\mathcal{V})} = 0 \quad \text{for all } f \in C^\infty_0(\mathcal{V}).
\]
Consequently, $\psi_m^{(2)} = 0$ on $\mathcal{V}$, and thus, by unique continuation, $\psi_m^{(2)} = 0$ globally on $M$, yielding a contradiction.

\textbf{Case 2:} The set
\[
S := \left\{k \in \mathbb{N} \,:\, \mu_m^{(2)} = \mu_{b_k}^{(2)}\right\}
\]
is non-empty but finite. In this case, it follows from \eqref{eigen_expansion} that
\begin{equation}\label{f_m_exp}
\left(f, \psi_m^{(2)}\right)_{L^2(\mathcal{V})} = \sum_{k \in S} \left(f, \psi^{(2)}_{b_k}\right)_{L^2(\mathcal{V})} \left(\psi_{b_k}^{(2)}, \psi_m^{(2)}\right)_{L^2(M)}.
\end{equation}
Since $\{\psi_k^{(2)}\}_{k=1}^\infty$ is a linearly independent set, and the functions $\psi_m^{(2)}$ and $\{\psi_{b_k}^{(2)}\}_{k \in S}$ all lie in the same eigenspace of $-\Delta_g + V_2$, it follows that their restrictions to $\mathcal{V}$ are linearly independent. Therefore, by an argument similar to that in Lemma~\ref{lem_set_F}, see also \cite[Lemma 5.2]{FKU24}, there exists a function $f \in C^{\infty}_0(\mathcal{V})$ such that
\[
(f, \psi_m^{(2)})_{L^2(\mathcal{V})} = 1 \quad \text{and} \quad (f, \psi_{b_k}^{(2)})_{L^2(\mathcal{V})} = 0 \quad \text{for all } k \in S.
\]
This contradicts equation \eqref{f_m_exp}. Thus, we have proved the claim that $b_k = k$ for all $k \in \mathbb{N}$.

Hence, equation \eqref{psi_eq} (recall that $N = 1$ in our setting) can now be rewritten as
\begin{equation}
	\label{psi_eq_alt}
	\mu_k^{(1)} = \mu_k^{(2)} \quad \text{and} \quad \psi_k^{(1)}(x) = \psi_k^{(2)}(x) \quad \text{for all } x \in O,\ k \in \mathbb{N}.
\end{equation}
Moreover, equation \eqref{eigen_expansion} may now be rewritten as follows: for any $f \in C^{\infty}_0(\mathcal{V})$, we have
\begin{equation}
	\label{eigen_expansion_1}
	f(x) = \sum_{k=1}^{\infty} \left(f, \psi^{(2)}_{k} \right)_{L^2(\mathcal{V})} \, \psi^{(2)}_k(x), \qquad x \in M,
\end{equation}

Recall from Remark~\ref{rmk_V_1_normal} that, without loss of generality, we have assumed $\{\psi_{k}^{(1)}\}_{k=1}^{\infty}$ to be an $L^2(M)$-orthonormal basis. Thus, in view of \eqref{eigen_expansion_1}, we are now in exactly the same setting as at the end of the proof of \cite[Theorem 1.11]{FKU24}. Arguing as in \cite[Theorem 1.11]{FKU24}, we conclude that the eigenfunctions $\{\psi_{k}^{(2)}\}_{k=1}^{\infty}$ form an orthonormal basis of $L^2(M)$; see Lemma~\ref{lem_app_orthonormal_basis} for details.

We have therefore arrived at the equality of two orthonormalized spectral data sets, observed on $O$, for the two operators $-\Delta_g + V_j$ on $M$, $j=1,2$. This completes the proof of Theorem~\ref{thm_spectral} in the case $N = 1$, via an application of Proposition~\ref{prop_Gelfand_standard}.
\end{proof}

We now proceed to prove Theorem~\ref{thm_spectral} in the case $N > 1$. To this end, we begin with some preliminary lemmas.

\begin{lemma}
\label{lem_indep}
Let $(M, g)$ be a smooth, closed, and connected Riemannian manifold, and let $\mathcal{V} \subset M$ be a nonempty open set. Let $V \in C^{\infty}(M)$ be real-valued, and let $m \in \mathbb{N}$. Suppose that $\{\phi_k\}_{k=1}^{m}\subset C^\infty(M)$ is a set of linearly independent eigenfunctions of the operator $-\Delta_g + V$ on $M$. Then the set $\{\phi_k|_{\mathcal{V}}\}_{k=1}^{m}$ is also linearly independent.
\end{lemma}

\begin{proof}
We begin with the observation that if a collection of eigenfunctions belongs to the same eigenspace and is linearly independent on $M$, then their restrictions to any nonempty open subset of $M$ are also linearly independent. This follows from the connectedness of $M$ and the unique continuation property for the operator $-\Delta_g + V - \lambda$, where $\lambda$ is the common eigenvalue.

We now consider the case where the eigenfunctions $\{\phi_k\}_{k=1}^m$ do not all belong to the same eigenspace. Let $\lambda_1, \dots, \lambda_L$, with $1 < L \le m$, denote the distinct eigenvalues associated to the collection $\{\phi_k\}_{k=1}^m$. For each $\ell = 1, \dots, L$, let $\{\phi_k^{(\ell)}\}_{k=1}^{d_\ell}$ denote the subset of eigenfunctions corresponding to the eigenvalue $\lambda_\ell$, so that
\[
\{\phi_k\}_{k=1}^m = \bigcup_{\ell=1}^L \{\phi_k^{(\ell)}\}_{k=1}^{d_\ell}, \qquad \sum_{\ell=1}^L d_\ell = m.
\]
Thus, the functions $\phi_k^{(\ell)}$ are simply a reindexing of the original set $\{\phi_k\}_{k=1}^m$, grouped according to eigenvalue. 

Suppose, for contradiction, that the set $\{\phi_k|_{\mathcal{V}}\}_{k=1}^m$ is linearly dependent. Then there exists a nonzero vector $c = (c_k^{(\ell)}) \in \mathbb{C}^m$ such that
\[
\sum_{\ell=1}^L \sum_{k=1}^{d_\ell} c_k^{(\ell)} \phi_k^{(\ell)} = 0 \quad \text{on } \mathcal{V}.
\]
Define
\[
f_\ell := \sum_{k=1}^{d_\ell} c_k^{(\ell)} \phi_k^{(\ell)}, \qquad \ell = 1, \dots, L,
\]
so that the above equation becomes
\begin{equation}\label{eq:sum_f_ell_zero}
\sum_{\ell=1}^L f_\ell = 0 \quad \text{on } \mathcal{V}.
\end{equation}

Since each $\phi_k^{(\ell)}$ is an eigenfunction of $-\Delta_g + V$ with eigenvalue $\lambda_\ell$, we have
\[
(-\Delta_g + V)^j f_\ell = \lambda_\ell^j f_\ell \quad \text{on } M \quad \text{for all } j \ge 0.
\]
Applying the operator $(-\Delta_g + V)^j$ to \eqref{eq:sum_f_ell_zero} for $j = 0, \dots, L - 1$, we obtain the system
\[
\sum_{\ell=1}^L \lambda_\ell^j f_\ell = 0 \quad \text{on } \mathcal{V}, \qquad j = 0, \dots, L - 1.
\]
This is a Vandermonde-type system for the unknowns $f_1, \dots, f_L$, evaluated pointwise. Since the $\lambda_\ell$ are distinct, the corresponding Vandermonde matrix is invertible, and hence
\[
f_\ell = 0 \quad \text{on } \mathcal{V} \quad \text{for all } \ell = 1, \dots, L.
\]

Since each $f_\ell = \sum_{k=1}^{d_\ell} c_k^{(\ell)} \phi_k^{(\ell)}$ is a linear combination of eigenfunctions in the same eigenspace and vanishes on the open set $\mathcal{V}$, it follows from the unique continuation principle (when $d_\ell = 1$) and from the earlier observation (when $d_\ell > 1$) that $c_k^{(\ell)} = 0$ for all $k$ and all $\ell$. This contradicts the assumption that $c \ne 0$.
\end{proof}

\begin{remark}
\label{rem_indep_1}
As an immediate consequence of Lemma~\ref{lem_indep}, we observe that if $\{\phi_k\}_{k=1}^\infty$ is a countable set of eigenfunctions of $-\Delta_g + V$ on $M$ that is linearly independent (in the sense that every finite subset is linearly independent), then the set of their restrictions $\{\phi_k|_{\mathcal{V}}\}_{k=1}^\infty$ is also linearly independent on $\mathcal{V}$.
\end{remark}

\begin{lemma}
\label{lem_V_eq}
Let $V_1, V_2 \in C^\infty(M)$, and let $\mu_k^{(j)}$ and $\psi_k^{(j)}$ (for $j = 1, 2$ and $k \in \mathbb{N}$) be as in Theorem~\ref{thm_spectral}. Suppose that $\mu_N^{(1)} = \mu_{b_N}^{(2)}$ and
\begin{equation}
\label{eq_16_4}
\psi_N^{(1)}(x) = \psi_{b_N}^{(2)}(x), \quad \text{for all } x \in O.
\end{equation}
Then $V_1 = V_2$ on $O$.
\end{lemma}

\begin{proof}
Applying the operator $-\Delta_g$ to both sides of \eqref{eq_16_4}, and using the fact that $\mu_N^{(1)} = \mu_{b_N}^{(2)}$, we obtain
\[
(V_1(x) - V_2(x))\, \psi_N^{(1)}(x) = 0 \quad \text{for all } x \in O.
\]
To conclude the proof, it suffices to show that the set
\[
\mathbb{D} := \{ x \in O : \psi_N^{(1)}(x) \neq 0 \}
\]
is dense in $O$. Suppose, for contradiction, that there exists a nonempty open set $U \subset O$ such that $\mathbb{D} \cap U = \emptyset$. Then $\psi_N^{(1)} = 0$ on $U$. Since
\[
(-\Delta_g + V_1 - \mu_N^{(1)}) \psi_N^{(1)} = 0 \quad \text{on } M,
\]
it follows from the unique continuation principle that $\psi_N^{(1)} \equiv 0$ on $M$, contradicting the assumption that $\psi_N^{(1)}$ is an eigenfunction. Therefore, $\mathbb{D}$ is dense in $O$, and we conclude that $V_1 = V_2$ on $O$.
\end{proof}

\begin{lemma}
\label{lem_indep_2_eigen}
Let $V_1, V_2 \in C^{\infty}(M)$ be real-valued functions. Let $\mathcal{V} \subset M$ be a nonempty open set, and assume that $V_1 = V_2$ on $\mathcal{V}$. Let $\phi_0^{(2)}$ be an eigenfunction of $-\Delta_g + V_2$ associated with the eigenvalue $\lambda_0^{(2)}$. Let $m \in \mathbb{N}$, and let $\{\phi_k^{(1)}\}_{k=1}^m$ be a set of linearly independent eigenfunctions of $-\Delta_g + V_1$ on $M$, associated respectively with eigenvalues $\{\lambda_k^{(1)}\}_{k=1}^m$. If $\lambda_0^{(2)} \notin \{\lambda_1^{(1)}, \ldots, \lambda_m^{(1)}\}$, then
\[
\phi_0^{(2)}|_{\mathcal{V}} \notin \mathrm{Span}\{\phi_1^{(1)}|_{\mathcal{V}}, \ldots, \phi_m^{(1)}|_{\mathcal{V}}\}.
\]
\end{lemma}

\begin{proof}
Suppose for contradiction that there exists a nonzero vector $c = (c_1, \ldots, c_m) \in \mathbb{C}^m$ such that
\begin{equation}
\label{eq_16_5}
\phi_0^{(2)} = \sum_{k=1}^m c_k \, \phi_k^{(1)} \quad \text{on} \quad \mathcal{V}.
\end{equation}
Applying the operator $-\Delta_g + V_2$ to both sides of \eqref{eq_16_5} and using that $V_1 = V_2$ on $\mathcal{V}$, we obtain
\begin{equation}
\label{eq_16_6}
\lambda_0^{(2)} \phi_0^{(2)} = \sum_{k=1}^m \lambda_k^{(1)} c_k \phi_k^{(1)} \quad \text{on} \quad \mathcal{V}.
\end{equation}
Combining equations~\eqref{eq_16_5} and~\eqref{eq_16_6}, we obtain
\begin{equation}
\label{eq_16_7}
\sum_{k=1}^m (\lambda_k^{(1)} - \lambda_0^{(2)}) c_k \phi_k^{(1)} = 0 \quad \text{on } \mathcal{V}.
\end{equation}
Since $\{\phi_k^{(1)}\}_{k=1}^m$ is linearly independent on $M$, Lemma~\ref{lem_indep} implies that the set $\{\phi_k^{(1)}|_{\mathcal{V}}\}_{k=1}^m$ is also linearly independent. Therefore, equation~\eqref{eq_16_7} implies that
\[
(\lambda_k^{(1)} - \lambda_0^{(2)}) c_k = 0 \quad \text{for all } k = 1, \ldots, m.
\]
Since $\lambda_0^{(2)} \notin \{\lambda_1^{(1)}, \ldots, \lambda_m^{(1)}\}$, we conclude that $c_k = 0$ for all $k = 1, \ldots, m$, which contradicts the assumption that $c \ne 0$. Hence, the claim follows.
\end{proof}

\begin{proof}[Proof of Theorem~\ref{thm_spectral} in the case $N>1$]
	We will let the points $p,q\in O$, the small number $\varepsilon\in (0,1)$, and the nonempty connected open sets $\mathcal V,\mathcal W\subset O$ to be defined analogously as in the proof for the case $N=1$. Recall that in particular, condition \eqref{balls_cond} and property \textbf{(P)} are simultaneously satisfied. We will assume without any loss in generality that there holds:
	\begin{equation}
		\label{N_cond}
		\mu_{N-1}^{(1)}<\mu_{N}^{(1)}.
	\end{equation}
	Note that this can always be achieved by increasing, if necessary, the value of $N$ in the statement of the theorem and going to a new subsequence $\{b_k\}_{k=N}^\infty$ to reduce the spectral data yet achieving \eqref{N_cond}. 

Let us now define
\[
F_{\mathcal V}^{(1)} = \mathrm{Span}\,\{\psi^{(1)}_k|_{\mathcal V} : k = 1, \ldots, N-1\}.
\]
Since the sequence $\{\psi_k^{(1)}\}_{k=1}^\infty$ is linearly independent on $\mathcal V$, it follows from Lemma~\ref{lem_indep} that the functions $\{\psi_k^{(1)}|_{\mathcal V}\}_{k=1}^{N-1}$ are also linearly independent in $L^2(\mathcal V)$, and thus $\dim(F_\mathcal{V}^{(1)}) = N - 1$. By the orthogonal decomposition theorem, we may write that 
	$$ L^2(\mathcal V) = F_{\mathcal V}^{(1)} \oplus (F_{\mathcal V}^{(1)})^\perp.$$
	Henceforth, we will identify elements in $L^2(\mathcal V)$ also as elements in $L^2(M)$ by setting them to be zero outside $\mathcal V$. 
	
Next, let us define
\begin{equation}
\label{def_G} 
G_{\mathcal V}^{(1)} := (F_{\mathcal V}^{(1)})^\perp \cap C_0^\infty(\mathcal V).
\end{equation}
Note that $G_{\mathcal V}^{(1)}$ is an infinite-dimensional vector space. Indeed, by Lemma~\ref{lem_set_F}, for any $n \in \mathbb{N}$, there exist functions $\theta_k \in C_0^\infty(\mathcal V)$, $k = 0, \dots, n$, such that
\[
(\theta_k, \psi_{N+l}^{(1)})_{L^2(\mathcal V)} = \delta_{kl}, \quad k, l = 0, \dots, n, \quad \text{and} \quad (\theta_k, \psi_j^{(1)})_{L^2(\mathcal V)} = 0, \quad j = 1, \dots, N-1.
\]
Hence, $\theta_k \in G_{\mathcal V}^{(1)}$ for all $k = 0, \dots, n$, and the family $\{\theta_k\}_{k=0}^n$ is linearly independent. This shows that $G_{\mathcal V}^{(1)}$ is infinite-dimensional.

Let $f \in G_\mathcal{V}^{(1)}$. Recalling from Remark~\ref{rmk_V_1_normal} that the set $\{\psi_k^{(1)}\}_{k=1}^\infty$ may be assumed to form an orthonormal basis for $L^2(M)$, we may write
\begin{equation}\label{f_exp}
	f(x) = \sum_{k=N}^\infty a_k \, \psi_k^{(1)}(x), \quad x \in M,
\end{equation}
where for each $k = N, N+1, \ldots$, the coefficients are given by
\begin{equation}\label{a_k_def_new}
	a_k = \left(f, \psi_k^{(1)}\right)_{L^2(\mathcal{V})} = \int_{\mathcal{V}} f(x) \, \overline{\psi_k^{(1)}(x)} \, dV_g = \int_{\mathcal{V}} f(x) \, \overline{\psi_{b_k}^{(2)}(x)} \, dV_g.
\end{equation}
The convergence in \eqref{f_exp} is to be understood in the $L^2(M)$-sense. Since $f \in H^s(M)$ for all $s \in \mathbb{R}$, it follows that the series in \eqref{f_exp} also converges in the $H^s(M)$-norm for any $s \in \mathbb{R}$. In particular, by Sobolev embedding, the convergence is uniform on $M$.

We may now repeat the argument used in the case $N = 1$ to show that the function
\[
	\tilde{f}(x) = \sum_{k=N}^{\infty} a_k \, \psi^{(2)}_{b_k}(x), \qquad x \in M,
\]
is well defined, smooth, and coincides with $f$ on all of $M$. Thus, for any $f \in G_{\mathcal{V}}^{(1)}$, we have 
\begin{equation}
	\label{f_key_iden}
	f(x) = \sum_{k=N}^\infty \left(f, \psi_{b_k}^{(2)}\right)_{L^2(\mathcal{V})} \, \psi^{(2)}_{b_k}(x),
\end{equation}
where the series converges uniformly on $M$.

Let us now define the set $\mathbb{A} \subset \mathbb{N}$ by
\begin{equation}
\label{def_set_A}
\mathbb{A} = \mathbb{N} \setminus \{b_k\}_{k=N}^\infty.
\end{equation}
We distinguish two cases depending on whether $\mathbb{A}$ is empty or not.

\textbf{Case 1.} $\mathbb{A} = \emptyset$. In this case, the identity \eqref{psi_eq} from the statement of Theorem~\ref{thm_spectral} can be rewritten as
\[
\mu_k^{(2)} = \mu_{c_k}^{(1)} \quad \text{and} \quad \psi_k^{(2)}(x) = \psi_{c_k}^{(1)}(x), \quad x \in O,\quad k=1,2,\dots, 
\]
for some strictly increasing sequence of indices $c_1 < c_2 < \cdots$ in $\mathbb{N}$. But then we are reduced to the case $N = 1$, with the roles of $V_1$ and $V_2$ interchanged. Hence, the argument already given for $N = 1$ applies directly, and the proof of Theorem~\ref{thm_spectral} is complete in this case.

\textbf{Case 2.} $\mathbb{A} \ne \emptyset$. In this case, we claim that  
\begin{equation}
	\label{b_0_claim}
	\forall\, l \in \mathbb A \qquad \mu_l^{(2)} \in \{\mu_1^{(1)}, \ldots, \mu_{N-1}^{(1)}\}.
\end{equation}
We prove this by contradiction. Suppose, on the contrary, that there exists some $l \in \mathbb A$ such that  
\begin{equation}
	\label{contrad_ass_0}
	\mu_l^{(2)} \notin \{\mu_1^{(1)}, \ldots, \mu_{N-1}^{(1)}\}.
\end{equation}
Define  
\[
\mathbb B_l := \{k \in \mathbb N \,:\, k \geq N,\ \mu_{b_k}^{(2)} = \mu_l^{(2)}\}.
\]
Note that $\mathbb B_l$ is either empty or a finite set of positive integers. 

Assume first that $\mathbb B_l = \emptyset$. In view of the assumption \eqref{contrad_ass_0}, Lemma~\ref{lem_indep_2_eigen} implies that the set 
\begin{equation}
\label{set_indep_total}
\{\psi^{(2)}_l|_{\mathcal V}\} \cup  \{\psi_k^{(1)}|_{\mathcal V}\}_{k=1}^{N-1}
\end{equation}
is linearly independent. Therefore, arguing as in the proof of Lemma~\ref{lem_set_F}, it follows that there exists some $\theta \in C^{\infty}_0(\mathcal V)$ such that
\begin{equation}
\label{def_theta}
\left(\theta, \psi_{l}^{(2)}\right)_{L^2(\mathcal V)} = 1 \quad \text{and} \quad \theta \in (F_\mathcal V^{(1)})^\perp.
\end{equation}
Identity \eqref{f_key_iden}, together with the fact that $\mathbb B_l = \emptyset$, implies that
\begin{equation}
\label{f_key_inner_product}
\left(f, \psi_l^{(2)}\right)_{L^2(\mathcal V)} = 0 \qquad \forall\, f \in G_{\mathcal V}^{(1)}.
\end{equation} 
Setting $f = \theta$ in \eqref{f_key_inner_product}, we obtain a contradiction, thus showing that $\mathbb B_l \ne \emptyset$. 

Hence, $\mathbb B_l$ must be a finite set of positive integers. We observe that the set  
\[
\{ \psi_{b_k}^{(2)}|_{\mathcal V} \}_{k \in \mathbb B_l} \cup \{ \psi^{(2)}_l|_{\mathcal V} \}
\]
is linearly independent, and every non-zero element in its span is the restriction to $\mathcal V$ of some eigenfunction of $-\Delta_g + V_2$ on $M$ corresponding to the eigenvalue $\mu_l^{(2)}$. Combining this with the assumption \eqref{contrad_ass_0} and Lemma~\ref{lem_indep_2_eigen}, we deduce that the set  
\[
\{ \psi_{b_k}^{(2)}|_{\mathcal V} \}_{k \in \mathbb B_l} \cup \{ \psi^{(2)}_l|_{\mathcal V} \} \cup \{ \psi_k^{(1)}|_{\mathcal V} \}_{k=1}^{N-1}
\]
is linearly independent. Hence, arguing as in the proof of Lemma~\ref{lem_set_F}, it follows that there exists some $\tilde \theta \in C_0^\infty(\mathcal V)$ such that  
\begin{equation}
\label{def_theta_0}
\left(\tilde \theta, \psi_{b_k}^{(2)}\right)_{L^2(\mathcal V)} = 0 \quad \forall\, k \in \mathbb B_l, \quad \left(\tilde \theta, \psi_l^{(2)}\right)_{L^2(\mathcal V)} = 1, \quad \text{and} \quad \tilde \theta \in (F_\mathcal V^{(1)})^\perp.
\end{equation}

Identity \eqref{f_key_iden} implies that  
\begin{equation}
\label{f_key_inner_product_0}
\left(f, \psi_l^{(2)}\right)_{L^2(\mathcal V)} = \sum_{k \in \mathbb B_l} \left(f, \psi_{b_k}^{(2)}\right)_{L^2(\mathcal V)}\, \left(\psi^{(2)}_{b_k}, \psi_l^{(2)}\right)_{L^2(M)} \qquad \forall\, f \in G_{\mathcal V}^{(1)}.
\end{equation}

Setting $f = \tilde \theta$ in \eqref{f_key_inner_product_0}, we again obtain a contradiction, thereby completing the proof of the claim \eqref{b_0_claim}.

With the proof of \eqref{b_0_claim} complete, we now claim that  
\begin{equation}
	\label{b_claim}
	\forall\, l \in \mathbb A \qquad \psi_l^{(2)}|_{\mathcal V} \in F_{\mathcal V}^{(1)}.
\end{equation}
To prove this claim, assume for contradiction that there exists $l \in \mathbb A$ such that  
\begin{equation}
	\label{contradiction_assumption}
	\psi_l^{(2)}|_{\mathcal V} \notin F_{\mathcal V}^{(1)}.
\end{equation}
Then the set in \eqref{set_indep_total} is linearly independent. Recalling \eqref{N_cond}, and in view of \eqref{b_0_claim} and \eqref{psi_eq}, we observe that the set $\mathbb B_l = \emptyset$. Hence, we are in the same situation as considered earlier when $\mathbb B_l$ was assumed to be empty. Proceeding as after \eqref{set_indep_total}, we obtain a contradiction. This completes the proof of the claim \eqref{b_claim}.

Our next crucial claim is that \eqref{b_0_claim} and \eqref{b_claim} imply the following: for any $l \in \mathbb{A}$, the equality
\begin{equation}\label{second_claim}
\psi^{(2)}_l(x) = \Phi_l^{(1)}(x) \quad \text{for all } x \in O
\end{equation}
holds for some eigenfunction 
\[
(-\Delta_g + V_1)\Phi_l^{(1)} = \mu_{l}^{(2)} \Phi_l^{(1)} \quad \text{on } M.
\]
To prove \eqref{second_claim}, we fix $l \in \mathbb{A}$.  Let $\lambda_1^{(1)}, \dots, \lambda_L^{(1)}$, with $1 \le L \le N-1$, denote the distinct eigenvalues in the set $\{\mu_1^{(1)}, \dots, \mu_{N-1}^{(1)}\}$. For each $m = 1, \dots, L$, let $\{\phi_{k,m}^{(1)}\}_{k=1}^{d_m}$ denote those eigenfunctions among $\{\psi_k^{(1)}\}_{k=1}^{N-1}$ corresponding to the eigenvalue $\lambda_m^{(1)}$, so that
\[
\{\psi_k^{(1)}\}_{k=1}^{N-1} = \bigcup_{m=1}^L \{\phi_{k,m}^{(1)}\}_{k=1}^{d_m}, \qquad \sum_{m=1}^L d_m = N-1.
\]
Thus, the $\phi_{k,m}^{(1)}$ are simply a reindexing of $\{\psi_k^{(1)}\}_{k=1}^{N-1}$, grouped by eigenvalue. In particular, the set $\{\phi_{k,m}^{(1)}|_{\mathcal{V}}\}$ remains linearly independent. 
Since $\psi_l^{(2)}|_{\mathcal{V}} \in F_{\mathcal{V}}^{(1)}$, there exists a nonzero vector $c = (c_{k,m}) \in \mathbb{C}^{N-1}$ such that
\begin{equation}
\label{eq_16_8}
\psi_l^{(2)} = \sum_{m=1}^L \sum_{k=1}^{d_m} c_{k,m} \phi_{k,m}^{(1)} \quad \text{on } \mathcal{V}.
\end{equation}
Applying $-\Delta_g + V_2$ to \eqref{eq_16_8} and using the fact that $V_1 = V_2$ on the set $O$ (see Lemma~\ref{lem_V_eq}), we obtain
\begin{equation}
\label{eq_16_9}
\mu_l^{(2)} \psi_l^{(2)} = \sum_{m=1}^L \sum_{k=1}^{d_m} \lambda_m^{(1)} c_{k,m} \phi_{k,m}^{(1)} \quad \text{on } \mathcal{V}.
\end{equation}
Combining \eqref{eq_16_8} and \eqref{eq_16_9}, we deduce that
\begin{equation}
\label{eq_16_10}
\sum_{m=1}^L \sum_{k=1}^{d_m} (\lambda_m^{(1)} - \mu_l^{(2)}) c_{k,m} \phi_{k,m}^{(1)} = 0 \quad \text{on } \mathcal{V}.
\end{equation}
Since the set $\{\phi_{k,m}^{(1)}|_{\mathcal{V}}\}$ is linearly independent, we conclude from \eqref{eq_16_10} that
\begin{equation}
\label{eq_16_11}
(\lambda_m^{(1)} - \mu_l^{(2)}) c_{k,m} = 0, \quad \text{for all } k = 1, \dots, d_m, \ m = 1, \dots, L.
\end{equation}
Recalling that $\mu_l^{(2)} \in \{\lambda_1^{(1)}, \dots, \lambda_L^{(1)}\}$ and that all $\lambda_m^{(1)}$ are distinct, there exists $m_0$ such that $\mu_l^{(2)} = \lambda_{m_0}^{(1)}$ and $\mu_l^{(2)} \ne \lambda_m^{(1)}$ for all $m \ne m_0$. Then \eqref{eq_16_11} implies that $c_{k,m} = 0$ for all $k = 1, \dots, d_m$ and $m \ne m_0$. Hence, \eqref{eq_16_8} yields
\begin{equation}
\label{eq_16_12}
\psi_l^{(2)} = \Phi_l^{(1)} \quad \text{on } \mathcal{V},
\end{equation}
where
\[
\Phi_l^{(1)} := \sum_{k=1}^{d_{m_0}} c_{k,m_0} \phi_{k,m_0}^{(1)} \quad \text{on } M
\]
is an eigenfunction of $-\Delta_g + V_1$ with eigenvalue $\mu_l^{(2)}$. Since $O$ is connected by assumption (see Theorem~\ref{thm_spectral}), and $V_1 = V_2$ on $O$, it follows from \eqref{eq_16_12} and the unique continuation property for solutions to the equation $(-\Delta_g + V_1 - \mu_l^{(2)}) u = 0$ on $O$ that the equality \eqref{eq_16_12} extends to all of $O$. This completes the proof of \eqref{second_claim}.

In view of \eqref{second_claim} together with \eqref{psi_eq}, we have shown that for all $k \in \mathbb{N}$,
\begin{equation}
\label{eq_16_12_1}
\mu_k^{(2)} = \mu_{c_k}^{(1)}, \quad \psi_k^{(2)}(x) = \Phi_{c_k}^{(1)}(x) \quad \text{for all } x \in O,
\end{equation}
where the sequence $(c_k)_{k=1}^\infty$ of positive integers can be chosen to be strictly increasing, and each $\Phi_{c_k}^{(1)}$ is an eigenfunction of $-\Delta_g + V_1$ corresponding to the eigenvalue $\mu_{c_k}^{(1)}$.  Since the set of eigenfunctions $\{\psi_k^{(2)}\}_{k=1}^\infty$ is linearly independent, it follows from \eqref{eq_16_12_1}, Remark~\ref{rem_indep_1}, and Lemma~\ref{lem_indep} that the set $\{\Phi_{c_k}^{(1)}\}_{k=1}^\infty$ is linearly independent on $O$, and hence also on $M$. In view of Remark~\ref{rem_extension_basis}, we have thus reduced the situation to one analogous to Theorem~\ref{thm_spectral}, but with $N = 1$ and the roles of $V_1$ and $V_2$ reversed. We can now complete the proof of Theorem~\ref{thm_spectral} by applying the same argument as in the case $N = 1$.
\end{proof}

\section{Proof of Theorem~\ref{thm_anosov}}

\label{proof_thm_2}

We begin with the following observation.
\begin{remark}
\label{rmk_proof_anosov}
When $\mathbb{N} \setminus \{a_k\}_{k=1}^\infty$ is either empty or finite, Theorem~\ref{thm_anosov} follows directly by applying the same arguments as in the proof of Theorem~\ref{thm_spectral}. Analogously, Theorem~\ref{thm_anosov} also holds when $\mathbb{N} \setminus \{b_k\}_{k=1}^\infty$ is either empty or finite. Therefore, for the remainder of this section, we assume that both sets $\mathbb{N} \setminus \{a_k\}_{k=1}^\infty$ and $\mathbb{N} \setminus \{b_k\}_{k=1}^\infty$ are infinite.
\end{remark}

Let us now recall a deep Paley--Wiener type interpolation result due to Beurling~\cite{BC89} and Kahane~\cite{K57}, stated here in a slightly weaker form (see also~\cite[Section~1.3]{OU09}). The relevant terminology is introduced in Definition~\ref{def_u_d}.
\begin{theorem}
\label{thm_PW_interpolation}
Let $\Gamma = \{\gamma_k\}_{k=1}^{\infty} \subset \mathbb{R}$ be a uniformly discrete set with zero upper uniform density, i.e., $D^+(\Gamma) = 0$. Let $J \subset \mathbb{R}$ be a nonempty, open, bounded interval. Then, for any sequence $c = (c_1, c_2, \ldots) \in \ell^{2}(\mathbb{N})$, there exists a function $h \in L^2(\mathbb{R})$ with $\operatorname{supp}(h) \subset J$ such that
\[
\mathcal{F}(h)(\gamma_k) = c_k \quad \text{for all } k \in \mathbb{N},
\]
and
\[
\|h\|_{L^2(\mathbb{R})} \leq C \|c\|_{\ell^2(\mathbb{N})},
\]
for some constant $C > 0$ independent of $c$.
\end{theorem}
In Theorem~\ref{thm_PW_interpolation}, we denote by $\mathcal{F}$ the Fourier transform defined for $f \in L^1(\mathbb{R})$ by
\[
\mathcal{F}(f)(\xi) = \int_{\mathbb{R}} e^{-i x \xi} f(x)\, dx,
\]
which extends to a bounded linear operator from $L^2(\mathbb{R})$ to $L^2(\mathbb{R})$. Since $h \in L^2(\mathbb{R})$ has compact support, it follows that $\mathcal{F}(h) \in C^\infty(\mathbb{R})$.

The interpolation result of Theorem~\ref{thm_PW_interpolation} plays a key role in the proofs of both the following proposition and Lemma~\ref{lem_ortho_nonempty} below.

\begin{proposition}
\label{prop_injectivity_anosov}
Let $(M, g)$ be a smooth, closed, connected Riemannian manifold of dimension $n \geq 2$, and let $U \subset M$ be a nonempty, open subset with smooth boundary such that $M\setminus\overline{U}\ne\emptyset$. Let $V \in C^{\infty}(M)$ be a real-valued potential, and suppose that assumption~\textbf{(UO)} holds for the pair $(U, V)$.
 Let $\{\phi_k\}_{k=1}^\infty \subset C^{\infty}(M)$ be an $L^2(M)$-orthonormal basis consisting of eigenfunctions:
\[
(-\Delta_g + V)\phi_k = \lambda_k \phi_k \quad \text{on } M,
\]
with $\lambda_1 \leq \lambda_2 \leq \cdots$.

Let $D = \{d_k\}_{k=1}^\infty \subset \mathbb{N}$ be a nonempty set with $d_1 < d_2 < \cdots$, and suppose that the set $\{\lambda_{d_k}\}_{k=1}^\infty$ is $\{\lambda_k\}_{k=1}^\infty$-sparse in the sense of Definition~\ref{def_adm}. Given any $a = (a_1, a_2, \ldots) \in \ell^2(\mathbb{N})$, define the linear map
\begin{equation}
\label{T_prop_1}
\mathcal{T}(a) := \left( \sum_{k=1}^\infty a_k\, \phi_{d_k} \right)\bigg|_U,
\end{equation}
where the series converges in $L^2(M)$, and thus $\mathcal{T}(a)\in L^2(U)$. Then the mapping $\mathcal{T} : \ell^2(\mathbb{N}) \to L^2(U)$ is continuous, injective, and has closed range. That is, there exist constants $C, C' > 0$ such that
\begin{equation}
\label{T_map}
C\, \|a\|_{\ell^2(\N)} \leq \|\mathcal{T}(a)\|_{L^2(U)} \leq C'\, \|a\|_{\ell^2(\N)} \quad \text{for all } a \in \ell^2(\mathbb{N}).
\end{equation}
\end{proposition}

\begin{proof}
The fact that the series in \eqref{T_prop_1} converges in $L^2(M)$ and that the linear map $\mathcal{T} : \ell^2(\mathbb{N}) \to L^2(U)$ is bounded are both straightforward consequences of the orthonormality of the basis $\{\phi_k\}_{k=1}^\infty$ in $L^2(M)$.

The remainder of the proof is devoted to establishing the lower bound in~\eqref{T_map}. To that end, it suffices to consider $a \in \ell^2(\mathbb{N}) \setminus \{0\}$, since the inequality is trivially satisfied when $a = 0$. We begin by noting that the operator
\[
P = -\Delta_{g} + (V(x) - \lambda_1 + 1)
\]
satisfies $P \ge I$ in the sense of self-adjoint operators. Let $f \in L^2(M)$, and consider the following initial value problem:
\begin{equation}
\label{wave_prop_anosov}
\begin{aligned}
\begin{cases}
\partial_t^2 v - \Delta_{g} v + (V(x) - \lambda_1 + 1)\, v = 0 
& \text{in } \mathcal{D'}(\mathbb{R} \times M), \\
v(0,\cdot) = f & \text{on } M, \\
\partial_t v(0,\cdot) = 0 & \text{on } M.
\end{cases}
\end{aligned}
\end{equation}
The problem \eqref{wave_prop_anosov} admits a unique solution in the sense of distributions,
\[
v \in C^0(\mathbb{R}; L^2(M)) \cap C^1(\mathbb{R}; H^{-1}(M)).
\]
and for any $\tilde T > 0$, there exists a constant $C = C(\tilde T) > 0$ such that
\begin{equation}
\label{eq_17_1}
\max_{0 \le t \le \tilde T} \left( \|v(t,\cdot)\|_{L^2(M)} + \|\partial_t v(t,\cdot)\|_{H^{-1}(M)} \right) \le C \|f\|_{L^2(M)},
\end{equation}
see, for example, \cite[Corollary 2.36]{KKL}.

Given any $a \in \ell^2(\mathbb{N}) \setminus \{0\}$, let
\begin{equation}
\label{wave_prop_anosov_spaces}
v_a \in C^0(\mathbb{R}; L^2(M)) \cap C^1(\mathbb{R}; H^{-1}(M))
\end{equation}
denote the unique weak solution to \eqref{wave_prop_anosov} with initial data
\[
f = \sum_{k=1}^\infty a_k \, \phi_{d_k} \in L^2(M).
\]
Then, the function $v_a$ admits the expansion
\[
v_a(t,\cdot) = \sum_{k=1}^\infty a_k \, \cos\left( \sqrt{\lambda_{d_k} - \lambda_1 + 1}\, t \right) \, \phi_{d_k},
\]
where the series converges in the topology of \eqref{wave_prop_anosov_spaces}. Note that $\mathcal{T}(a) = v_a(0,\cdot)|_U$.
Fix a point $x_0 \in U$, and consider a sufficiently small open connected geodesic ball $\widetilde{U} = \mathbb{B}_r(x_0)$ with smooth boundary such that $\widetilde{U} \Subset U$. Define
\begin{equation}
\label{def_prop_time}
T = \frac{1}{2} \,\mathrm{dist}_g(\partial \widetilde{U}, M \setminus U) > 0.
\end{equation}

We next claim that
\begin{equation}
\label{energy_est}
\|v_a\|_{L^{\infty}(0,T; L^2(\widetilde{U}))} \leq C_0\, \|\mathcal{T}(a)\|_{L^2(U)},
\end{equation}
where $T$ is as defined in \eqref{def_prop_time}, and $C_0 > 0$ is a constant independent of $a$.
To show this, define
\begin{equation}
\label{eq_17_2}
\tilde f = 
\begin{cases} 
\mathcal T(a), & \text{in } U,\\
0, & \text{in } M \setminus U,
\end{cases}
\end{equation}
and let $w_a \in C^0([0,\infty); L^2(M)) \cap C^1([0,\infty); H^{-1}(M))$ be the unique solution to \eqref{wave_prop_anosov} with initial data $f = \tilde f$. Then the difference $v_a - w_a$ satisfies \eqref{wave_prop_anosov} with initial data
\[
f = 
\begin{cases} 
0, & \text{in } U,\\
\sum\limits_{k=1}^\infty a_k\, \phi_{d_k}, & \text{in } M \setminus U.
\end{cases}
\]

By the finite propagation speed property (with $T$ as in \eqref{def_prop_time}), we obtain
\begin{equation}
\label{eq_17_3}
v_a(t,x) = w_a(t,x), \quad \text{for all } t \in [0, T] \text{ and for a.e. } x \in \widetilde{U},
\end{equation}
see \cite[Chapter 6, Proposition 1.3 and Lemma 1.2, p. 487]{Taylor_book_I} and \cite[Lemma 2.15]{Saksala_Shedlock_2025}.

Using \eqref{eq_17_3} together with the energy estimate \eqref{eq_17_1} for $w_a$ (with $\tilde T=T$ as in \eqref{def_prop_time} and $\tilde f$ as in \eqref{eq_17_2}), we obtain
\[
\max_{0 \le t \le T} \|v_a(t,\cdot)\|_{L^2(\widetilde{U})} = \max_{0 \le t \le T} \|w_a(t,\cdot)\|_{L^2(\widetilde{U})} \le C \|\tilde f\|_{L^2(M)} = C \|\mathcal T(a)\|_{L^2(U)},
\]
which establishes \eqref{energy_est}.

Next, we define the set
\begin{equation}
\label{ef_Gamma_prop}
\Gamma = \left\{
\begin{array}{l}
\sqrt{\lambda_{d_j} - \lambda_1 + 1} - \sqrt{\lambda_{d_k} - \lambda_1 + 1}, \\
\pm\left( \sqrt{\lambda_{d_j} - \lambda_1 + 1} + \sqrt{\lambda_{d_k} - \lambda_1 + 1} \right)
\end{array}
: j, k \in \mathbb{N} \right\}.
\end{equation}
Note that $0 \in \Gamma$. Since $\{\lambda_{d_k}\}_{k=1}^\infty$ is $\{\lambda_k\}_{k=1}^\infty$-sparse, the set $\Gamma$ is uniformly discrete and has zero upper uniform density; see Definition~\ref{def_adm}.

Therefore, by Theorem~\ref{thm_PW_interpolation} with $\Gamma$ as in~\eqref{ef_Gamma_prop}, there exists a function
\begin{equation}
\label{h_prop_1}
h \in L^2(\mathbb{R}) \quad \text{with} \quad \mathrm{supp}\, h \subset (0, T),
\end{equation}
such that
\begin{equation}
\label{interpolation_eq}
\begin{cases}
\mathcal{F}(h)(\xi) = 0 & \text{for all } \xi \in \Gamma \setminus \{0\}, \\
\mathcal{F}(h)(0) = 1,
\end{cases}
\end{equation}
and
\begin{equation}
\label{eq_17_4}
\|h\|_{L^2(\mathbb{R})} \le C.
\end{equation}
Using \eqref{h_prop_1} and \eqref{eq_17_4}, we obtain
\begin{equation}
\label{h_norm}
\|h\|_{L^1(\mathbb{R})} \le \sqrt{T}\, \|h\|_{L^2(\mathbb{R})} \le C,
\end{equation}
where $C = C(T) > 0$.

Let us now consider
\begin{equation}
\label{def_Integral}
I = \int_\mathbb{R} \int_{\widetilde{U}} h(t)\, |v_a(t,x)|^2 \, dV_{g}(x)\, dt.
\end{equation}
In view of \eqref{h_prop_1}, \eqref{energy_est}, and \eqref{h_norm}, we have
\begin{equation}
\label{I_est_1}
|I| \leq C_1\, \|\mathcal{T}(a)\|^2_{L^2(U)},
\end{equation}
where $C_1 > 0$ is independent of $a$.

We now proceed to obtain a lower bound for $|I|$. To that end, we fix $\epsilon > 0$ and choose $N_0 \in \mathbb{N}$ sufficiently large so that
\begin{equation}
\label{remainder_est_0}
\sum_{k = N + 1}^{\infty} |a_k|^2 \leq \epsilon^2 \qquad \text{for all } N > N_0.
\end{equation}
Note that \eqref{remainder_est_0}, together with the fact that $\{\phi_k\}_{k=1}^\infty$ is an orthonormal basis of $L^2(M)$, implies that
\begin{equation}
\label{remainder_est_1}
\left\| \sum_{k = N+1}^\infty a_k\, \cos\left( \sqrt{\lambda_{d_k} - \lambda_1 + 1}\, t \right)\, \phi_{d_k} \right\|_{L^\infty(0,T; L^2(M))} \leq \epsilon \quad \text{for all } N > N_0.
\end{equation}

Returning to \eqref{def_Integral}, we now estimate the contribution of the truncated sum. Using \eqref{remainder_est_1}, together with \eqref{I_est_1} and \eqref{h_norm}, we obtain
\begin{equation}
\label{I_est_2}
\begin{aligned}
&\left| \int_0^T \int_{\widetilde{U}} h(t) \left| \sum_{k=1}^N a_k\, \cos\left( \sqrt{\lambda_{d_k} - \lambda_1 + 1}\, t \right)\, \phi_{d_k}(x) \right|^2 dV_g(x)\, dt \right| \\
&\le \int_0^T \int_{\widetilde{U}} |h(t)| \left( |v_a(t,x)| + \left| \sum_{k = N+1}^\infty a_k\, \cos\left( \sqrt{\lambda_{d_k} - \lambda_1 + 1}\, t \right)\, \phi_{d_k}(x) \right| \right)^2 dV_g(x)\, dt \\
&\le 2 C_1\, \|\mathcal{T}(a)\|^2_{L^2(U)} + 2 \int_0^T |h(t)| \left\| \sum_{k = N+1}^\infty a_k\, \cos\left( \sqrt{\lambda_{d_k} - \lambda_1 + 1}\, t \right)\, \phi_{d_k} \right\|^2_{L^2(M)} dt \\
&\le 2 C_1\, \|\mathcal{T}(a)\|^2_{L^2(U)} + 2 \epsilon^2 C \quad \text{for all } N > N_0,
\end{aligned}
\end{equation}
where $C_1 > 0$ and $C > 0$ are constants independent of $a$.

Let us now give a lower bound for the left-hand side of the above estimate. Denote the distinct elements of the set $\{\lambda_{d_k}\}_{k=1}^N$ by $\mu_{1}<\ldots<\mu_m$, written in strictly increasing order for some $m=m(N) \leq N$, and define for each $k=1,\ldots,m$ the nonempty finite set
\[
\mathbb{D}_k = \{l \in \{1,2,\dots, N\} \,:\, \lambda_{d_l} = \mu_k\}.
\]

Using the identity
\begin{equation}
\label{eq_cos_identity_new}
\cos\left(\sqrt{\lambda_{d_k} - \lambda_1 + 1}\, t\right)
= \frac{1}{2} \left( e^{\mathrm{i} \sqrt{\lambda_{d_k} - \lambda_1 + 1}\, t}
+ e^{-\mathrm{i} \sqrt{\lambda_{d_k} - \lambda_1 + 1}\, t} \right), \quad k \in \mathbb{N},
\end{equation}
together with the interpolation property \eqref{interpolation_eq}, we deduce that for all $N > N_0$,
\begin{align}
\label{I_est_3}
&\left| \int_0^T \int_{\widetilde{U}} h(t) \left| \sum_{k=1}^N a_k\, \cos\left(\sqrt{\lambda_{d_k} - \lambda_1 + 1}\, t\right)\, \phi_{d_k}(x) \right|^2\, dV_g(x)\, dt \right| \notag \\
&= \frac{1}{4} \Bigg| \int_{\widetilde{U}} 
\sum_{k=1}^N \sum_{j=1}^N a_k \overline{a_j} \bigg( \int_{\mathbb{R}} h(t)\, \Big[
\notag \\
&\qquad e^{-\mathrm{i}(\sqrt{\lambda_{d_j} - \lambda_1 + 1} - \sqrt{\lambda_{d_k} - \lambda_1 + 1}) t}
+ e^{-\mathrm{i}(\sqrt{\lambda_{d_j} - \lambda_1 + 1} + \sqrt{\lambda_{d_k} - \lambda_1 + 1}) t}
\notag \\
&\qquad + e^{-\mathrm{i}(-\sqrt{\lambda_{d_j} - \lambda_1 + 1} - \sqrt{\lambda_{d_k} - \lambda_1 + 1}) t}
+ e^{-\mathrm{i}(\sqrt{\lambda_{d_k} - \lambda_1 + 1} - \sqrt{\lambda_{d_j} - \lambda_1 + 1}) t}
\Big] dt \bigg)
\notag \\
&\qquad \times\, \phi_{d_k}(x)\, \overline{\phi_{d_j}(x)}\, dV_g(x) \Bigg| \notag \\
&= \frac{1}{2} \left| \int_{\widetilde{U}} 
\sum_{\ell=1}^{m} \sum_{k, j \in \mathbb{D}_\ell}
a_k \overline{a_j}\, \phi_{d_k}(x)\, \overline{\phi_{d_j}(x)}\, dV_g(x) \right|  = \frac{1}{2} \sum_{\ell=1}^{m(N)} \int_{\widetilde{U}} \left| \sum_{k \in \mathbb{D}_\ell} a_k\, \phi_{d_k}(x) \right|^2\, dV_g(x).
\end{align}

Letting $l = 1, \ldots, m(N)$, and defining $w_l \in C^{\infty}(M)$ by
\[
w_l(x) = \sum_{k \in \mathbb{D}_l} a_k\, \phi_{d_k}(x), \quad x \in M,
\]
we have
\[
(-\Delta_g + V)\, w_l = \mu_l\, w_l \quad \text{on } M.
\]
If $w_l \ne 0$, then applying the observability estimate~\eqref{eq_observability_eigenfunctions_def}, as ensured by assumption~\textbf{(UO)}, yields
\begin{equation}
\label{eq_17_5}
\left\| \sum_{k \in \mathbb{D}_l} a_k\, \phi_{d_k} \right\|_{L^2(\widetilde{U})}
\geq C_2 \left\| \sum_{k \in \mathbb{D}_l} a_k\, \phi_{d_k} \right\|_{L^2(M)}
= C_2 \left( \sum_{k \in \mathbb{D}_l} |a_k|^2 \right)^{1/2},
\end{equation}
where $C_2 > 0$ is a constant independent of $a$. Here, we used the fact that $\{\phi_k\}_{k=1}^\infty$ forms an orthonormal basis of $L^2(M)$ to obtain the final equality. The bound~\eqref{eq_17_5} holds trivially in the case $w_l = 0$.

Using the lower bound \eqref{eq_17_5} in \eqref{I_est_3} and combining it with \eqref{I_est_2}, we obtain
\[
\frac{(C_2)^2}{2} \left( \sum_{k=1}^N |a_k|^2 \right) \leq 2 C_1\, \|\mathcal{T}(a)\|^2_{L^2(U)} + 2 \epsilon^2 C,
\]
for all $N > N_0$. Letting $N \to \infty$ and then $\epsilon \to 0$, we conclude that
\[
\frac{(C_2)^2}{2} \|a\|_{\ell^2}^2 \leq 2 C_1\, \|\mathcal{T}(a)\|^2_{L^2(U)},
\]
which establishes the desired lower bound in~\eqref{T_map}.
\end{proof}

The following corollary of Proposition~\ref{prop_injectivity_anosov} establishes an interpolation result that may be of independent interest.
\begin{corollary}
Let the notation be as in Proposition~\ref{prop_injectivity_anosov}. Then for any sequence $c = (c_1, c_2, \dots) \in \ell^2(\mathbb{N})$, there exists a function $f \in L^2(U)$ such that
\[
(f, \phi_{d_k})_{L^2(U)} = c_k \quad \text{for all } k \in \mathbb{N}.
\]
\end{corollary}
\begin{proof}
The adjoint $\mathcal{T}^* : L^2(U) \to \ell^2(\mathbb{N})$ of the operator $\mathcal{T}$ defined in~\eqref{T_prop_1} is given by
\[
\mathcal{T}^*(f) = \big( (f, \phi_{d_1})_{L^2(\mathcal{V})},\, (f, \phi_{d_2})_{L^2(U)},\, \dots \big).
\]
By Proposition~\ref{prop_injectivity_anosov}, the operator $\mathcal{T}^*$ is surjective, and the result follows.
\end{proof}

\begin{lemma}
	\label{lem_ortho_nonempty}
Let us adopt the same notations as in Proposition~\ref{prop_injectivity_anosov}. Let $\widetilde{V} \in C^{\infty}(M)$ be real-valued, and assume that $\widetilde{V}|_{U} = V|_{U}$. Suppose that $\Phi$ is an eigenfunction of the operator $-\Delta_g + \widetilde{V}$ on $M$ corresponding to an eigenvalue $\lambda$, and assume that
\begin{equation}
	\label{Phi_cond}
	\Phi|_{U} = \mathcal{T}(a) \quad \text{for some } a \in \ell^2(\mathbb{N}).
\end{equation}
Then $\lambda \in \{ \lambda_{d_k} \}_{k=1}^\infty$, and there exists an eigenfunction $\phi$ of $-\Delta_g + V$ on $M$, associated to the eigenvalue $\lambda$, such that
\begin{equation}
\label{claim_lem_ortho_nonempty}
\Phi|_{U} = \phi|_{U}.
\end{equation}
\end{lemma}

\begin{proof}
Let $a \in \ell^2(\mathbb{N})$ be such that \eqref{Phi_cond} holds, and recall that
\[
\mathcal{T}(a) = \left(\sum_{k=1}^\infty a_k \phi_{d_k}\right)\big|_{U}.
\]
Recalling that the series $\sum_{k=1}^\infty a_k \phi_{d_k}$ converges in $L^2(M)$, we have
\[
(-\Delta_g + V)\left(\sum_{k=1}^\infty a_k \phi_{d_k}\right) = \sum_{k=1}^\infty a_k \lambda_{d_k} \phi_{d_k} \quad \text{in the } H^{-2}(M) \text{ sense}.
\]
Restricting both sides to the set $U$ and using \eqref{Phi_cond} along with the identity $V = \widetilde{V}$ on $U$, we deduce that
\begin{equation}\label{T_a_eq}
\left(\sum_{k=1}^\infty a_k\,(\lambda_{d_k} - \lambda)\, \phi_{d_k} \right)\big|_{U} = 0 \quad \text{in the $H^{-2}(U)$ sense}.
\end{equation}

Define $f$ by
\[
f = \sum_{k=1}^\infty a_k\,(\lambda_{d_k} - \lambda)\, \phi_{d_k},
\]
where the series converges in $H^{-2}(M)$, so in particular $f \in H^{-2}(M)$. Then \eqref{T_a_eq} implies
\begin{equation}
\label{T_a_eq_f_new}
f\big|_{U} = 0 \quad \text{in the $H^{-2}(U)$ sense}.
\end{equation}

Define
\begin{equation}
\label{eq_19_1}
w^f \in C^0(\mathbb{R}; H^{-2}(M)) \cap C^1(\mathbb{R}; H^{-3}(M))
\end{equation}
by
\[
w^f(t, \cdot) := \sum_{k=1}^\infty a_k\,(\lambda_{d_k} - \lambda)\, \cos\left(\sqrt{\lambda_{d_k} - \lambda_1 + 1}\, t\right)\, \phi_{d_k},
\]
where the series is convergent in the topology of the space in \eqref{eq_19_1}. It is straightforward to verify that \( w^f \) satisfies the following initial value problem in the sense of distributions:
\begin{equation}
\label{wave_prop_anosov_w}
\begin{cases}
\partial_t^2 w - \Delta_g w + (V(x) - \lambda_1 + 1)\, w = 0 
&\text{in } \mathcal{D}'(\mathbb{R} \times M), \\
w(0,x) = f &\text{in } H^{-2}(M), \\
\partial_t w(0,x) = 0 &\text{in } H^{-3}(M).
\end{cases}
\end{equation}
Moreover, the problem \eqref{wave_prop_anosov_w} admits a unique solution in the space \eqref{eq_19_1}. To see this, suppose that $w_1$ and $w_2$ are two solutions to \eqref{wave_prop_anosov_w} with the same initial data $f \in H^{-2}(M)$. Then their difference $\tilde{w} := w_1 - w_2$ satisfies \eqref{wave_prop_anosov_w} with homogeneous initial data $f = 0$. By \cite[Theorem 26.1.4]{Hor_book_IV}, it follows that $\tilde{w} \in C^\infty(\mathbb{R} \times M)$. Since the only smooth solution of \eqref{wave_prop_anosov_w} with vanishing initial data is identically zero, see \cite[Theorem 2.13]{Saksala_Shedlock_2025}, we conclude that $\tilde{w} = 0$, and hence $w_1 = w_2$, see also \cite[Lemma 6.3]{AFO_2022}.   

Fix a point $x_0 \in U$, and let $\widetilde{U} = \mathbb{B}_r(x_0)$ be a sufficiently small, open, connected geodesic ball with smooth boundary such that $\widetilde{U} \Subset U$.  Let $T > 0$ be defined by \eqref{def_prop_time}. In view of \eqref{T_a_eq_f_new}, and by the finite speed of propagation for wave equations, we obtain that for each $t \in [0,T]$,
\begin{equation}
	\label{w_vanish_1}
	w^f(t,\cdot)\big|_{\widetilde{U}} 
	= \left( \sum_{k=1}^\infty a_k (\lambda_{d_k} - \lambda) \cos\left( \sqrt{\lambda_{d_k} - \lambda_1 + 1}\, t \right) \phi_{d_k} \right)\bigg|_{\widetilde{U}} = 0 \quad \text{in }\mathcal{D}' (\widetilde{U}),
\end{equation}
see \cite[Chapter 6, Proposition 1.3, p.~487]{Taylor_book_I} and \cite[Lemma 2.15]{Saksala_Shedlock_2025}.

Let us denote by $\mu_1 < \mu_2 < \ldots$ the distinct elements of the sequence $\{\lambda_{d_k}\}_{k=1}^\infty$, and define for each $l \in \mathbb{N}$ the nonempty finite set
\[
\mathbb{D}_l := \{k \in \mathbb{N} \,:\, \lambda_{d_k} = \mu_l\}.
\]
We also define
\[
\tilde \Gamma := \left\{ \pm \sqrt{\lambda_{d_k} - \lambda_1 + 1} : \, k \in \mathbb{N} \right\}.
\]
We claim that $\tilde\Gamma$ is uniformly discrete and satisfies $D^+(\tilde\Gamma) = 0$; see Definition~\ref{def_u_d}. To prove this, we first note that since the sequence $\{\lambda_{d_k}\}_{k=1}^\infty$ is $\{\lambda_k\}_{k=1}^\infty$-sparse, the set $\Gamma$ defined in~\eqref{ef_Gamma_prop} is uniformly discrete and has zero upper uniform density (see Definition~\ref{def_adm}).  Now, observe that for any $\tilde\gamma_j, \tilde\gamma_k \in \tilde\Gamma$, the pairwise difference $\tilde\gamma_j - \tilde\gamma_k =: \gamma_{jk} \in \Gamma$, and since $0 \in \Gamma$, the uniform discreteness of $\Gamma$ implies that the distances $|\tilde\gamma_j - \tilde\gamma_k| = |\gamma_{jk} - 0|$ are uniformly bounded below for $j \ne k$. Hence, $\tilde\Gamma$ is uniformly discrete. 

To show that $D^+(\tilde\Gamma) = 0$, consider the translated set
\[
\tilde\Gamma' := \{\pm \sqrt{\lambda_{d_k} - \lambda_1 + 1} +  \sqrt{\lambda_{d_1} - \lambda_1 + 1} : \, k \in \mathbb{N}\}.
\]
It is clear that $\tilde\Gamma' \subset \Gamma$, so by Definition~\ref{def_u_d}, we have
\[
0 \le D^+(\tilde\Gamma) = D^+(\tilde\Gamma') \le D^+(\Gamma) = 0,
\]
which completes the proof of the claim.

Therefore, the set $\tilde\Gamma$ satisfies the hypotheses of Theorem~\ref{thm_PW_interpolation}, and by applying that theorem, we conclude that for any $l \in \mathbb{N}$, there exists a function
\begin{equation}
\label{eq_17_6}
h_l \in L^2(\mathbb{R}) \quad \text{with} \quad \operatorname{supp}(h_l) \subset (0, T),
\end{equation}
such that
\begin{equation}
\label{h_eq_interpolation}
\begin{cases}
\mathcal{F}(h_l)\big(-\sqrt{\mu_k - \lambda_1 + 1}\big) = 0 & \text{for all } k \in \mathbb{N}, \\
\mathcal{F}(h_l)\big(\sqrt{\mu_k - \lambda_1 + 1}\big) = 0 & \text{for all } k \in \mathbb{N} \setminus \{l\}, \\
\mathcal{F}(h_l)\big(\sqrt{\mu_l - \lambda_1 + 1}\big) = 1. &
\end{cases}
\end{equation}

Letting $\chi \in C_0^\infty(\widetilde U)$, we obtain from \eqref{w_vanish_1} that for all $t \in [0,T]$,
\begin{equation}
\label{eq_17_7}
\begin{aligned}
0 &= \left\langle \sum_{k=1}^\infty a_k (\lambda_{d_k} - \lambda) \cos\left( \sqrt{\lambda_{d_k} - \lambda_1 + 1}\, t \right) \phi_{d_k},\, \chi \right\rangle_{\mathcal{D}'(M), C^\infty(M)} \\
&= \sum_{k=1}^\infty \cos\left( \sqrt{\lambda_{d_k} - \lambda_1 + 1}\, t \right) a_k (\lambda_{d_k} - \lambda) (\phi_{d_k}, \chi)_{L^2(M)},
\end{aligned}
\end{equation}
where $\langle \cdot, \cdot \rangle_{\mathcal{D}'(M), C^\infty(M)}$ denotes the dual pairing between distributions and test functions, induced by the $L^2$ inner product on $M$ with respect to the Riemannian volume form $dV_g$.

Multiplying \eqref{eq_17_7} by $h_l(t)$, integrating over $t \in \mathbb{R}$, and using \eqref{eq_17_6}, we obtain
\begin{equation}
\label{eq_17_8}
\sum_{k=1}^\infty \int_{\mathbb{R}} h_l(t) \cos\left( \sqrt{\lambda_{d_k} - \lambda_1 + 1}\, t \right) a_k (\lambda_{d_k} - \lambda) (\phi_{d_k}, \chi)_{L^2(M)}\, dt = 0,
\end{equation}
for all $l = 1, 2, \dots$. Here, the exchange of the sum and the integral is justified by the fact that $h_l \in L^1(\mathbb{R})$, the bound $|(\phi_{d_k}, \chi)_{L^2(M)}| \le C (\lambda_{d_k} - \lambda_1 + 1)^{-L}$ for all $L > 0$, as well as the following consequence of Weyl's law for the operator $-\Delta_g + V - \lambda_1 + 1$:
\[
\lambda_k - \lambda_1 + 1 \ge C k^{2/n}, \quad \text{for } k \gg 1,
\]
for some constant $C > 0$. In addition, we use that $d_k \ge k$ and that $a \in \ell^2(\mathbb{N})$.

Using \eqref{eq_cos_identity_new} and \eqref{h_eq_interpolation}, we obtain from \eqref{eq_17_8} that
\[
\sum_{k\in \mathbb{D}_l} a_k (\mu_l - \lambda) (\phi_{d_k}, \chi)_{L^2(\widetilde U)} = 0,
\]
for all $l = 1, 2, \dots$. As $\chi \in C^\infty_0(\widetilde U)$ is arbitrary, it follows that
\begin{equation}
\label{eq_17_9}
(\mu_l - \lambda) \sum_{k \in \mathbb{D}_l} a_k \phi_{d_k}(x) = 0, \quad x \in \widetilde U,
\end{equation}
for all $l = 1, 2, \dots$.

Now, if $\lambda \ne \mu_l$ for any $l = 1, 2, \dots$, then \eqref{eq_17_9} implies that $\sum_{k \in \mathbb{D}_l} a_k \phi_{d_k} = 0$ on $\widetilde U$ for all $l = 1, 2, \dots$. Since $\{\phi_k\}_{k=1}^\infty$ are linearly independent, it follows that $a_k = 0$ for all $k \in \mathbb{N}$, which contradicts \eqref{Phi_cond}. Hence, there exists a unique $l_0$ such that $\lambda = \mu_{l_0}$, and it follows from \eqref{eq_17_9} that $a_k = 0$ for all $k \in \mathbb{N} \setminus \mathbb{D}_{l_0}$.

Therefore, the function
\[
\phi := \sum_{k \in \mathbb{D}_{l_0}} a_k \phi_{d_k} \ne 0
\]
is an eigenfunction of $-\Delta_g + V$ associated with the eigenvalue $\lambda$, and we have
\[
\Phi|_U = \mathcal{T}(a) = \phi|_U,
\]
which proves the claim \eqref{claim_lem_ortho_nonempty}.
\end{proof}

\begin{remark}\label{remark_nontrivial_anosov}
Let us emphasize that, in view of Lemma~\ref{lem_ortho_nonempty}, the mapping $\mathcal{T} : \ell^2(\mathbb{N}) \to L^2(U)$ introduced above is \textbf{not} surjective. Indeed, since the set $\{\lambda_{d_k}\}_{k=1}^\infty$ is a proper subset of $\{\lambda_k\}_{k=1}^\infty$ (see Definition~\ref{def_adm}), there exists an eigenvalue $\lambda_{k_0} \notin \{\lambda_{d_k}\}_{k=1}^\infty$. By applying Lemma~\ref{lem_ortho_nonempty} with $\tilde{V} = V$, it follows that for any eigenfunction $\Phi$ of $-\Delta_g + V$ associated with the eigenvalue $\lambda_{k_0}$, one has
\[
\Phi|_U \notin \operatorname{Im}\, \mathcal{T}.
\]
In particular, by the orthogonal decomposition theorem, we have
\begin{equation}\label{T_ortho}
L^2(U) = \operatorname{Im}\, \mathcal{T} \oplus \left( \operatorname{Im}\, \mathcal{T} \right)^\perp,
\end{equation}
where $\left( \operatorname{Im}\, \mathcal{T} \right)^\perp$ is a \textbf{nontrivial} closed subspace of $L^2(U)$.
\end{remark}

We are ready to prove Theorem~\ref{thm_anosov}. 

\begin{proof}[Proof of Theorem~\ref{thm_anosov}]
Let us begin with a few remarks about the setup. Akin to Remark~\ref{rmk_V_1_normal}, we may assume without loss of generality that $\{\psi_k^{(1)}\}_{k=1}^\infty$ forms an $L^2(M)$-orthonormal basis. 

Let $p \in O$ be as in the statement of Theorem~\ref{thm_anosov}, and fix a point $q \in \mathcal{A}_{M,g}(p) \subset O$. Let $\varepsilon \in (0,1)$ be sufficiently small so that condition~\eqref{balls_cond} and property~\textbf{(P)} in~\eqref{eq_condition_P} are satisfied, and so that the boundaries of the open connected sets $\mathcal{V}, \mathcal{W} \subset O$ defined in~\eqref{eq_W_V} are smooth. As a final preliminary remark, we note that, analogously to the proof of Lemma~\ref{lem_V_eq}, we have
\begin{equation}\label{V_eq_O}
	V_1 = V_2 \quad \text{on } O.
\end{equation}

In view of Remark~\ref{rmk_proof_anosov}, we will assume, as we may, that $D := \mathbb{N} \setminus \mathbb{A}$ is countably infinite. For the remainder of this proof, we define the map $\mathcal{T}_1 : \ell^2(\mathbb{N}) \to L^2(\mathcal{V})$ analogously to Proposition~\ref{prop_injectivity_anosov}, taking $V = V_1$, $U = \mathcal{V}$, and $\{\phi_k\}_{k=1}^\infty = \{\psi_k^{(1)}\}_{k=1}^\infty$. We also write $D = \{d_k\}_{k=1}^\infty$ with $d_1 < d_2 < \cdots$. That is,
\[
\mathcal{T}_1(a) = \left( \sum_{k=1}^\infty a_k \psi_{d_k}^{(1)} \right) \bigg|_{\mathcal{V}}, \quad a \in \ell^2(\mathbb{N}),
\]
where the series converges in $L^2(M)$.

We regard functions in $L^2(\mathcal{V})$ as elements of $L^2(M)$ via the standard extension by zero. By Remark~\ref{remark_nontrivial_anosov}, we have the following orthogonal decomposition:
\begin{equation}
\label{eq_orthog_dec_T_1}
L^2(\mathcal{V}) = \operatorname{Im} \mathcal{T}_1 \oplus (\operatorname{Im} \mathcal{T}_1)^\perp,
\end{equation}
where the subspace $(\operatorname{Im} \mathcal{T}_1)^\perp$ of $L^2(\mathcal{V})$ (and hence of $L^2(M)$) is nontrivial.

Given these considerations, let us fix an arbitrary
\[
f \in (\operatorname{Im} \mathcal{T}_1)^\perp.
\]
Since $\{\psi_k^{(1)}\}_{k=1}^\infty$ forms an $L^2(M)$-orthonormal basis, recalling that $D = \mathbb{N} \setminus \mathbb{A}$, and using that $\psi_{d_k}^{(1)}|_{\mathcal{V}} \in \operatorname{Im} \mathcal{T}_1$ for all $k \in \mathbb{N}$, we may write
\begin{equation}\label{f_1_exp_anosov}
	f = \sum_{k=1}^\infty c_k \psi_{a_k}^{(1)} \quad \text{on } M,
\end{equation}
where the coefficients $c_k$ are given by
\begin{equation}\label{c_k_def_anosov}
	c_k = (f, \psi_{a_k}^{(1)})_{L^2(\mathcal{V})} = \int_{\mathcal{V}} f(x)\, \overline{\psi_{a_k}^{(1)}(x)}\, dV_g = \int_{\mathcal{V}} f(x)\, \overline{\psi_{b_k}^{(2)}(x)}\, dV_g.
\end{equation}
The convergence in \eqref{f_1_exp_anosov} is to be understood in the $L^2(M)$-sense.

We claim that
\begin{equation}
\label{eigen_expansion_anosov}
	f = \sum_{k=1}^{\infty} \left(f, \psi^{(2)}_{b_k} \right)_{L^2(\mathcal{V})}\, \psi_{b_k}^{(2)} \qquad \text{in } H^{-n-1}(M),
\end{equation}
where the convergence on the right-hand side is to be understood in the $H^{-n-1}(M)$ sense. Here, we recall that $n \ge 2$ denotes the dimension of the manifold $M$. Although the proof of \eqref{eigen_expansion_anosov} is similar to that of \eqref{eigen_expansion} in Theorem~\ref{thm_spectral}, we include it here for the reader’s convenience, as we are working in the setting of distributions rather than smooth functions.

First, let us show that the right-hand side of \eqref{eigen_expansion_anosov} is well-defined in the $H^{-n-1}(M)$ topology. To that end, we first note that
\begin{equation}
\label{c_k_bound}
	|c_k| \leq \|f\|_{L^{2}(\mathcal V)} \, \|\psi_{a_k}^{(1)}\|_{L^{2}(M)} = \|f\|_{L^2(\mathcal V)},
\end{equation}
for all $k \in \mathbb{N}$.  Using the observability estimate~\eqref{eq_observability_eigenfunctions_def} for eigenfunctions of $-\Delta_g + V_2$ on $O$, as ensured by assumption~\textbf{(UO)}, together with~\eqref{psi_eq_1}, we obtain
\begin{equation}
\label{eq_18_1}
	\|\psi^{(2)}_{b_k}\|_{L^2(M)} \leq C\|\psi^{(2)}_{b_k}\|_{L^2(O)} = C\|\psi^{(1)}_{a_k}\|_{L^2(O)} 
	\leq C\|\psi^{(1)}_{a_k}\|_{L^2(M)} = C,
\end{equation}
where $C > 0$ is a constant independent of $k$; see also~\eqref{eq_500_1} and~\eqref{eq_500_2}.

Let $\tau > 0$ be such that, for both operators $P_j := -\Delta_g + V_j + \tau$, $j = 1, 2$, we have $P_j \ge I$ in the sense of self-adjoint operators. Fix $m \in \mathbb{N}$. For any $\varphi \in C^\infty(M)$, using \eqref{eq_18_1} and \eqref{eq_12_4} with $P = P_2$, we obtain
\[
|(\psi^{(2)}_{b_k}, \varphi)_{L^2(M)}| \le (\mu_{b_k}^{(2)}+\tau)^{-m} |(\psi^{(2)}_{b_k}, P_2^m \varphi)_{L^2(M)}| 
\le C (\mu_{b_k}^{(2)}+\tau)^{-m} \|\varphi\|_{H^{2m}(M)},
\]
and therefore,
\begin{equation}
\label{eq_18_2}
	\|\psi^{(2)}_{b_k}\|_{H^{-2m}(M)} \le C (\mu_{b_k}^{(2)}+\tau)^{-m},
\end{equation}
where $C > 0$ is independent of $k$.  Interpolating between the bounds \eqref{eq_18_2} and \eqref{eq_18_1}, we deduce that for any $s \ge 0$,
\begin{equation}
\label{eq_18_3}
	\|\psi^{(2)}_{b_k}\|_{H^{-s}(M)} \le C (\mu_{b_k}^{(2)}+\tau)^{-s/2},
\end{equation}
with $C > 0$ independent of $k$.

Now let $s > n$. From \eqref{c_k_bound}, \eqref{eq_18_3}, and the analogue of \eqref{eq_11_5} for the operator $P_2$, we obtain
\begin{equation}
\label{eq_18_4}
	\|c_k \psi_{b_k}^{(2)}\|_{H^{-s}(M)} \le C \|f\|_{L^2(\mathcal V)} (\mu_{b_k}^{(2)} + \tau)^{-s/2} 
	\le C \|f\|_{L^2(\mathcal V)} (\mu_k^{(2)} + \tau)^{-s/2} \le C \|f\|_{L^2(\mathcal V)} k^{-s/n},
\end{equation}
for all sufficiently large $k$. Here, we also used that $b_k \ge k$, and therefore $\mu_{b_k}^{(2)} \ge \mu_k^{(2)}$. The estimate \eqref{eq_18_4} implies that the series in \eqref{eigen_expansion_anosov} converges in $H^{-s}(M)$ for all $s > n$, and in particular in $H^{-n-1}(M)$. This shows that the right-hand side of \eqref{eigen_expansion_anosov} is indeed well-defined in the $H^{-n-1}(M)$ topology.

For now, we define $\tilde f \in H^{-n-1}(M)$ by
\[
\tilde f := \sum_{k=1}^\infty \left(f, \psi_{b_k}^{(2)}\right)_{L^2(\mathcal{V})} \, \psi_{b_k}^{(2)},
\]
where the convergence is understood in the $H^{-n-1}(M)$ topology.

Let us also emphasize that the equality~\eqref{psi_eq_1} implies
\begin{equation}\label{tildefonO_anosov}
	\tilde f|_{O} = f|_{O} \quad \text{in } \mathcal{D}'(O).
\end{equation}

Let us next define
\begin{equation}\label{norm_1}
U_1 \in C^1(\mathbb{R}; H^{-1}(M)) \cap C^0(\mathbb{R}; L^2(M))
\end{equation}
via
\begin{equation}\label{U_1_def_anosov}
U_1(t,\cdot) := \sum_{k=1}^{\infty} c_k \cos\left( \sqrt{\mu_{a_k}^{(1)} + \tau}\, t \right) \psi_{a_k}^{(1)}, \qquad t \in \mathbb{R},
\end{equation}
where the convergence is understood in the topology of the space in~\eqref{norm_1}.

Analogously, we define
\begin{equation}\label{norm_2}
U_2 \in C^1(\mathbb{R}; H^{-n-2}(M)) \cap C^0(\mathbb{R}; H^{-n-1}(M))
\end{equation}
via
\begin{equation}\label{U_2_def_anosov}
U_2(t, \cdot) := \sum_{k=1}^{\infty} c_k \cos\left( \sqrt{\mu_{b_k}^{(2)} + \tau}\, t \right) \psi_{b_k}^{(2)}, \qquad t \in \mathbb{R},
\end{equation}
with convergence understood in the topology of the space in~\eqref{norm_2}. The convergence in both cases can be proved analogously to that of the series defining $f$ and $\tilde f$, and is therefore omitted.

It is straightforward to verify that for each $j=1,2$, the function $U_j$ defined above satisfies the following initial value problem in the sense of distributions on $\mathbb{R} \times M$:
\begin{equation}\label{wave_pf_anosov}
\begin{cases}
\partial_t^2 U_j - \Delta_g U_j + (V_j(x) + \tau) U_j = 0 & \text{in } \mathcal{D}'(\mathbb{R} \times M), \\
U_j(0,\cdot) = f_j, & \\
\partial_t U_j(0, \cdot) = 0, & 
\end{cases}
\end{equation}
where $f_1 := f$ and $f_2 := \tilde f$ and the initial conditions are understood in the topologies of the respective function spaces in~\eqref{norm_1} and~\eqref{norm_2}. Moreover, $U_j$ is the unique solution to this problem in the space specified in~\eqref{norm_1} or~\eqref{norm_2}, respectively.

Next, we observe that the following equality holds for $t \in \mathbb{R}$:
\begin{equation}
	\label{U_j_equality_anosov}
	U_1(t,\cdot)|_{O} = U_2(t,\cdot)|_{O} \quad \text{in } \mathcal{D}'(O),
\end{equation}
thanks to the definitions \eqref{norm_1}--\eqref{U_2_def_anosov}, together with the identity \eqref{psi_eq_1}. We recall that $f \in L^2(\mathcal{V})$. Therefore, using the inequality
\[
\operatorname{dist}_{g}(\partial \mathcal{V}, \partial \mathcal{W}) \geq \operatorname{dist}_{g}(p,q) - \varepsilon,
\]
and the finite speed of propagation for the wave equation \eqref{wave_pf} with $j = 1$ (see \cite[Chapter 6, Proposition 1.3, p. 487]{Taylor_book_I} and \cite[Lemma 2.15]{Saksala_Shedlock_2025}), it follows that for $t \in (-T, T)$,
\begin{equation}
	\label{U_1_vanish_anosov}
	U_1(t,\cdot)|_{\mathcal{W}} = 0 \quad \text{in } \mathcal{D}'(\mathcal{W}),
\end{equation}
where
\[
T := \operatorname{dist}_{g}(\partial \mathcal{V}, \partial \mathcal{W}) \geq \operatorname{dist}_{g}(p,q) - \varepsilon.
\]
In view of \eqref{U_1_vanish_anosov} and \eqref{U_j_equality_anosov}, we deduce that for $t \in (-T, T)$,
\begin{equation}
	\label{U_2_vanish_anosov}
	U_2(t,\cdot)|_{\mathcal{W}} = 0 \quad \text{in } \mathcal{D}'(\mathcal{W}).
\end{equation}
Applying the global version of the sharp unique continuation result of Tataru for distributional solutions (see \cite{Tataru2} for the local version, \cite[Remark, p. 205]{Hor97}, and the arguments in the proof of the global version in \cite[Theorem 3.21, Theorem 3.24]{LL23}), we conclude that
\begin{equation}
\label{eq_18_5}
\tilde{f}|_{\widetilde{O}} = U_2(0,\cdot)|_{\widetilde{O}} = 0 \quad \text{in } \mathcal{D}'(\tilde{\mathcal{O}}),
\end{equation}
where
\[
\widetilde{O} = \{ x \in M : \operatorname{dist}_{g}(x,p) < \operatorname{dist}_{g}(p,q) - \varepsilon \}.
\]
Note that the property {\bf (P)} (see \eqref{eq_condition_P}) implies that $M \setminus O \subset \widetilde{O}$, and hence $M = O \cup \widetilde{O}$. Since both $O$ and $\widetilde{O}$ are nonempty open sets, we can use a partition of unity subordinate to the open cover $O \cup \widetilde{O}$ to conclude that \eqref{tildefonO_anosov}, together with the fact that $f = \tilde{f}$ in $\mathcal{D}'(\widetilde{\mathcal{O}})$ (which follows from \eqref{eq_18_5}), implies that $f = \tilde{f}$ in $\mathcal{D}'(M)$. This completes the proof of \eqref{eigen_expansion_anosov}.

Recalling Remark~\ref{rmk_proof_anosov}, we may assume that $\mathbb{N} \setminus \{b_k\}_{k=1}^\infty$ is nonempty. Let $\ell \in \mathbb{N} \setminus \{b_k\}_{k=1}^\infty$, and define
\[
\mathbb{B}_\ell := \{k \in \mathbb{N} \,:\, \mu_{b_k}^{(2)} = \mu_\ell^{(2)}\}.
\]
Note that the set $\mathbb{B}_\ell$ is either empty or finite.

\textbf{Case 1.} $\mathbb{B}_\ell = \emptyset$. Using $\psi^{(2)}_\ell \in C^\infty(M)$ as a test function, and noting that the series in \eqref{eigen_expansion_anosov} converges in $H^{-n-1}(M)$, we obtain
\begin{equation}\label{case_1}
\left(f, \psi^{(2)}_\ell \right)_{L^2(\mathcal{V})} = 0 \quad \text{for all } f \in (\operatorname{Im} \mathcal{T}_1)^\perp.
\end{equation}
The orthogonal decomposition \eqref{eq_orthog_dec_T_1} then implies that $\psi^{(2)}_\ell|_{\mathcal{V}} \in \operatorname{Im} \mathcal{T}_1$. Moreover, using \eqref{V_eq_O} and Lemma~\ref{lem_ortho_nonempty}, we deduce that there exists $k_{\ell} \in \mathbb{N} \setminus \mathbb{A}$ such that $\mu_\ell^{(2)} = \mu_{k_\ell}^{(1)}$, and that there exists an eigenfunction $\tilde \psi_{k_\ell}^{(1)}$ of $-\Delta_g + V_1$ on $M$, associated with the eigenvalue $\mu_\ell^{(2)}$, such that
\begin{equation}\label{case_1_1}
\psi^{(2)}_\ell|_{\mathcal{V}} = \tilde \psi^{(1)}_{k_\ell}|_{\mathcal{V}}.
\end{equation}
Given \eqref{V_eq_O} and the fact that $O$ is connected, the unique continuation principle for the operator $-\Delta_g + V_1 - \mu_\ell^{(2)}$ on the set $O$ implies that
\begin{equation}\label{case_1_2}
\psi^{(2)}_\ell|_{O} = \tilde \psi^{(1)}_{k_\ell}|_{O}.
\end{equation}

\textbf{Case 2.} $\mathbb{B}_\ell$ is finite. Arguing as in Case~1, we first obtain from \eqref{eigen_expansion_anosov} that
\begin{equation}\label{case_2}
\left(f, \psi^{(2)}_\ell - \sum_{k \in \mathbb{B}_\ell} (\psi^{(2)}_\ell, \psi_{b_k}^{(2)})_{L^2(M)}\, \psi_{b_k}^{(2)} \right)_{L^2(\mathcal{V})} = 0 \quad \text{for all } f \in (\operatorname{Im} \mathcal{T}_1)^\perp.
\end{equation}
Here, we define 
\begin{equation}
\label{eq_18_7}
\tilde{\psi}^{(2)}_\ell := \psi^{(2)}_\ell - \sum_{k \in \mathbb{B}_\ell} (\psi^{(2)}_\ell, \psi_{b_k}^{(2)})_{L^2(M)}\, \psi_{b_k}^{(2)},
\end{equation}
 which is an eigenfunction of $-\Delta_g + V_2$ on $M$, associated with the eigenvalue $\mu_\ell^{(2)}$, since the eigenfunctions $\psi^{(2)}_\ell$ and $\{\psi_{b_k}^{(2)}\}_{k \in \mathbb{B}_\ell}$ are linearly independent. Continuing as in Case~1, we conclude that there exists $k_{\ell} \in \mathbb{N} \setminus \mathbb{A}$ such that $\mu_\ell^{(2)} = \mu_{k_\ell}^{(1)}$, and that there exists an eigenfunction $\tilde \psi_{k_\ell}^{(1)}$ of $-\Delta_g + V_1$ on $M$, associated with the eigenvalue $\mu_\ell^{(2)}$, such that
\[
\tilde{\psi}^{(2)}_\ell|_{O} = \tilde \psi^{(1)}_{k_\ell}|_{O}.
\]

Hence, summarizing, we have shown that for any $\ell \in \mathbb{L} := \mathbb{N} \setminus \{b_k\}_{k=1}^\infty$, there exists $k_{\ell} \in \mathbb{N} \setminus \mathbb{A}$ such that
\begin{equation}
\label{eq_18_8}
\mu_\ell^{(2)} = \mu_{k_\ell}^{(1)}, \quad \tilde{\psi}^{(2)}_\ell|_{O} = \tilde \psi^{(1)}_{k_\ell}|_{O},
\end{equation}
where $\tilde{\psi}^{(2)}_\ell$ is an eigenfunction of $-\Delta_g + V_2$ on $M$, associated with the eigenvalue $\mu_\ell^{(2)}$, and such that either $\tilde{\psi}^{(2)}_\ell = \psi^{(2)}_\ell$, or $\tilde{\psi}^{(2)}_\ell$ is given by \eqref{eq_18_7}. Moreover, $\tilde \psi^{(1)}_{k_\ell}$ is an eigenfunction of $-\Delta_g + V_1$ on $M$, associated with the eigenvalue $\mu_\ell^{(2)}$.

Since the eigenfunctions $\{\psi_k^{(2)}\}_{k=1}^\infty$ are linearly independent, it follows that the combined collection $\{\tilde{\psi}^{(2)}_\ell\}_{\ell \in \mathbb{L}} \cup \{\psi_{b_k}^{(2)}\}_{k=1}^\infty$ is also linearly independent. Thus, by Lemma~\ref{lem_indep} and Remark~\ref{rem_indep_1}, their restrictions to $O$ remain linearly independent on $O$.  It therefore follows from \eqref{eq_18_8} and \eqref{psi_eq_1} that the eigenfunctions $\{\tilde{\psi}^{(1)}_{k_\ell}\}_{\ell \in \mathbb{L}}$ and $\{\psi_{a_k}^{(1)}\}_{k=1}^\infty$ are linearly independent.

We have thus reduced the proof of Theorem~\ref{thm_anosov} to the case where $\{b_k\}_{k=1}^\infty = \mathbb{N}$, by re-enumerating the new eigenfunctions of the operator $-\Delta_g + V_2$ according to its spectrum $\mu_1^{(2)} \leq \mu_2^{(2)} \leq \cdots$, and re-indexing the corresponding eigenvalues and eigenfunctions of $-\Delta_g + V_1$ accordingly. The proof is now complete, by Remark~\ref{rmk_proof_anosov}.
\end{proof}

\section{Proofs of Theorems~\ref{thm_main_heat}, \ref{thm_main_heat_metric}, \ref{thm_main_Schrodinger}, and~\ref{thm_main_wave}}

\label{proof_passive_measurements}

\subsection{Proof of Theorem \ref{thm_main_heat}}
Let $f_1, f_2 \in C^\infty(M)$ and $V_1, V_2 \in C^\infty(M)$. We begin by recalling the standard spectral representation of the solution $u_j^{f_j}$ to \eqref{int_heat_1} (with $V = V_j$ and $f = f_j$); see, for example, \cite[Chapter 6, pp.~482–483]{Taylor_book_I}.  To that end, for each $j = 1, 2$, let  
\[
\mu_1^{(j)} \le \mu_2^{(j)} \le \cdots \nearrow +\infty
\]
denote the eigenvalues of the operator $-\Delta_g + V_j$ on $M$, listed in increasing order and counted with multiplicity. Let $\{\psi_k^{(j)}\}_{k=1}^\infty \subset C^\infty(M)$ be the corresponding $L^2(M)$-orthonormal basis of eigenfunctions, satisfying
\[
(-\Delta_g + V_j)\,\psi_k^{(j)} = \mu_k^{(j)}\,\psi_k^{(j)} \quad \text{on } M.
\]
Let $\tau > 0$ be chosen such that, for both elliptic operators $P_j := -\Delta_g + V_j + \tau$, $j = 1, 2$, we have $P_j \ge I$ in the sense of self-adjoint operators.

The solution $u_j^{f_j}$ to \eqref{int_heat_1} then admits the spectral expansion
\begin{equation}
\label{eq_11_1}
u_j^{f_j}(t,x) = \sum_{k=1}^\infty e^{-\mu_k^{(j)} t}\, (f_j, \psi_k^{(j)})_{L^2(M)}\, \psi_k^{(j)}(x), \quad j = 1, 2.
\end{equation}
Since $f_j \in C^\infty(M)$, its Fourier coefficients satisfy the estimate
\begin{equation}
\label{eq_11_1_1}
\left|(f_j, \psi_k^{(j)})_{L^2(M)}\right| \le C_N\, (\mu_k^{(j)} + \tau)^{-N}, \quad \text{for all } N \in \mathbb{N}.
\end{equation}
Hence, the series \eqref{eq_11_1}, as well as the corresponding series of termwise $t$-derivatives, converge uniformly on $[0,\infty)$ in $H^s(M)$ for every $s \ge 0$. Consequently,
\[
u_j^{f_j} \in C^\ell\big([0,\infty); H^s(M)\big) \quad \text{for all } \ell, s \ge 0,
\]
and by Sobolev embedding, $u_j^{f_j} \in C^\infty\big([0,\infty) \times M\big)$.  

Similarly, one can establish the standard fact that the spectral representation \eqref{eq_11_1} defines a holomorphic function of $t$ in the right half-plane $\{ t \in \mathbb{C} : \operatorname{Re}(t) > 0 \}$, with values in $C^\infty(M)$. In particular, $u_j^{f_j}(t, x)$ is real-analytic in $t > 0$, with values in $C^\infty(M)$. 

Therefore, combining \eqref{int_heat_3} with the spectral expansion \eqref{eq_11_1}, we obtain
\begin{equation}
\label{eq_11_2}
\sum_{k=1}^\infty e^{-\mu_k^{(1)} t}\, (f_1, \psi_k^{(1)})_{L^2(M)}\, \psi_k^{(1)}(x)
= \sum_{k=1}^\infty e^{-\mu_k^{(2)} t}\, (f_2, \psi_k^{(2)})_{L^2(M)}\, \psi_k^{(2)}(x),
\end{equation}
for all $t>0$ and $x\in O$.

We now proceed to recover from \eqref{eq_11_2} the partial spectral data on $O$ required to apply Theorem~\ref{thm_spectral}, from which it follows that $V_1 = V_2$ on $M$. While this procedure is standard, we include the details for completeness and for the reader’s convenience; see also \cite{KKL, KKLM}.

We begin by taking the Laplace transform
\[
\mathcal{L}(h)(z) = \int_0^\infty h(t)\, e^{-z t} \, dt, \quad \text{for } \operatorname{Re}(z) > -\min\{\mu_1^{(1)}, \mu_1^{(2)}\},
\]
of both sides of \eqref{eq_11_2}. This yields
\begin{equation}
\label{eq_11_3}
\sum_{k=1}^\infty \frac{(f_1, \psi_k^{(1)})_{L^2(M)}\, \psi_k^{(1)}(x)}{\mu_k^{(1)} + z}
= \sum_{k=1}^\infty \frac{(f_2, \psi_k^{(2)})_{L^2(M)}\, \psi_k^{(2)}(x)}{\mu_k^{(2)} + z},
\end{equation}
valid for all $x \in O$ and all $z \in \mathbb{C}$ with $\operatorname{Re}(z) > -\min\{\mu_1^{(1)}, \mu_1^{(2)}\}$.

Let $\Omega_j := \mathbb{C} \setminus \bigcup_{k=1}^\infty \{-\mu_k^{(j)}\}$, and define
\begin{equation}
\label{eq_11_3_0}
\mathcal{R}_j(z;x) := \sum_{k=1}^\infty \frac{(f_j, \psi_k^{(j)})_{L^2(M)}\, \psi_k^{(j)}(x)}{\mu_k^{(j)} + z}, \quad z \in \Omega_j,\ x \in M,
\end{equation}
for $j = 1, 2$. We claim that the map $z \mapsto \mathcal{R}_j(z,\cdot)$ is holomorphic on $\Omega_j$, with possibly simple poles at the points $z=-\mu_k^{(j)}$, with values in $H^s(M)$ for all $s \ge 0$. In particular, by Sobolev embedding, $\mathcal{R}_j(z,\cdot)$ takes values in $C^0(M)$. The proof of the claim is similar to the argument used in \cite[Proposition 4.1]{FKU24}. For completeness, and since similar arguments will be used in the proofs of Theorems \ref{thm_main_heat_metric}, \ref{thm_main_Schrodinger}, and \ref{thm_main_wave} below, we present the proof for $\mathcal{R}_1$, as the case of $\mathcal{R}_2$ is analogous.

Let $K \subset \Omega_1$ be compact. Since each term in the definition of $\mathcal{R}_1(z,\cdot)$ is holomorphic in $z\in \Omega_1$ with values in $H^s(M)$, it suffices to show that the series converges uniformly in $z \in K$ in the $H^s(M)$-norm. First, there exists a constant $c > 0$ such that
\begin{equation}
\label{eq_11_3_1}
|\mu_k^{(1)} + z| \ge c \quad \text{for all } z \in K,\ \text{and all } k \in \mathbb{N},
\end{equation}
see \cite[Proposition 4.1]{FKU24}.  Next, by \eqref{eq_12_4} with $P=P_1$, we have
\begin{equation}
\label{eq_11_6}
\|\psi_k^{(1)}\|_{H^s(M)} \le C (\mu_k^{(1)} + \tau)^{s/2}.
\end{equation}
Finally, using \eqref{eq_11_1_1}, \eqref{eq_11_3_1}, \eqref{eq_11_5}, and \eqref{eq_11_6}, and choosing $N > s/2 + n/2$, we obtain for all $z \in K$ and all sufficiently large $k$,
\begin{equation}
\label{eq_11_7}
\frac{\|(f_1, \psi_k^{(1)})_{L^2(M)}\, \psi_k^{(1)}\|_{H^s(M)}}{|\mu_k^{(1)} + z|}
\le \frac{C_N}{( \mu_k^{(1)} + \tau)^{N - s/2}}
\le \frac{C_N}{k^{\frac{2}{n}(N - s/2)}}.
\end{equation}
Hence, the series in \eqref{eq_11_3_0} converges uniformly on $K$ in the $H^s(M)$-norm. Thus, $z \mapsto \mathcal{R}_1(z, \cdot)$ is holomorphic on $\Omega_1$ with values in $H^s(M)$. By analytic continuation, the identity \eqref{eq_11_3} remains valid for all $z \in \mathbb{C} \setminus \bigcup_{k=1}^\infty \{-\mu_k^{(1)}, -\mu_k^{(2)}\}$ and all $x \in O$.

Recall that there exists $k_0 \in \mathbb{N}$ such that the eigenvalues $\mu_k^{(1)}$ of the operator $-\Delta_g + V_1$ are simple for all $k \geq k_0$. Then, from \eqref{eq_11_3}, we obtain that for all $k \geq k_0$,
\begin{equation}
\label{eq_11_8}
\begin{aligned}
(f_1, &\psi_k^{(1)})_{L^2(M)}\, \psi_k^{(1)}(x)
= \lim_{z \to -\mu_k^{(1)}} (z+\mu_k^{(1)}) \mathcal{R}_1(z;x) = \lim_{z \to -\mu_k^{(1)}} (z+\mu_k^{(1)}) \mathcal{R}_2(z;x) \\
&= 
\begin{cases}
0, & \text{if } \mu_k^{(1)} \neq \mu_l^{(2)} \text{ for all } l \in \mathbb{N}, \\
\sum\limits_{l \in \mathbb{N} :\, \mu_k^{(1)} = \mu_l^{(2)}} (f_2, \psi_l^{(2)})_{L^2(M)}\, \psi_l^{(2)}(x), & \text{if } \mu_k^{(1)} = \mu_l^{(2)} \text{ for some } l,
\end{cases}
\end{aligned}
\end{equation}
for all $x \in O$. Since $\psi_k^{(1)}|_O \neq 0$ and by \eqref{int_heat_2_1}, it follows from \eqref{eq_11_8} that for each $k \geq k_0$, there exists $l_k \in \mathbb{N}$ such that 
\[
\mu_k^{(1)} = \mu_{l_k}^{(2)}, \quad \text{and} \quad \tilde \psi_k^{(1)}(x) = \tilde \psi_{l_k}^{(2)}(x), \quad x \in O,
\]
where
\[
\tilde \psi_k^{(1)} := (f_1, \psi_k^{(1)})_{L^2(M)}\, \psi_k^{(1)}, \quad 
\tilde \psi_{l_k}^{(2)} := \sum\limits_{l \in \mathbb{N} :\, \mu_k^{(1)} = \mu_l^{(2)}} (f_2, \psi_l^{(2)})_{L^2(M)}\, \psi_l^{(2)}.
\]
These functions are eigenfunctions of the operators $-\Delta_g + V_1$ and $-\Delta_g + V_2$, respectively, corresponding to the eigenvalues $\mu_k^{(1)} = \mu_{l_k}^{(2)}$. 
The systems $\{\tilde \psi_k^{(1)}\}_{k \ge k_0}$ and $\{\tilde \psi_{l_k}^{(2)}\}_{k \ge k_0}$ are linearly independent, since the eigenfunctions in these sequences correspond to distinct eigenvalues, see also Lemma \ref{lem_indep}. Finally, we note that the index sequence $(l_k)_{k \ge k_0}$ is strictly increasing. In view of Remark~\ref{rem_extension_basis}, we may apply Theorem~\ref{thm_spectral} to conclude that $V_1 = V_2$ on $M$.

We now show that $f_1 = f_2$ on $M$. Since $V_1 = V_2$, it follows that the difference $u_1^{f_1} - u_2^{f_2} \in C^\infty([0, \infty) \times M)$ satisfies the heat equation
\[
(\partial_t - \Delta_g + V_1)(u_1^{f_1} - u_2^{f_2}) = 0 \quad \text{on } (0, \infty) \times M.
\]
Moreover, by \eqref{int_heat_3}, we have $u_1^{f_1} - u_2^{f_2} = 0$ on $(0, \infty) \times O$. Since $M$ is connected, the unique continuation principle for heat equations (see, e.g., \cite[Sections 1 and 4]{Lin_F-H_1990}) implies that $u_1^{f_1} - u_2^{f_2} = 0$ on $(0, \infty) \times M$. In particular, evaluating at $t = 0$ yields $f_1 = f_2$ on $M$.
This concludes the proof of Theorem~\ref{thm_main_heat}.

\subsection{Proof of Theorem \ref{thm_main_heat_metric}}
Arguing as in the proof of Theorem~\ref{thm_main_heat}, we obtain the identity
\begin{equation}
\label{eq_11_9}
\sum_{k=0}^\infty \frac{(f_1, \psi_k^{(1)})_{L^2(M_1)}\, \psi_k^{(1)}(x)}{\mu_k^{(1)} + z}
= \sum_{k=0}^\infty \frac{(f_2, \psi_k^{(2)})_{L^2(M_2)}\, \psi_k^{(2)}(x)}{\mu_k^{(2)} + z},
\end{equation}
valid for all $z \in \Omega := \mathbb{C} \setminus \bigcup_{k=0}^\infty \{-\mu_k^{(1)}, -\mu_k^{(2)}\}$ and all $x \in O$. Both sides of \eqref{eq_11_9} define holomorphic functions of $z \in \Omega$, with simple poles at $z = -\mu_k^{(1)}$ (on the left-hand side) and $z = -\mu_k^{(2)}$ (on the right-hand side), for $k = 0, 1, 2, \dots$, taking values in $H^s(M_j)$ for any $s \ge 0$, for $j = 1, 2$, respectively. Therefore, by Sobolev embedding, they also take values in $C^0(M_j)$. 

For any $k = 0, 1, 2, \dots$, multiplying both sides of \eqref{eq_11_9} by $(z + \mu_k^{(1)})$ and taking the limit $z \to -\mu_k^{(1)}$, we obtain
\begin{equation}
\label{eq_11_10}
(f_1, \psi_k^{(1)})_{L^2(M_1)}\, \psi_k^{(1)}(x) = 
\begin{cases}
0, & \text{if } \mu_k^{(1)} \neq \mu_l^{(2)} \text{ for all } l \in \mathbb{N}, \\
(f_2, \psi_l^{(2)})_{L^2(M_2)}\, \psi_l^{(2)}(x), & \text{if } \mu_k^{(1)} = \mu_l^{(2)} \text{ for some } l,
\end{cases}
\end{equation}
for all $x \in O$. Here we used that all eigenvalues of the Schr\"odinger operator $-\Delta_{g_j} + V_j$ are simple, for $j = 1, 2$. Since $\psi_k^{(1)}|_O \neq 0$, and using \eqref{int_heat_2_1_metric}, it follows from \eqref{eq_11_10} that for each $k = 0, 1, 2, \dots$, there exists $l \in \mathbb{N} \cup \{0\}$ such that
\begin{equation}
\label{eq_11_11}
\mu_k^{(1)} = \mu_l^{(2)}, \quad \text{and} \quad \tilde \psi_k^{(1)}(x) = \tilde \psi_l^{(2)}(x), \quad x \in O.
\end{equation}
Here
\[
\tilde \psi_k^{(1)} := (f_1, \psi_k^{(1)})_{L^2(M_1)}\, \psi_k^{(1)}, \quad 
\tilde \psi_l^{(2)} := (f_2, \psi_l^{(2)})_{L^2(M_2)}\, \psi_l^{(2)}
\]
are eigenfunctions of the operators $-\Delta_{g_1}$ on $M_1$ and $-\Delta_{g_2}$ on $M_2$, respectively, corresponding to the eigenvalue $\mu_k^{(1)} = \mu_l^{(2)}$. Repeating the argument by multiplying both sides of \eqref{eq_11_9} by $(z + \mu_l^{(2)})$ and letting $z \to -\mu_l^{(2)}$, and using \eqref{int_heat_2_1_metric_eigen}, we conclude that $l = k$, hence 
\[
\mu_k^{(1)} = \mu_k^{(2)}, \qquad \tilde \psi_k^{(1)}|_O = \tilde \psi_k^{(2)}|_O \quad \text{for all } k \ge 0.
\]
The systems $\{\tilde \psi_k^{(1)}\}_{k=0}^\infty$ and $\{\tilde \psi_k^{(2)}\}_{k=0}^\infty$ are linearly independent, since the eigenfunctions in these sequences correspond to distinct eigenvalues. An application of \cite[Theorem 1.11]{FKU24} gives that there exists a smooth diffeomorphism $\Phi : M_1 \to M_2$ such that $\Phi|_O = \mathrm{Id}$ and $g_1 = \Phi^* g_2$.

We will next show that $f_1 = f_2 \circ \Phi$ on $M_1$. To that end, we first note that since $\Phi$ is a Riemannian isometry, we have
\begin{equation}
\label{eq_11_12}
(-\Delta_{g_1})(u \circ \Phi) = (-\Delta_{g_2} u) \circ \Phi,
\end{equation}
for all $u \in C^\infty(M_2)$; see \cite[pages 99--100]{GPR_book}. Let $u_2 = u_2^{f_2} \in C^\infty([0,\infty) \times M_2)$ be the solution to \eqref{int_heat_1} with $g = g_2$, $V = 0$, and initial data $f = f_2$. We define the map
\[
\tilde \Phi : [0,\infty) \times M_1 \to [0,\infty) \times M_2, \quad \tilde \Phi(t,x) := (t, \Phi(x)).
\]
Then $\tilde u_2 := u_2 \circ \tilde \Phi \in C^\infty([0,\infty) \times M_1)$ satisfies the following initial value problem for the heat equation:
\begin{equation}
\label{eq_11_13}
\begin{cases}
\partial_t \tilde u_2 - \Delta_{g_1} \tilde u_2 = 0 & \quad \text{on } (0,\infty) \times M_1, \\
\tilde u_2|_{t=0} = f_2 \circ \Phi.
\end{cases}
\end{equation}
Indeed, for $(t,x) \in (0,\infty) \times M_1$, we have
\[
\partial_t \tilde u_2(t,x) = (\partial_t u_2)(t, \Phi(x)), \quad \text{and} \quad -\Delta_{g_1} \tilde u_2(t,x) = (-\Delta_{g_2} u_2)(t, \Phi(x)),
\]
where the second identity follows from \eqref{eq_11_12}. Thus, \eqref{eq_11_13} follows.

Hence, $u_1^{f_1} - \tilde u_2 \in C^\infty([0,\infty) \times M_1)$ satisfies the heat equation
\[
(\partial_t - \Delta_{g_1})(u_1^{f_1} - \tilde u_2) = 0 \quad \text{on } (0,\infty) \times M_1.
\]
Moreover, by \eqref{int_heat_3_metric} and the fact that $\Phi|_O = \mathrm{Id}$, we have $u_1^{f_1} = \tilde u_2$ on $(0, \infty) \times O$. Since $M_1$ is connected, the unique continuation principle for the heat equation (see, e.g., \cite[Sections~1 and 4]{Lin_F-H_1990}) implies that $u_1^{f_1} = \tilde u_2$ on $(0, \infty) \times M_1$. In particular, evaluating at $t = 0$ yields $f_1 = f_2 \circ \Phi$ on $M_1$. This concludes the proof of Theorem~\ref{thm_main_heat_metric}.

\subsection{Proof of Theorem \ref{thm_main_Schrodinger}}
We will follow the proof of Theorem~\ref{thm_main_heat} and adopt the notation from that proof. The solution $u_j^{f_j}$ of \eqref{schrodinger_ip_pf} (with $V = V_j$ and $f = f_j$) then admits the spectral expansion
\begin{equation}
\label{eq_14_1}
u_j^{f_j}(t,x) = \sum_{k=1}^\infty e^{-\textrm{i} \mu_k^{(j)} t}\, (f_j, \psi_k^{(j)})_{L^2(M)}\, \psi_k^{(j)}(x), \quad j = 1, 2.
\end{equation}
Since $f_j \in C^\infty(M)$, the series \eqref{eq_14_1}, as well as the corresponding series of termwise $t$-derivatives, converge uniformly on $[0,\infty)$ in $H^s(M)$ for every $s \ge 0$. Consequently,
\[
u_j^{f_j} \in C^\ell\big([0,\infty); H^s(M)\big) \quad \text{for all } \ell, s \ge 0,
\]
and by Sobolev embedding, $u_j^{f_j} \in C^\infty\big([0,\infty) \times M\big)$.

It follows from \eqref{int_schrodinder_3} and \eqref{eq_14_1} that
\begin{equation}
\label{eq_14_2}
\sum_{k=1}^\infty e^{-\textrm{i}\mu_k^{(1)} t}\, (f_1, \psi_k^{(1)})_{L^2(M)}\, \psi_k^{(1)}(x)
= \sum_{k=1}^\infty e^{-\textrm{i}\mu_k^{(2)} t}\, (f_2, \psi_k^{(2)})_{L^2(M)}\, \psi_k^{(2)}(x),
\quad t > 0,\; x \in O.
\end{equation}
Taking the Laplace transform of both sides of \eqref{eq_14_2} for $ \operatorname{Re}(z) > 0$, we obtain that 
\begin{equation}
\label{eq_14_3}
\sum_{k=1}^\infty \frac{(f_1, \psi_k^{(1)})_{L^2(M)}\, \psi_k^{(1)}(x)}{z + \textrm{i}\mu_k^{(1)}}
= \sum_{k=1}^\infty \frac{(f_2, \psi_k^{(2)})_{L^2(M)}\, \psi_k^{(2)}(x)}{z + \textrm{i}\mu_k^{(2)}},
\end{equation}
valid for all $x \in O$ and all $z \in \mathbb{C}$ with $\operatorname{Re}(z) > 0$. 

Similar to the proof of Theorem~\ref{thm_main_heat}, one can show that both sides of \eqref{eq_14_3} define holomorphic functions of $z \in \mathbb{C} \setminus \bigcup_{k=1}^\infty \{-\textrm{i}\mu_k^{(j)}\}$, with possible simple poles at $z = -\textrm{i}\mu_k^{(1)}$ on the left-hand side and at $z = -\textrm{i}\mu_k^{(2)}$ on the right-hand side, for $k = 1, 2, \dots$, taking values in $H^s(M)$ for any $s \ge 0$. Therefore, by Sobolev embedding, they also take values in $C^0(M)$. Hence, it follows by analytic continuation that the identity \eqref{eq_14_3} extends to all $z \in \mathbb{C} \setminus \bigcup_{k=1}^\infty \{-\textrm{i}\mu_k^{(1)}, -\textrm{i}\mu_k^{(2)}\}$ and all $x \in O$.

Proceeding as in the proof of Theorem~\ref{thm_main_heat}, but using \eqref{eq_14_3} in place of \eqref{eq_11_3}, and applying Theorem~\ref{thm_spectral}, we conclude that $V_1 = V_2$ on $M$.

Finally, to show that $f_1 = f_2$ on $M$, we note that since $V_1 = V_2$, the difference $u_1^{f_1} - u_2^{f_2} \in C^\infty([0, \infty) \times M)$ satisfies the initial value problem for the Schr\"odinger equation
\[
\begin{cases}
\textrm{i}\,\partial_t (u_1^{f_1} - u_2^{f_2}) = (-\Delta_g  + V_1)(u_1^{f_1} - u_2^{f_2})
& \text{on } (0,\infty)\times M, \\
(u_1^{f_1} - u_2^{f_2})|_{t=0} = f_1 - f_2 & \text{on } M.
\end{cases}
\]
Since $O \subset M$ satisfies the geometric control condition, there exists a time $T > 0$ such that solutions to the Schr\"odinger equation are observable from $O$ over the time interval $[0, T]$.  Applying the observability estimate \eqref{eq_10_2}, we obtain that there exists a constant $C > 0$, depending only on $(M, g)$, $O$, $T$, and $V_1$, such that
\[
\|f_1 - f_2\|_{L^2(M)}^2 \leq C \int_0^T \|u_1^{f_1}(t,\cdot) - u_2^{f_2}(t,\cdot)\|_{L^2(O)}^2 \, dt.
\]
Combining this with the identity \eqref{int_schrodinder_3}, we conclude that $f_1 = f_2$ on $M$. This completes the proof of Theorem~\ref{thm_main_Schrodinger}.

\subsection{Proof of Theorem \ref{thm_main_wave}} We will follow the proof of Theorem~\ref{thm_main_heat} and adopt the notation from that proof. 
The solution $u_j^{f_j, h_j}$ of \eqref{wave_ip_pf} (with $V = V_j$, $f = f_j$, and $h = h_j$) admits the spectral expansion
\begin{equation}
\label{eq_15_1}
u_j^{f_j,h_j}(t,x) = \sum_{k=1}^\infty  \left[(f_j, \psi_k^{(j)})_{L^2(M)} \, \partial_t s_k^{(j)}(t) + (h_j, \psi_k^{(j)})_{L^2(M)} \, s_k^{(j)}(t)\right] \psi_k^{(j)}(x),
\end{equation}
where 
\begin{equation}
\label{eq_16_3}
s_k^{(j)}(t) =
\begin{cases}
 \frac{\sin(\sqrt{\mu_k^{(j)}} \, t)}{\sqrt{\mu_k^{(j)}}}, & \mu_k^{(j)} > 0, \\[1ex]
t, & \mu_k^{(j)} = 0, \\[1ex]
 \frac{\sinh(\sqrt{|\mu_k^{(j)}|} \, t)}{\sqrt{|\mu_k^{(j)}|}}, & \mu_k^{(j)} < 0,
\end{cases}
\end{equation}
for $k\in \N$ and $j = 1, 2$.

Since $f_j, h_j \in C^\infty(M)$, the series \eqref{eq_15_1}, as well as the corresponding series of termwise $t$-derivatives, converge uniformly on $[0,\infty)$ in $H^s(M)$ for every $s \ge 0$. Consequently,
\[
u_j^{f_j, h_j} \in C^\ell\big([0,\infty); H^s(M)\big) \quad \text{for all } \ell, s \ge 0,
\]
and by Sobolev embedding, $u_j^{f_j,h_j} \in C^\infty\big([0,\infty) \times M\big)$.

Now \eqref{int_wave_3} and \eqref{eq_15_1} imply that
\begin{equation}
\label{eq_15_2}
\begin{aligned}
&\sum_{k=1}^\infty  \left[(f_1, \psi_k^{(1)})_{L^2(M)} \, \partial_t s_k^{(1)}(t) + (h_1, \psi_k^{(1)})_{L^2(M)} \, s_k^{(1)}(t)\right] \psi_k^{(1)}(x)\\
&= \sum_{k=1}^\infty  \left[(f_2, \psi_k^{(2)})_{L^2(M)} \, \partial_t s_k^{(2)}(t) + (h_2, \psi_k^{(2)})_{L^2(M)} \, s_k^{(2)}(t)\right] \psi_k^{(2)}(x),
\quad  t > 0,\, x \in O.
\end{aligned}
\end{equation}

Taking the Laplace transform of \eqref{eq_15_2} for $\operatorname{Re}(z) > \max\left\{ \sqrt{|\mu_1^{(1)}|}, \sqrt{|\mu_1^{(2)}|} \right\}$, and using the identities
\[
\mathcal{L}(s_k^{(j)})(z) =
\begin{cases}
\frac{1}{z^2+\mu_k^{(j)}}, & \mu_k^{(j)} > 0, \\
\frac{1}{z^2}, & \mu_k^{(j)} = 0, \\
\frac{1}{z^2 - |\mu_k^{(j)}|}, & \mu_k^{(j)} < 0,
\end{cases}
\quad 
\mathcal{L}(\partial_t s_k^{(j)})(z) =
\begin{cases}
\frac{z}{z^2+\mu_k^{(j)}}, & \mu_k^{(j)} > 0, \\
\frac{z}{z^2}, & \mu_k^{(j)} = 0, \\
\frac{z}{z^2 - |\mu_k^{(j)}|}, & \mu_k^{(j)} < 0,
\end{cases}
\]
we obtain
\begin{equation}
\label{eq_15_3}
\begin{aligned}
&\sum_{k=1}^\infty  \frac{\left[(f_1, \psi_k^{(1)})_{L^2(M)} \, z + (h_1, \psi_k^{(1)})_{L^2(M)} \right] \psi_k^{(1)}(x)}{z^2 + \mu_k^{(1)}}\\
&= \sum_{k=1}^\infty  \frac{\left[(f_2, \psi_k^{(2)})_{L^2(M)} \, z + (h_2, \psi_k^{(2)})_{L^2(M)} \right] \psi_k^{(2)}(x)}{z^2 + \mu_k^{(2)}},
\end{aligned}
\end{equation}
valid for all $x \in O$ and all $z \in \mathbb{C}$ with $\operatorname{Re}(z) > \max\left\{ \sqrt{|\mu_1^{(1)}|}, \sqrt{|\mu_1^{(2)}|} \right\}$.

Let $\Sigma_j := \left( \bigcup_{\mu_k^{(j)} \ge 0} \{ \pm i\sqrt{\mu_k^{(j)}} \} \cup \bigcup_{\mu_k^{(j)} < 0} \{ \pm \sqrt{|\mu_k^{(j)}|} \} \right)$, $j=1,2$. 
Arguing as in the proof of Theorem~\ref{thm_main_heat}, one can show that both sides of~\eqref{eq_15_3} define holomorphic functions of $z \in \mathbb{C} \setminus \Sigma_j$, with possible poles at $z \in \Sigma_1$ on the left-hand side and at $z \in \Sigma_2$ on the right-hand side, taking values in $H^s(M)$ for any $s \ge 0$. Therefore, by the Sobolev embedding theorem, they also take values in $C^0(M)$. It then follows by analytic continuation that the identity~\eqref{eq_15_3} extends to all $z \in \mathbb{C} \setminus (\Sigma_1 \cup \Sigma_2)$ and all $x \in O$. 

Let $\tilde k_0 \in \mathbb{N}$ be such that $\tilde k_0 \ge k_0$ and $\mu_k^{(1)} > 0$ for all $k \ge \tilde k_0$, so that the eigenvalues $\mu_k^{(1)}$ of the operator $-\Delta_g + V_1$ are simple and positive in this range. In view of the assumption~\eqref{int_wave_2_1}, for any $k \ge \tilde k_0$, either
\begin{equation}
\label{eq_15_4}
(f_1, \psi_k^{(1)})_{L^2(M)} \, (i\sqrt{\mu_k^{(1)}}) + (h_1, \psi_k^{(1)})_{L^2(M)} \ne 0, 
\end{equation}
or
\begin{equation}
\label{eq_15_5}
(f_1, \psi_k^{(1)})_{L^2(M)} \, (-i\sqrt{\mu_k^{(1)}}) + (h_1, \psi_k^{(1)})_{L^2(M)} \ne 0. 
\end{equation}

Fix $k \ge \tilde k_0$ and assume that~\eqref{eq_15_4} holds; the case~\eqref{eq_15_5} can be treated similarly. Multiplying both sides of~\eqref{eq_15_3} by $(z - i\sqrt{\mu_k^{(1)}})$ and letting $z \to i\sqrt{\mu_k^{(1)}}$, we obtain that there exists $l_k \in \mathbb{N}$ such that
\[
\mu_k^{(1)} = \mu_{l_k}^{(2)}, \quad \text{and} \quad \tilde \psi_k^{(1)}(x) = \tilde \psi_{l_k}^{(2)}(x), \quad x \in O,
\]
where
\begin{align*}
\tilde \psi_k^{(1)} &:= \left[(f_1, \psi_k^{(1)})_{L^2(M)} \, (i\sqrt{\mu_k^{(1)}}) + (h_1, \psi_k^{(1)})_{L^2(M)} \right] \psi_k^{(1)},\\ 
\tilde \psi_{l_k}^{(2)} &:= \sum_{l \in \mathbb{N} \,:\, \mu_k^{(1)} = \mu_l^{(2)}} \left[(f_2, \psi_l^{(2)})_{L^2(M)} \, (i\sqrt{\mu_k^{(1)}}) + (h_2, \psi_l^{(2)})_{L^2(M)} \right] \psi_l^{(2)}.
\end{align*}
These functions are eigenfunctions of the operators $-\Delta_g + V_1$ and $-\Delta_g + V_2$, respectively, corresponding to the eigenvalue $\mu_k^{(1)} = \mu_{l_k}^{(2)}$. Since $k \ge \tilde k_0$ is arbitrary, Theorem~\ref{thm_spectral} implies that $V_1 = V_2$ on $M$.

To show that $f_1=f_2$ and $h_1=h_2$ on $M$, we first observe that $u_1^{f_1,h_1}-u_2^{f_2,h_2}\in C^\infty([0,\infty)\times M)$ safisfies the following initial value problem, 
\begin{equation}
	\label{eq_15_6}
	\begin{aligned}
		\begin{cases}
			(\partial_t^2 - \Delta_g + V_1)(u_1^{f_1,h_1}-u_2^{f_2,h_2}) = 0
			& \text{on } (0,\infty) \times M, \\
			(u_1^{f_1,h_1}-u_2^{f_2,h_2})|_{t=0} = f_1-f_2 & \text{on } M, \\
			\partial_t (u_1^{f_1,h_1}-u_2^{f_2,h_2})|_{t=0} = h_1-h_2 & \text{on } M.
		\end{cases}
	\end{aligned}
\end{equation}

Since $O \subset M$ satisfies the geometric control condition, there exists a time $T > 0$ such that the solution of \eqref{eq_15_6} satisfies an observability estimate from the open set $O$ over the time interval $[0, T]$. That is, there exists a constant $C > 0$, depending only on $(M, g)$, $O$, $T$, and $V_1$, such that
\begin{equation}
	\label{eq_15_7}
\|f_1 - f_2\|_{H^1(M)}^2 + \|h_1 - h_2\|_{L^2(M)}^2 \le C \int_0^T \|u_1^{f_1, h_1}(t,\cdot) - u_2^{f_2, h_2}(t,\cdot)\|_{H^1(O)}^2 \, dt,
\end{equation}
see \cite{Rauch_Taylor_1974}, \cite[Theorem 1.5]{Laurent_Leautaud_2016}, and also \cite[Proposition 1.2]{Le_Rousseau_Lebeau_TT} for the use of other energy spaces.

Combining \eqref{eq_15_7} with the identity \eqref{int_wave_3}, we conclude that $f_1 = f_2$ and $h_1 = h_2$ on $M$. This completes the proof of Theorem~\ref{thm_main_wave}.

\begin{appendix}

\section{Auxiliary results for the proof of Theorem~\ref{thm_spectral}}
\label{appendix_first_thm}

The following result is needed for the proof of Theorem~\ref{thm_spectral} in the case $N = 1$. Its proof is similar to that of \cite[Theorem 1.11]{FKU24} and is included here both because similar arguments are used in the proofs of Theorem~\ref{thm_spectral} for $N > 1$ and Theorem~\ref{thm_anosov}, and for completeness and the reader's convenience.

\begin{lemma}
\label{lemma_app_eq_f}
Let the notation be as in Theorem~\ref{thm_spectral} with $N = 1$ and its proof. Let $f_1 \in C^{\infty}_0(\mathcal{V})$ be defined by \eqref{f_1_exp}, and let $f_2 \in C^\infty(M)$ be defined by \eqref{tilde_f}. Then $f_1(x) = f_2(x)$ for all $x \in M$.
\end{lemma}

\begin{proof}
We follow the proof of \cite[Theorem 1.11]{FKU24}. First, the identity \eqref{psi_eq} with $N = 1$, together with the uniform convergence of the series in \eqref{f_1_exp} and \eqref{tilde_f}, implies that
\begin{equation}\label{tildefonO}
	f_2(x) = f_1(x), \qquad x \in O.
\end{equation}

We now show that $f_2$ vanishes identically on $M \setminus O$, so that the equality \eqref{tildefonO} extends to all of $M$. To this end, consider the functions
\begin{equation}\label{U_j_def}
	\begin{aligned}
		U_1(t,x) &= \sum_{k=1}^{\infty} a_k \cos\left(\sqrt{\mu_k^{(1)} + \tau}\, t\right) \,\psi_k^{(1)}(x), \qquad t \in \mathbb{R}, \quad x \in M, \\
		U_2(t,x) &= \sum_{k=1}^{\infty} a_k \cos\left(\sqrt{\mu_{b_k}^{(2)} + \tau}\, t\right)\, \psi_{b_k}^{(2)}(x), \qquad t \in \mathbb{R}, \quad x \in M,
	\end{aligned}
\end{equation}
where the coefficients $\{a_k\}_{k=1}^\infty$ are defined in \eqref{a_k_def}. These definitions are well-posed in the sense that the series converge in the norm of $C^l(\mathbb{R}; H^s(M))$ for all $l \ge 0$ and $s \ge 0$. This follows analogously to the justification given for the well-posedness of \eqref{tilde_f}, and we omit the details here.

It is straightforward to verify that for each $j = 1, 2$, the function $U_j \in C^\infty(\mathbb{R} \times M)$ is the unique solution to the initial value problem:
\begin{equation}\label{wave_pf}
	\begin{cases}
		\partial_t^2 U_j - \Delta_{g} U_j + (\tau + V_j(x)) U_j = 0 
		& \text{on } \mathbb{R} \times M, \\
		U_j(0,x) = f_j(x) & \text{on } M, \\
		\partial_t U_j(0,x) = 0 & \text{on } M.
	\end{cases}
\end{equation}
For the existence and uniqueness of solutions to \eqref{wave_pf}, see \cite[Theorem 2.13]{Saksala_Shedlock_2025}; see also \cite[Chapter 6, pp.\ 485–486]{Taylor_book_I} and \cite[Section 2.3]{KKL}.

In particular, it follows from \eqref{U_j_def} and \eqref{psi_eq} that
\begin{equation}
	\label{U_j_equality}
	U_1(t,x) = U_2(t,x), \qquad t \in \mathbb{R}, \quad x \in O.
\end{equation}
Since $f_1 \in C^{\infty}_0(\mathcal{V})$, we apply the finite speed of propagation for the wave equation \eqref{wave_pf} with $j = 1$. Using
\[
\mathrm{dist}_{g}(\partial \mathcal{V}, \partial \mathcal{W}) \geq \mathrm{dist}_{g}(p,q) - \varepsilon,
\]
we obtain
\begin{equation}
	\label{U_1_vanish}
	U_1(t,x) = 0, \quad t \in [-T, T], \quad x \in \mathcal{W},
\end{equation}
where
\[
T := \mathrm{dist}_{g}(\partial \mathcal{V}, \partial \mathcal{W}) \geq \mathrm{dist}_{g}(p,q) - \varepsilon,
\]
see \cite[Corollary 2.15]{Saksala_Shedlock_2025}. In view of \eqref{U_1_vanish} and \eqref{U_j_equality}, we deduce that
\begin{equation}
	\label{U_2_vanish}
	U_2(t,x) = 0, \quad t \in [-T, T], \quad x \in \mathcal{W}.
\end{equation}
By the global version of Tataru’s unique continuation theorem (\cite[Theorem 3.24]{LL23}, \cite{KKL}; see also \cite{Tataru} for the local version), it follows that
\[
f_2(x) = U_2(0,x) = 0 \quad \text{for all } x \in M \text{ such that } \mathrm{dist}_{g}(x,p) \leq \mathrm{dist}_{g}(p,q) - \varepsilon.
\]
Combining this with property~\textbf{(P)} (see \eqref{eq_condition_P}), we conclude that
\begin{equation}
	\label{f_2_vanish}
	f_2(x) = 0, \qquad x \in M \setminus O.
\end{equation}

It now follows from \eqref{f_2_vanish} and \eqref{tildefonO} that $f_2(x) = f_1(x)$ for all $x \in M$, completing the proof.
\end{proof}

The following result is needed for the proof of Theorem~\ref{thm_spectral} in the case $N = 1$. Its proof is the same as that of \cite[Theorem 1.11]{FKU24} and is included here for completeness and the reader’s convenience.

\begin{lemma}
\label{lem_app_orthonormal_basis}
Let the notation be as in Theorem~\ref{thm_spectral} with $N = 1$ and its proof. The identity \eqref{eigen_expansion_1} implies that the eigenfunctions $\{\psi_{k}^{(2)}\}_{k=1}^{\infty}$ form an orthonormal basis for $L^2(M)$.  
\end{lemma}

\begin{proof}
We follow the end of the proof of \cite[Theorem 1.11]{FKU24}.  We begin by renaming the ordered set $\{\psi_k^{(2)}\}_{k=1}^{\infty}$ according to the multiplicity $d_k$ of the distinct eigenvalues, rewriting it as the ordered set
\begin{equation}\label{new_basis}
\{\theta_{k,1},\ldots,\theta_{k,d_k}\}_{k=1}^{\infty}.
\end{equation}
We emphasize that this is merely a relabeling of the elements of the sequence, without altering their order. Equation \eqref{eigen_expansion_1} can now be restated as follows: for any $f \in C^{\infty}_0(\mathcal{V})$, we have
\begin{equation}\label{new_eigen_expansion}
f(x) = \sum_{k=1}^{\infty} \sum_{\ell=1}^{d_k} \left(f, \theta_{k,\ell} \right)_{L^2(\mathcal{V})} \, \theta_{k,\ell}(x), \qquad x \in M.
\end{equation}
Our goal is to show that \eqref{new_eigen_expansion} implies that, for any fixed $k \in \mathbb{N}$, the set $\{\theta_{k,1},\ldots,\theta_{k,d_k}\}$ is orthonormal in $L^2(M)$. To this end, we observe that, in view of \eqref{new_eigen_expansion} and the fact that eigenfunctions of $-\Delta_g + V_2$ corresponding to distinct eigenvalues are orthogonal in $L^2(M)$, we have
\begin{equation}
	\label{final_eigen_expansion}
	\left(f, \theta_{k,m} \right)_{L^2(\mathcal{V})} = \sum_{\ell=1}^{d_k} \left(f, \theta_{k,\ell} \right)_{L^2(\mathcal{V})} \, \left(\theta_{k,\ell}, \theta_{k,m} \right)_{L^2(M)},
\end{equation}
for all $f \in C^{\infty}_0(\mathcal{V})$, $k \in \mathbb{N}$, and $m = 1, \ldots, d_k$. It follows from this identity, together with \cite[Lemma 5.2]{FKU24}, that
\[
\left( \theta_{k,m}, \theta_{k,\ell} \right)_{L^2(M)} = \delta_{m\ell}, \quad \text{for all } k \in \mathbb{N}, \quad m, \ell = 1, \ldots, d_k,
\]
where $\delta_{m\ell}$ denotes the Kronecker delta. This shows that the set \eqref{new_basis}, and thus the sequence $\{\psi_k^{(2)}\}_{k=1}^{\infty}$, forms an orthonormal basis for $L^2(M)$. 
\end{proof}

The following result, which is needed for the proof of Theorem~\ref{thm_spectral} in the case $N = 1$, is a direct consequence of \cite[Theorem 1.1]{Saksala_Shedlock_2025}; see also \cite{Lu_Jinpeng_2025}. It addresses the classical Gel'fand inverse spectral problem for Schrödinger operators $-\Delta_g + V_j$, $j = 1, 2$, on a closed manifold, given orthonormalized spectral data. We refer to \cite[Corollary 2]{HLOS_2018} for a similar result in the case of the Laplacian on a closed manifold; see also \cite{KrKaLa} and \cite[Section 3]{Fei21}. For the case of manifolds with boundary, we refer to \cite{KKL}. The proof is included here for completeness and the reader’s convenience.

\begin{proposition}
\label{prop_Gelfand_standard}
Let $(M, g)$ be a smooth, closed, and connected Riemannian manifold of dimension $n \geq 2$. Let $O \subset M$ be a nonempty open set. Let $V_j \in C^\infty(M)$ be real-valued potentials, and denote by $\mu_1^{(j)} \leq \mu_2^{(j)} \leq \dots$ the eigenvalues of the operator $-\Delta_g + V_j$ on $M$, listed in increasing order and counted with multiplicity, for $j = 1, 2$. Let $\{\psi_k^{(j)}\}_{k=1}^\infty \subset C^\infty(M)$ be the corresponding $L^2(M)$-orthonormal basis of eigenfunctions for $-\Delta_g + V_j$, for $j = 1, 2$. 
Assume that
\begin{equation}
\label{eq_16_1}
\mu_k^{(1)} = \mu_k^{(2)} \quad \text{and} \quad \psi_k^{(1)}(x) = \psi_k^{(2)}(x), \quad \text{for all } x \in O \text{ and all } k \in \mathbb{N}.
\end{equation}
Then $V_1 = V_2$ on $M$.
\end{proposition}

\begin{proof}
We follow the proof of \cite[Corollary 2]{HLOS_2018}. Let $F \in C^\infty_0((0,\infty) \times O)$, and for $j = 1, 2$, consider the initial value problem for the wave equation:
\begin{equation}
\label{eq_16_2}
\begin{aligned}
\begin{cases}
\partial_t^2 u_j(t,x) - \Delta_g u_j(t,x) + V_j(x) u_j(t,x) = F(t,x) & \text{on } (0,\infty) \times M, \\
u_j(0,x) = 0 & \text{on } M, \\
\partial_t u_j(0,x) = 0 & \text{on } M.
\end{cases}
\end{aligned}
\end{equation}

It is well known that the problem \eqref{eq_16_2} admits a unique solution $u_j = u_j^{F} \in C^\infty([0,\infty) \times M)$ for $j = 1, 2$; see \cite[Theorem 2.13]{Saksala_Shedlock_2025}. The solution $u_j^F$ admits the following spectral representation:
\begin{equation}
\label{eq_16_2_1}
u_j^F(t,x) = \sum_{k=1}^{\infty} \left( \int_0^t (F(\tau, \cdot), \psi_k^{(j)})_{L^2(O)} \, s_k^{(j)}(t - \tau) \, d\tau \right) \psi_k^{(j)}(x),
\end{equation}
where $s_k^{(j)}(t)$ is given by \eqref{eq_16_3}, for all $k \in \mathbb{N}$ and $j = 1, 2$.

Since $F \in C^\infty_0((0,\infty) \times O)$, the series \eqref{eq_16_2_1}, as well as the series of its termwise $t$-derivatives, converge uniformly on $[0,\infty)$ in $H^s(M)$ for every $s \geq 0$. Hence, for any $F \in C^\infty_0((0,\infty) \times O)$, the equality \eqref{eq_16_1} implies that
\[
u_1^F(t,x) = u_2^F(t,x), \quad (t,x) \in (0,\infty) \times O,
\]
showing the equality of the local source-to-solution maps for the problems \eqref{eq_16_2}. An application of \cite[Theorem 1.1]{Saksala_Shedlock_2025} then yields that $V_1 = V_2$ on $M$.
\end{proof}

In the proof of Theorem~\ref{thm_spectral}, we will need the following result, which is similar to \cite[Lemma 5.2]{FKU24}. We include the proof for completeness and the reader’s convenience.
\begin{lemma}
\label{lem_set_F}
Let $\{\psi_k^{(1)}\}_{k=1}^\infty$ be as defined and satisfying the assumptions of Theorem~\ref{thm_spectral}, and let $L \geq 1$. Let $\mathcal{V} \subset M$ be a nonempty open set. Then for any $c = (c_1, \dots, c_{L}) \in \mathbb{C}^L$, there exists $\theta \in C^\infty_0(\mathcal{V})$ such that
\[
( \theta, \psi_k^{(1)} )_{L^2(\mathcal{V})} = c_k \quad \text{for all } k = 1, \dots, L.
\]
\end{lemma}

\begin{proof}
Since the set $\{\psi_k^{(1)}\}_{k=1}^\infty$ is linearly independent, it follows from Lemma~\ref{lem_indep} that $\{\psi_k^{(1)}|_{\mathcal{V}}\}_{k=1}^{L}$ is also linearly independent in $L^2(\mathcal{V})$. To prove the claim, it suffices to show that the linear map
\[
T : C_0^\infty(\mathcal{V}) \to \mathbb{C}^{L}, \quad T(\theta) := \left( (\theta, \psi_1^{(1)})_{L^2(\mathcal{V})}, \ldots, (\theta, \psi_{L}^{(1)})_{L^2(\mathcal{V})} \right)
\]
is surjective. Suppose, for contradiction, that $T$ is not surjective. Then there exists a nonzero vector $c = (c_1, \ldots, c_L) \in \mathbb{C}^L$ such that $T(\theta) \cdot c = 0$ for all  $\theta \in C_0^\infty(\mathcal{V})$,
where $\cdot$ denotes the standard Hermitian inner product on $\mathbb{C}^L$, defined by
$ v \cdot w := \sum_{k=1}^{L} v_k \overline{w_k}$.
This implies
\[
\left( \theta, \sum_{k=1}^{L} c_k \psi_k^{(1)} \right)_{L^2(\mathcal{V})} = 0 \quad \text{for all } \theta \in C_0^\infty(\mathcal{V}),
\]
and hence
\[
\sum_{k=1}^{L} c_k \psi_k^{(1)} = 0 \quad \text{on}\quad \mathcal{V},
\]
which contradicts the linear independence of the set $\{\psi_k^{(1)}|_{\mathcal{V}}\}_{k=1}^{L}$.
\end{proof}

\end{appendix}

\section*{Acknowledgements}

We are very grateful to Svetlana Jitomirskaya, Matti Lassas, Plamen Stefanov, Rakesh, and Gunther Uhlmann for helpful discussions. The research of K.K. is partially supported by the National Science Foundation (DMS 2408793).

\end{document}